\documentclass[11pt,twoside,leqno]{aomamlt2e} 
  \pageno{1243}
\received{March 31, 2003}
 \allowdisplaybreaks

\newcommand{\nn}{\nonumber}
\newcommand{\C}{{\mathbb C}}       
\newcommand{\R}{{\mathbb R}}       
\newcommand{\Z}{{\mathbb Z}}       
\newcommand{\DD}{{\mathcal D}}
\newcommand{\HH}{{\mathcal H}}

\newcommand{\CC}{{\mathcal C}}
\newcommand{\CE}{{{\mathcal C}_{\varepsilon}}}

\newcommand{\diam}{{\rm diam}}
\newcommand{\dist}{{\rm dist}}

\newcommand{\fiproof}{{\hspace*{\fill} $\square$ \vspace{2pt}}}

\newcommand{\lra}{{\longrightarrow}}

\newcommand{\real}{{\rm Re}}
\newcommand{\rf}[1]{{{\rm (\ref{#1})}}}
\newcommand{\supp}{{\rm supp}}

\newcommand{\vphi}{{\varphi}}
\newcommand{\ve}{{\varepsilon}}

\newcommand{\wt}[1]{{\widetilde{#1}}}
\newcommand{\wh}[1]{{\widehat{#1}}}

\newcommand{\sss}{{\rm Stop}}

\newcommand{\ttt}{{\rm Top}}
\newcommand{\ttd}{{\rm Top_{dy}}}
\newcommand{\tree}{{\rm Tree}}
\newcommand{\treeg}{{\rm Tree^{Reg}}}
\newcommand{\reg}{{\rm Reg_{dy}}}
\newcommand{\hhh}{{\HH^1_{\Gamma_R}}}
\newcommand{\wttt}{{\wh{\ttt}}}
\newcommand{\roo}{{\rm Root}}
\newcommand{\bal}{{\rm Bal}}
\newcommand{\qsss}{{\rm Qstp}}
\newcommand{\QS}{{\wh{\wh{S} \hspace{1mm}} \hspace{-1mm}}}
\newcommand{\maxbad}{{\rm Bad}}

\newtheorem{theorem}{Theorem}[section]
\newtheorem{lemma}[theorem]{Lemma}
\newtheorem{mlemma}[theorem]{Main Lemma}

\newtheorem{propo}[theorem]{Proposition}

\theorembodyfont{\upshape}

\newtheorem{remark}[theorem]{{\it Remark}}

\numberwithin{equation}{section}

\begin{document}
\currannalsline{162}{2005} 

\title{Bilipschitz maps, analytic capacity,\\ and the Cauchy
integral}

 \acknowledgements{Partially supported by the program Ram\'{o}n y Cajal (Spain)  
and by grants
BFM2000-0361 and MTM2004-00519 (Spain), 2001-SGR-00431 (Generalitat
de Catalunya), and HPRN-2000-0116 (European Union).}
 \author{Xavier Tolsa}

 \institution{Instituci\'o Catalana de Recerca i Estudis Avan\c{c}ats  
(ICREA) and \\ Universitat Aut\`onoma de Bar\-ce\-lo\-na, Spain\\
\email{xtolsa@mat.uab.es}}



 \shorttitle{Bilipschitz maps, analytic capacity,  and the Cauchy
integral}

\centerline{\bf Abstract}
\vglue12pt

Let $\vphi:\C\rightarrow \C$ be a bilipschitz map. We prove that
if $E\subset\C$ is compact, and $\gamma(E)$, $\alpha(E)$ stand for
its analytic and continuous analytic capacity respectively, then
$C^{-1}\gamma(E)\leq \gamma(\vphi(E)) \leq C\gamma(E)$ and
$C^{-1}\alpha(E)\leq \alpha(\vphi(E)) \leq C\alpha(E)$, where $C$
depends only on the bilipschitz constant of $\vphi$. Further, we
show that if $\mu$ is a Radon measure on $\C$ and the Cauchy
transform is bounded on $L^2(\mu)$, then the Cauchy transform is
also bounded on $L^2(\vphi_\sharp\mu)$, where $\vphi_\sharp\mu$ is
the image measure of $\mu$ by $\vphi$. To obtain these results, we
estimate the curvature of $\vphi_\sharp\mu$ by means of a corona
type decomposition.

\section{Introduction}

A compact set $E\subset \C$ is said to be removable for bounded analytic
functions if for any open set $\Omega$ containing $E$, every
bounded function analytic on $\Omega\setminus E$ has an analytic
extension to $\Omega$. In order to study removability,  in the 1940's  
Ahlfors \cite{Ahlfors}
introduced the notion of analytic capacity.
The {\em analytic capacity} of a compact set $E\subset\C$ is
$$\gamma(E) = \sup|f'(\infty)|,$$
where the supremum is taken over all analytic functions
$f:\C\setminus E\lra \C$ with $|f|\leq1$ on $\C\setminus E$, and
$f'(\infty)=\lim_{z\to\infty} z(f(z)-f(\infty))$.

In \cite{Ahlfors}, Ahlfors proved that $E$ is removable for bounded  
analytic functions if and only if
$\gamma(E)=0$.

Painlev\'{e}'s problem consists of characterizing removable
singularities for bounded analytic functions in a metric/geometric
way. By Ahlfors' result this is equivalent to describing compact
sets with positive analytic capacity in metric/geometric terms.

Vitushkin in the 1950's and 1960's showed that analytic capacity
plays a central role in problems of uniform rational approximation
on compact sets of the complex plane. Further, he introduced the
continuous analytic capacity $\alpha$, defined as
$$\alpha(E) = \sup|f'(\infty)|,$$
where the supremum is taken over all {\em continuous} functions
$f:\C\lra \C$ which are analytic on $\C\setminus E$, and uniformly
bounded by $1$ on $\C$. Many results obtained by Vitushkin in
connection with uniform rational approximation are stated in terms
of $\alpha$ and $\gamma$. See \cite{Vitushkin}, for example.

Until quite recently it was not known if removability is preserved
by an affine map such as $\vphi(x,y)=(x,2y)$ (with $x,y\in\R$).
 From the results of \cite{Tolsa-sem} (see Theorem A below) it
easily follows that this is true even for $\CC^{1+\ve}$
diffeomorphisms. In the present paper we show that this also
holds for bilipschitz maps. Remember that a map $\vphi:\C\lra\C$
is bilipschitz if it is bijective and there exists some constant
$L>0$ such that $$L^{-1}|z-w| \leq |\vphi(z)-\vphi(w)| \leq
L\,|z-w|$$ for all $z,w\in\C$. The precise result that we will prove
is the following.

\begin{theorem} \label{teogam}
Let $E\subset \C$ be a compact set and $\vphi:\C\rightarrow \C$ a
bilipschitz map. There exists a positive constants $C$ depending
only on $\vphi$ such that
\begin{equation} \label{compgam}
C^{-1}\gamma(E)\leq
\gamma(\vphi(E))\leq C\gamma(E)
\end{equation}
and
\begin{equation} \label{compalf}
C^{-1}\alpha(E)\leq \alpha(\vphi(E))\leq
C\alpha(E).
\end{equation}
\end{theorem}

As far as we know, the question on the behaviour of analytic
capacity under bilipschitz maps was first raised by Verdera
\cite[p.435]{Verdera-nato}. See also \cite[p.113]{Pajot} for a
more recent reference to the problem.

At first glance, the results stated in Theorem \ref{teogam} may
seem surprising, since $f$ and $f\circ\vphi$ are rarely both
analytic simultaneously. However, by the results of G. David
\cite{David}, it turns out that if $E$ is compact with finite
length (i.e. $\HH^1(E)<\infty$, where $\HH^1$ stands for the
$1$-dimensional Hausdorff measure), then $\gamma(E)>0$ if and only
if $\gamma(\vphi(E))>0$. Moreover, Garnett and Verdera \cite{GV}
proved recently that $\gamma(E)$ and $\gamma(\vphi(E))$ are
comparable for a large class of Cantor sets $E$ which may have non
$\sigma$-finite length.

Let us remark that the assumption that $\vphi$ is bilipschitz in  
Theorem \ref{teogam} is
necessary for \rf{compgam} or \rf{compalf} to hold. The precise  
statement reads as follows.

\begin{propo} \label{propconv}
Let $\vphi:\C\lra\C$ be a homeomorphism such that either \rf{compgam}  
holds for all compact sets
$E\subset \C${\rm ,} or \rf{compalf} holds for all compact sets
$E\subset \C$ \/{\rm (}\/in both cases with $C$ independent of $E${\rm )}. Then $\vphi$  
is bilipschitz.
\end{propo}


We introduce now some additional notation. A positive Radon
measure $\mu$ is said to have linear growth if there exists some
constant $C$ such that $\mu(B(x,r))\leq Cr$ for all $x\in\C$,
$r>0$. The linear density of $\mu$ at $x\in\C$ is (if it exists)
$$\Theta_\mu(x) = \lim_{r\to0}\frac{\mu(B(x,r))}r.$$

Given three pairwise different points $x,y,z\in\C$, their {\em
Menger curvature} is
$$c(x,y,z) = \frac{1}{R(x,y,z)},$$
where $R(x,y,z)$ is the radius of the circumference passing
through $x,y,z$ (with $R(x,y,z)=\infty$, $c(x,y,z)=0$ if $x,y,z$
lie on a same line). If two among these points coincide, we set
$c(x,y,z)=0$. For a positive Radon measure $\mu$, we define the
{\em curvature of $\mu$} as
\begin{equation} \label{defcurv}
c^2(\mu) = \int\!\!\int\!\!\int c(x,y,z)^2\, d\mu(x) d\mu(y)
d\mu(z).
\end{equation}
The notion of curvature of measures was introduced by Melnikov
\cite{Melnikov} when he was studying a discrete version of
analytic capacity, and it is one of the ideas which is responsible
for the recent advances in connection with analytic capacity.

Given a complex Radon measure $\nu$ on $\C$, the {\em Cauchy
transform} of $\nu$ is
$$\CC \nu(z) =\int \frac{1}{\xi-z}\, d\nu(\xi).$$
This definition does not make sense, in general, for
$z\in\supp(\nu)$, although one can easily see that the integral
above is convergent at a.e.\ $z\in\C$ (with respect to Lebesgue
measure). This is the reason why one considers the $\ve$-{\em
truncated Cauchy transform} of $\nu$, which is defined as
$$\CE \nu(z) =\int_{|\xi-z|>\ve} {\frac1{\xi-z}}\, d\nu(\xi),$$
for any $\ve>0$ and $z\in\C$. Given a $\mu$-measurable function
$f$ on $\C$ (where $\mu$ is some fixed positive Radon measure on
$\C$), we write $\CC_\mu f \equiv \CC (f \,d\mu)$ and
$\CC_{\mu,\ve} f \equiv \CE (f \,d\mu)$ for any $\ve>0$. It is
said that the Cauchy transform is bounded on $L^2(\mu)$ if the
operators $\CC_{\mu,\ve}$ are bounded on $L^2(\mu)$ uniformly on
$\ve>0$.

The relationship between the Cauchy transform and curvature of
measures was found by Melnikov and Verdera \cite{MV}. They proved
that if $\mu$ has linear growth, then
\begin{equation}\label{eqmv}
\|\CE\mu\|_{L^2(\mu)}^2 = \frac{1}{6} c^2_\ve(\mu) + O(\mu(\C)),
\end{equation}
where $c^2_\ve(\mu)$ is an $\ve$-truncated version of $c^2(\mu)$
(defined as in the right-hand side of \rf{defcurv}, but with the
triple integral over  $\{x,y,z\!\in\!\C:\! |x-y|,|y-z|,|x-z|\!>\!\ve\}$).
Moreover, there is also a strong connection (see \cite{Pajot})
between the notion of curvature of measures and the $\beta$'s from
Jones' travelling salesman theorem \cite{Jones}. The relationship
with Favard length is an open problem (see Section 6 of the
excellent survey paper \cite{Mattila-proj}, for example).


The proof of Theorem \ref{teogam}, as well as the one of the
result of Garnett and Verdera \cite{GV}, use the following
characterization of analytic capacity in terms of curvature of
measures obtained recently by the author.

\demo{\scshape Theorem A {\rm (\cite{Tolsa-sem})}} 
{\it For any compact $E\subset\C${\rm ,}
  $$\gamma(E)\simeq \sup\mu(E),$$
where the supremum is taken over all Borel measures $\mu$
supported on $E$ such that $\mu(B(x,r))\leq r$ for all $x\in E${\rm ,}
$r>0$ and $c^2(\mu)\leq \mu(E)$.}
\Enddemo

The notation $A\simeq B$ in the theorem means that there exists an
absolute constant $C>0$ such that $C^{-1}A\leq B\leq CA$.

The corresponding result for $\alpha$ is the following.

\demo{\scshape Theorem B {\rm (\cite{Tolsa-alfa})}} 
{\it For any compact $E\subset\C${\rm ,}
$$\alpha(E)\simeq \sup\mu(E),$$
where the supremum is taken over the Borel measures $\mu$
supported on $E$ such that $\Theta_\mu(x)=0$ for all $x\in E${\rm ,}
$\mu(B(x,r))\leq r$ for all $x\in E${\rm ,} $r>0${\rm ,} and $c^2(\mu)\leq
\mu(E)$.}
\Enddemo

Although the notion of curvature of a measure has a definite geometric  
flavor, it is not
clear if the characterizations of $\gamma$ and $\alpha$ in Theorems A  
and B can
be considered as purely metric/geometric. Nevertheless,
Theorem \ref{teogam} asserts that $\gamma$ and $\alpha$ have a
metric nature, in a sense.

Theorem \ref{teogam} is a direct consequence of the next result and of
Theorems A and B.

\begin{theorem} \label{teocurv}
Let $\mu$ be a Radon measure supported on a compact $E\subset \C${\rm ,}
such that $\mu(B(x,r))\leq r$ for all $x\in E,\,r>0$ and
$c^2(\mu)<\infty$. Let $\vphi:\C\rightarrow \C$ be a bilipschitz
mapping. There exists a positive constant $C$ depending only on
$\vphi$ such that
$$c^2(\vphi_{\sharp}\mu) \leq C\bigl(\mu(E) + c^2(\mu)\bigr).$$
\end{theorem}

In the inequality above, $\vphi_{\sharp}\mu$ stands for the image
measure of $\mu$ by $\vphi$. That is to say, $\vphi_{\sharp}\mu(A) =
\mu(\vphi^{-1}(A))$ for $A\subset \C$.

We will prove Theorem \ref{teocurv} using a corona type
decomposition, analogous to the one used by David and Semmes in
\cite{DS1} and \cite{DS2} for {\rm {\rm AD}} regular sets (i.e. for sets $E$
such that $\HH^1(E\cap B(x,r))\simeq r$ for all $x\in E$, $r>0$).
The ideas go back to Carleson's corona construction. See
\cite{AHMTT} for a recent survey on similar techniques. In our
situation, the measures $\mu$ that we will consider do not satisfy
any doubling or homogeneity condition. This fact is responsible
for most of the technical difficulties that appear in the proof of
Theorem \ref{teocurv}.

By the relationship \rf{eqmv} between curvature and the Cauchy
integral, the results in \cite{Tolsa-duke} (or in \cite{NTV}),
and Theorem \ref{teocurv}, we also deduce the next result.

\begin{theorem} \label{teocauchy}
Let $\vphi:\C\lra\C$ be a bilipschitz map and $\mu$ a Radon
measure on $\C$ without atoms. Set $\sigma = \vphi_\sharp\mu$. If
$\CC_\mu$ is bounded on $L^2(\mu)${\rm ,} then $\CC_\sigma$ is bounded
on $L^2(\sigma)$.
\end{theorem}

Notice that the theorem by Coifman, McIntosh and Meyer \cite{CMM}
concerning the $L^2$ boundedness of the Cauchy transform on
Lipschitz graphs (with respect to arc length measure) can be
considered as a particular case of\break Theorem \ref{teocauchy}.
Indeed, if $x\mapsto A(x)$ defines a Lipschitz graph on $\C$, then
the map $\vphi(x,y) = (x,y+A(x))$ is bilipschitz. Since $\vphi$
sends the real line to the\break Lipschitz graph defined by $A$ and the
Cauchy transform is bounded on $L^2(dx)$ on the real line (because
it coincides with the Hilbert transform), from Theorem
\ref{teocauchy} we infer that it is also bounded on the Lipschitz
graph.

The plan of the paper is the following. In Section \ref{secprelim} we  
prove (the easy)
Proposition \ref{propconv} and
introduce additional notation and definitions. The rest of the paper
is devoted to the proof of Theorem \ref{teocurv}, which we have
split into two main lemmas. The first one, Main Lemma
\ref{lemcorona}, deals with the construction of a suitable corona
type decomposition of $E$, and it is proved in Sections
\ref{seccorona}--\ref{seclds}. The second one, Main Lemma
\ref{lemafi}, is proved in Section \ref{seccurv}, and it shows how
one can estimate the curvature of a measure by means of a corona
type decomposition. So the proof of Theorem \ref{teocurv} works as
follows. In Main Lemma \ref{lemcorona} we construct a corona type
decomposition of $E$, which is stable under bilipschitz maps. That
is to say, $\vphi$ sends the corona decomposition of $E$ (perhaps
we should say of the pair $(E,\mu)$) to another corona
decomposition of $\vphi(E)$ (i.e.\ of the pair
$(\vphi(E),\vphi_{\sharp}\mu)$). Then, Main Lemma \ref{lemafi}
yields the required estimates for $c^2(\vphi_{\sharp}\mu)$.

\section{Preliminaries} \label{secprelim}

\Subsec{Proof of Proposition {\rm \ref{propconv}}}
Let $\vphi:\C\lra\C$ be a homeomorphism and suppose that  
$\gamma(\vphi(E))\simeq \gamma(E)$ for all
compact sets $E\subset \C$. Given $x,y\in\C$, consider the segment  
$E=[x,y]$. Then $\vphi(E)$ is a curve and its
analytic capacity is comparable to its diameter. Thus,
$$|\vphi(x)- \vphi(y)| \leq \diam(\vphi(E)) \simeq \gamma(\vphi(E))  
\simeq \gamma(E) \simeq |x-y|.$$
The converse inequality, $|x-y|\lesssim |\vphi(x)- \vphi(y)|$, follows  
by application of the previous argument
to $\vphi^{-1}$.

If instead of $\gamma(\vphi(E))\simeq \gamma(E)$ we assume now that  
with $\alpha(\vphi(E))\simeq \alpha(E)$ for all
compact sets $E$, a similar argument works. For example, given  
$x,y\in\C$, one can take $E$ to be the
closed ball centered at $x$ with radius $2|x-y|$, and then one can  
argue as above. \hfill\qed

\Subsec{Two remarks}
There are bijections $\vphi:\C\lra\C$ such that $\gamma(\vphi(E))\simeq  
\gamma(E)$ and
$\alpha(\vphi(E))\simeq \alpha(E)$, for any compact $E\subset \C$,  
which are not homeomorphisms. For example,
set $\vphi(z) = z$ if $\real(z) \geq 0$ and $\vphi(z) = z + i$ if  
$\real(z)<0$. Using the semiadditivity of
$\gamma$ and $\alpha$ one easily sees that $\gamma(\vphi(E))\simeq  
\gamma(E)$ and
$\alpha(\vphi(E))\simeq \alpha(E)$.

If the map $\vphi:\C\lra\C$ is assumed to be only Lipschitz, then none  
of the inequalities
$\gamma(\vphi(E))\gtrsim\gamma(E)$ or  
$\gamma(\vphi(E))\lesssim\gamma(E)$ holds, in general. To check
this, for the first inequality consider a constant map and $E$  
arbitrary with $\gamma(E)>0$. For the second
inequality, one only has to take into account that there are purely  
unrectifiable sets with finite length
which project orthogonally onto a segment (with positive length) in  
some direction.

\Subsec{Additional notation and definitions}
An Ahlfors-David regular curve (or {\rm {\rm AD}} regular curve) is a curve
$\Gamma$ such that $\HH^1(\Gamma\cap B(x,r))\leq C_3 r$ for all
$x\in \Gamma$, $r>0$, and some fixed $C_3>0$. If we want to
specify the constant $C_3$, we will say that $\Gamma$ is
``$C_3$-{\rm {\rm AD}} regular''.

In connection with the definition of $c^2(\mu)$, we also set
$$c^2_\mu(x) = \int\!\!\int c(x,y,z)^2\, d\mu(y) d\mu(z).$$
Thus, $c^2(\mu) = \int c^2_\mu(x)\, d\mu(x).$ If
$A\subset\C$ is {$\mu$}-measurable,
$$c^2_\mu(x,y,A)=\int_{A} c(x,y,z)^2 d\mu(z),\hspace{8mm} x,y\in\C,$$
and, if $A,B,C\subset\C$ are {$\mu$}-measurable,
$$c^2_\mu(x,A,B)=\int_{A}\int_{B} c(x,y,z)^2
d\mu(y) d\mu(z),\hspace{8mm}x\in\C,$$ and
$$c^2_\mu(A,B,C)=\int_{A}\int_{B}\int_{C} c(x,y,z)^2
d\mu(x) d\mu(y)d\mu(z).$$

The curvature operator $K_\mu$ is
$$K_\mu(f)(x)=\int k_\mu(x,y) f(y) d\mu(y),\hspace{8mm}
x\in\C,f\in L^1_{\rm loc}(\mu),$$ where $k_\mu(x,y)$ is the kernel
$$k_\mu(x,y)=\int c(x,y,z)^2 d\mu(z) = c^2_\mu(x,y,\C),\hspace{8mm}
x,y\in\C.$$
  For $j\in\Z$, the truncated operators
$K_{\mu,j}$, $j\in\Z$, are defined as
$$K_{\mu,j} f(x) = \int_{|x-y|>2^{-j}}
k_\mu(x,y)\,f(y)\,d\mu(y),\hspace{8mm} x\in\C,f\in
L^1_{\rm loc}(\mu).$$

In this paper, by a square we mean a square with sides parallel to
the axes. Moreover, we assume the squares to be half closed - half
open. The side length of a square $Q$ is denoted by $\ell(Q)$.
Given a square $Q$ and $a>0$, $aQ$ denotes the square concentric
with $Q$ with side length $a\ell(Q)$. The average (linear) density
of a Radon measure $\mu$ on $Q$ is
\begin{equation} \label{eqdens}
\theta_\mu(Q) := \frac{\mu(Q)}{\ell(Q)}.
\end{equation}

  A square $Q\subset \C$ is called $4$-dyadic if it is of the form
$[j2^{-n},\, (j+4)2^{-n}) \times [k2^{-n},\, (k+4)2^{-n})$, with
$j,k,n\in \Z$. So a $4$-dyadic square with side length
$4\cdot2^{-n}$ is made up of $16$ dyadic squares with side length
$2^{-n}$. We will work quite often with $4$-dyadic squares. All
our arguments would also work with other variants of this type of
square, such as squares $5Q$ with $Q$ dyadic, say. However, our
choice of $4$-dyadic squares has some advantages. For example, if
$Q$ is $4$-dyadic, $\frac12Q$ is another square made up of $4$
dyadic squares, and some calculations may be a little simpler.

Given a square $Q$ (which may be nondyadic) with side length
$2^{-n}$, we denote $J(Q):=n$. Given $a,b>1$, we say that $Q$ is
$(a,b)$-doubling if $\mu(aQ)\leq b\mu(Q)$. If we do not want to
specify the constant $b$, we say that $Q$ is $a$-doubling.

\begin{remark}
If $b>a^2$, then it easily follows that for $\mu$-a.e. $x\in \C$
there exists a sequence of $(a,b)$-doubling squares $\{Q_n\}_n$
centered at $x$ with $\ell(Q_n)\to0$ (and with
$\ell(Q_n)=2^{-k_n}$ for some $k_n\in\Z$ if necessary).
\end{remark}

As usual, in this paper the letter `$C$' stands for an absolute
constant which may change its value at different occurrences. On
the other hand, constants with subscripts, such as $C_1$, retain
their value at different occurrences. The notation $A\lesssim B$
means that there is a positive absolute constant $C$ such that
$A\leq CB$. So $A\simeq B$ is equivalent to $A\lesssim B \lesssim
A$.

\section{The corona decomposition}  \label{seccorona}

This section deals with the corona construction. In the next lemma we
will introduce a family $\ttt(E)$ of $4$-dyadic squares (the top
squares) satisfying some precise properties. Given any square
$Q\in\ttt(E)$, we denote by $\sss(Q)$ the subfamily of the squares
$P\in \ttt(E)$ satisfying
\begin{itemize}
\item[(a)] $P\cap3Q\neq\varnothing$,
\item[(b)] $\ell(P)\leq\frac1{8}\ell(Q)$,
\item[(c)] $P$ is maximal, in the sense that there does not exist
another square $P'\in \ttt(E)$ satisfying (a) and (b) which
contains $P$.
\end{itemize}
We also denote by  $Z(\mu)$ the set of points $x\in\C$ such that
there does not exist a sequence of $(70,5000)$-doubling squares
$\{Q_n\}_n$ centered at $x$ with $\ell(Q_n)\to0$ as $n\to\infty$,
so that moreover $\ell(Q_n)=2^{-k_n}$ for some $k_n\in\Z$. By the
preceding remark we have $\mu(Z(\mu))=0$.

The set of good points for $Q$ is defined as
$$G(Q):= 3Q\cap\supp(\mu)\setminus  
\Bigl[Z(\mu)\cup\bigcup_{P\in\sss(Q)}P\Bigr].$$

Given two squares $Q\subset R$, we set
$$\delta_\mu(Q,R) := \int_{R_Q\setminus Q} \frac1{|y-x_Q|}\,d\mu(y),$$
where $x_Q$ stands for the center of $Q$, and $R_Q$ is the
smallest square concentric with $Q$ that contains $R$.

\begin{mlemma}[The corona decomposition] \label{lemcorona}
Let $\mu$ be a Radon measure supported on $E\subset\C$ such that
$\mu(B(x,r))\leq C_0r$ for all $x\in\C, \,r>0$ and
$c^2(\mu)<\infty$.
  There exists a family ${\rm Top}(E)$ of $4$-dyadic
$(16,5000)$-doubling squares \/{\rm (}\/called {\em top squares)} which
satisfy the packing condition
\begin{equation} \label{pack}
\sum_{Q\in \ttt(E)} \theta_\mu(Q)^2 \mu(Q) \leq C \bigl(\mu(E) +
c^2(\mu)\bigr),
\end{equation}
  and such that for each square $Q\in \ttt(E)$ there exists a
  $C_3$-{\rm {\rm AD}} regular curve $\Gamma_Q$ such that\/{\rm :}\/
\begin{itemize}
\item[{\rm (a)}] $G(Q)\subset\Gamma_Q$.

\item[{\rm (b)}] For each $P\in \sss(Q)$ there exists some square
$\wt{P}$ containing $P$ such that $\delta_\mu(P,\wt{P})\leq
C\theta_\mu(Q)$ and $\wt{P}\cap \Gamma_Q\neq \varnothing$.

\item[{\rm (c)}] If $P$ is a square with $\ell(P)\leq \ell(Q)$
such that either $P\cap G(Q)\neq\varnothing$ or there is another
square $P'\in\sss(Q)$ such that $P\cap P'\neq\varnothing$ and
$\ell(P')\leq\ell(P)${\rm ,} then $\mu(P)\leq
C\,\theta_\mu(Q)\,\ell(P).$
\end{itemize}
Moreover{\rm ,} $\ttt(E)$ contains some $4$-dyadic square $R_0$ such
that $E\subset R_0$.
\end{mlemma}

Notice that the {\rm AD} regularity constant of the curves $\Gamma_Q$ in
the lemma is uniformly bounded above by the constant $C_3$.

  In Subsections \ref{secpreselec}, \ref{secstop} and
\ref{secelim} we explain how the $4$-dyadic squares in $\ttt(E)$
are chosen. Section \ref{secgam} deals with the construction of
the curves $\Gamma_Q$. The packing condition \rf{pack} is proved
in Sections \ref{secpack}--\ref{seclds}

The squares in $\ttt(E)$ are obtained by stopping-time arguments.
The first step consists of choosing a family $\ttt_0(E)$ which is
a kind of pre-selection of the $4$-dyadic squares which are
candidates to be in $\ttt(E)$. In the second step, some
unnecessary squares in $\ttt_0(E)$ are eliminated. The remaining
family  of squares is $\ttt(E)$.

 \Subsec{Pre-selection of the top squares}
\label{secpreselec}
To prove the Main Lemma \ref{lemcorona}, we will assume that $E$
is contained in a dyadic square with side length comparable to
$\diam(E)$. It easy to check that the lemma follows from this
particular case.

All the squares in $\ttt_0(E)$ will be chosen to be
$(16,5000)$-doubling. We define the family $\ttt_0(E)$ by
induction. Let $R_0$ be a $4$-dyadic square with $\ell(R_0)\simeq
\diam(E)$ such that $E$ is contained in one of the four dyadic
squares in $\frac12R_0$ with side length $\ell(R_0)/4$. Then, we
set $R_0\in\ttt_0(E)$. Suppose now that we have already decided
that some squares belong to $\ttt_0(E)$. If $Q$ is one of them,
then it generates a (finite or countable) family of ``bad''
$(16,5000)$-doubling $4$-dyadic squares, called
$\maxbad(Q)$. We will explain precisely below how this family is
constructed. For the moment, let us say that if $P\in\maxbad(Q)$,
then $P\subset 4Q$ and $\ell(P)\leq\ell(Q)/8$. One should think
that, in a sense, $\supp(\mu_{|3Q})$ can be well approximated by a
``nice'' curve $\Gamma_Q$ up to the scale of the squares in
$\maxbad(Q)$. All the squares in $\maxbad(Q)$ become also elements
of the family $\ttt_0(E)$.

In other words, we start the construction of $\ttt_0(E)$ by $R_0$.
The next squares that we choose as elements of $\ttt_0(E)$ are the  
squares from the family
$\maxbad(R_0)$. And the following ones are those generated as bad
squares of some square which is already in $\maxbad(R_0)$, and so on.  
The
family $\ttt_0(E)$ is at most countable. Moreover, in this process
of generation of squares of $\ttt_0(E)$, {\it a~priori\/}, it may  
happen that
some bad square $P$ is generated by two different squares
$Q_1,Q_2\in\ttt_0(E)$ (i.e. $P\in \maxbad(Q_1)\cap\maxbad(Q_2)$). We  
do not care about this fact.

\Subsec{The family $\maxbad(R)$} \label{secstop}
Let $R$ be some fixed $(16,5000)$-doubling $4$-dyadic square. We
will show now how we construct $\maxbad(R)$. Roughly speaking, a
square $Q$ with center in $3R$ and $\ell(Q)\leq\ell(R)/32$ is not good  
(we prefer to reserve
the terminology ``bad'' for the final choice)
for the approximation of $\mu_{|3R}$ by an Ahlfors regular curve
$\Gamma_R$ if either:
\begin{itemize}
\item[(a)] $\theta_\mu(Q)\gg \theta_\mu(R)$ (i.e.\ too high density), or
\item[(b)] $K_{\mu,J(Q)+10}\chi_E(x) -
K_{\mu,J(R)-4}\chi_E(x)$ is too big for ``many'' points $x\in Q$
(i.e. too high curvature), or
\item[(c)] $\theta_\mu(Q)\ll \theta_\mu(R)$ (i.e. too low
density).
\end{itemize}

A first attempt to construct $\maxbad(R)$ might consist of
choosing some kind of maximal family of squares satisfying
(a), (b) or (c). However, we want the squares from $\maxbad(R)$ to
be doubling, and so the arguments for the construction will become
somewhat more involved.

Let $A>0$ be some big constant (to be chosen below, in Subsection
\ref{sub53}), $\delta>0$ be some small constant (which will be fixed
in Section \ref{seclds}, depending on $A$, besides other
things), and $\ve_0>0$ be another small constant (to be chosen also in
Section \ref{seclds}, depending on $A$ and $\delta$). Let $Q$ be
some {\em $(70,5000)$-doubling} square centered at some point in
$3R\cap\supp(\mu)$, with $\ell(Q) =2^{-n} \ell(R)$, $n\geq5$. We
introduce the following notation:
\begin{itemize}
\item[(a)] If $\theta_\mu(Q) \geq A\theta_\mu(R)$, then we write $Q\in
{\rm HD}_{c,0}(R)$ (high density).

\item[(b)] If $Q\not\in {\rm HD}_{c,0}(R)$ and
$$\mu\bigl\{x\in Q:\,K_{\mu,J(Q)+10}\chi_E(x) - K_{\mu,
J(R)-4}\chi_E(x)\geq \ve_0\theta_\mu(R)^2\bigr\} \geq
\frac12\,\mu(Q),$$ then  $Q\in {\rm HC}_{c,0}(R)$ (high
curvature).

\item[(c)] If $Q\not\in {\rm HD}_{c,0}(R)\cup {\rm HC}_{c,0}(R)$ and
there exists some square $S_Q$ such that $Q\subset \frac1{100}
S_Q$, with $\ell(S_Q)\leq \ell(R)/8$ and $\theta_\mu(S_Q)\leq
\delta\, \theta_\mu(R)$, then we set $Q\in {\rm LD}_{c,0}(R)$ (low
density).
\end{itemize}
The subindex $c$ in ${\rm HD}_{c,0}$, ${\rm LD}_{c,0}$, and ${\rm HC}_{c,0}$ refers
to the fact that these families contain squares whose {\em
centers} belong to $\supp(\mu)$.

For each point $x\in 3R\cap\supp(\mu)$ which belongs to some
square from ${\rm HD}_{c,0}(R)\cup {\rm HC}_{c,0}(R)\cup {\rm LD}_{c,0}(R)$ consider the
largest square $Q_x\in {\rm HD}_{c,0}(R)\cup {\rm HC}_{c,0}(R)\cup
{\rm LD}_{c,0}(R)$ which contains $x$. Let $\wh{Q}_x$ be a $4$-dyadic
square with side length $4\ell(Q_x)$ such that
$Q_x\subset\frac12\wh{Q}_x$. Now we apply Besicovitch's covering
theorem to the family $\{\wh{Q}_x\}_x$ (notice that this theorem
can be applied because $x\in\frac12 \wh{Q}_x$), and we obtain a
family of $4$-dyadic squares $\{\wh{Q}_{x_i}\}_i$ with finite
overlap such that the union of the squares from ${\rm HD}_{c,0}(R)\cup  
{\rm HC}_{c,0}(R)\cup {\rm LD}_{c,0}(R)$
is contained (as a set in $\C$) in
$\bigcup_i \wh{Q}_{x_i}$. We define
$$\maxbad(R):=\{\wh{Q}_{x_i}\}_i.$$

If $Q_{x_i}\in {\rm HD}_{c,0}(R)$, then we write $\wh{Q}_{x_i}\in
{\rm HD}_0(R)$, and analogously with ${\rm HC}_{c,0}(R)$, ${\rm LD}_{c,0}(R)$ and
${\rm HC}_0(R)$, ${\rm LD}_0(R)$.

\begin{remark}
The constants that we denote by $C$ (with or without subindex) in
the rest of the proof of Main Lemma \ref{lemcorona} do not depend on
$A$, $\delta$, or $\ve_0$, unless stated otherwise.
\end{remark}

In the next two lemmas we show some properties fulfilled by the family
$\maxbad(R)$.

\begin{lemma}  \label{lemprop01}
Given $R\in\ttt_0(E)${\rm ,} the following properties hold for every
$Q\in\maxbad(R)$\/{\rm :}\/
\begin{itemize}
\item[{\rm (a)}] $Q$ is $(16,5000)$-doubling and $\frac12Q$ is
$(32,5000)$-doubling.
\item[{\rm (b)}] If $Q\in {\rm HD}_0(R)${\rm ,} then $\theta_\mu(Q) \gtrsim  
A\theta_\mu(R)$.
\item[{\rm (c)}] If $Q\in {\rm HC}_0(R)${\rm ,} then
$$\hskip-18pt \mu\bigl\{x\in \frac12 Q:\,K_{\mu,J(Q)+12}\chi_E(x) - K_{\mu,
J(R)-4}\chi_E(x)\geq C^{-1}\ve_0\theta_\mu(R)^2\bigr\} \gtrsim
\frac12 \mu(Q).$$
\item[{\rm (d)}] If $Q\in {\rm LD}_0(R)${\rm ,} then there exists some square $S_Q$ such  
that $Q\subset\frac1{20}S_Q${\rm ,} $\ell(S_Q)\leq \ell(R)/8${\rm ,} with
$\theta_\mu(S_Q)\lesssim \delta\,\theta_\mu(R).$
\end{itemize}
\end{lemma}

\Proof 
The doubling properties of $Q$ and $\frac12Q$ follow easily.
Let $x\in 3R\cap\supp(\mu)$ be such that $Q=\wh{Q}_x$, by the
notation above. Since $\frac12\wh{Q}_x \supset Q_x$, $Q_x$ is
$(70,5000)$-doubling, and $70Q_x\supset16\wh{Q}_x$, we get
$$\mu(\wh{Q}_x) \geq \mu(\tfrac12\wh{Q}_x) \geq \mu(Q_x) \geq
\frac1{5000}\mu(70Q_x) \geq \frac1{5000}\mu(16\wh{Q}_x).$$

The statements (b), (c) and (d) are a direct consequence of the  
definitions and of the fact that $\theta_\mu(Q_x)\simeq\theta_\mu(
\wh{Q}_x )$.
\hfill\qed

\begin{lemma}  \label{lemprop1}
Given $R\in\ttt_0(E)${\rm ,} the following properties hold for every
$Q\in\maxbad(R)$\/{\rm :}\/
\begin{itemize}
\item[{\rm (a)}] If $P$ is a square such that $P\cap Q\neq\varnothing$
and $\ell(Q)\leq\ell(P)\leq\ell(R)${\rm ,} then $$\mu(P)\leq C_4
A\,\theta_\mu(R)\,\ell(P).$$
\item[{\rm (b)}] If $P$ is a square concentric with $Q${\rm ,}
$\ell(P)\leq\ell(R)/8${\rm ,} and $\delta_\mu(Q,P)\geq C_5
A\theta_\mu(R)$ \/{\rm (}\/where $C_5>0$ is big enough\/{\rm ),}\/ then \end{itemize}\vglue-15pt
\begin{equation} \label{qaz1}
\mu(P)\geq \delta\,\theta_\mu(R)\,\ell(P)
\end{equation} \begin{itemize} \item[] and \end{itemize}\vglue-15pt
\begin{equation} \label{qaz2}
K_{\mu,J(P)+12}\chi_E(x) - K_{\mu,J(R)-2}\chi_E(x) \lesssim
A^2\theta_\mu(R)^2 \quad\mbox{for all $x\in P$.}
\end{equation}
\end{lemma}
\vglue8pt

Before proving the lemma we recall the following result, whose
proof follows by standard arguments (see \cite[Lemma
2.4]{Tolsa-duke}).

\begin{lemma}\label{curpert}
Let $x,y,z\in\C$ be three pairwise different points{\rm ,} and let
$x'\in\C$ be such that $C_6^{-1}|x-y|\leq |x'-y| \leq C_6|x-y|$.
Then{\rm ,}
$$|c(x,y,z)-c(x',y,z)|\leq (4+2C_6)\,\frac{|x-x'|}{|x-y|\,|x-z|}.$$
\end{lemma}

{\it Proof of Lemma} \ref{lemprop1}.
Let $x\in 3R\cap\supp(\mu)$ be such that $Q=\wh{Q}_x$, by  the
notation above.

Let us prove (a). If $\ell(R)/8<\ell(P)\leq \ell(R)$, then
$$\mu(P)\leq \mu(5R)\lesssim\mu(R) =
\theta_\mu(R)\ell(R)\lesssim\theta_\mu(R)\ell(P).$$ If $P$ is of
the form $2^nQ_x$, $n\geq 1$, with $\ell(P)\leq \ell(R)/8$, and
$P$ is $(70,5000)$-doubling, then
$$\mu(2^nQ_x)=\theta_\mu(2^nQ_x)\,\ell(2^nQ_x)\leq
A\theta_\mu(R)\ell(2^nQ_x),$$ by the definition of $Q_x$. If $P$
is of the same type but it is not doubling, then we take the
smallest $(70,5000)$-doubling square $\wt{P} :=2^mQ_x$ such that
$P\subset \wt{P}$ and $\ell(P)\leq \ell(R)/8$, in case that it
exists. If all the squares $2^mQ_x$ containing $P$ with
$\ell(2^mQ)\leq \ell(R)/8$ are non-$(70,5000)$-doubling, we set
$\wt{P}:=2^mQ_x$, with $m$ such that $\ell(P)=\ell(R)/8$. In any
case we have $$\theta_\mu(P) \leq \theta_\mu(2P) \leq
\theta_\mu(2^{2}P) \leq \cdots \leq \theta_\mu(\wt{P}) \leq
CA\theta_\mu(R).$$ The statement (a) for an arbitrary square $P$
such that $P\cap \wh{Q}_x\neq\varnothing$ and
$\ell(\wh{Q}_x)\leq\ell(P)\leq\ell(R)$ follows easily from the
preceding instances.

Now we turn our attention to (b). Let $P\supset Q$ be a square
concentric with $Q$ such that $\ell(P)\leq \ell(R)/8$ and
$\delta_\mu(Q,P)\geq C_5A\theta_\mu(R)$. It is easy to check (by
estimates analogous to the ones of \cite[Lemma 2.1]{Tolsa-bmo})
that if $C_5$ is chosen big enough, then there exists some
$(70,5000)$-doubling square $2^nQ_x$ such that $2Q\subset
2^nQ_x\subset \frac1{100}P$. Then $2^nQ_x\not\in {\rm LD}_{c,0}(R)$, and
by construction $\theta_\mu(P)\geq\delta\theta_\mu(R)$.

On the other hand, also by construction, there exists some $y\in
2^nQ_x$ such that
\begin{multline*}
K_{\mu,J(P)+12}\chi_E(y) - K_{\mu,J(R)-2}\chi_E(y)\\
\leq \; K_{\mu,J(2^nQ_x)+12}\chi_E(y) - K_{\mu,J(R)-4}\chi_E(y)
\;\leq\; \ve_0\theta_\mu(R)^2.
\end{multline*}
Then, for any $x\in P$,
\begin{eqnarray} \label{qaz33}
&&K_{\mu,J(P)+12} \chi_E(x) - K_{\mu,J(R)-2}\chi_E(x) \\   
&&\qquad =\iint_{2^{-12}\ell(P) < |x-t|\leq 4\ell(R)}
c(x,t,z)^2\,d\mu(t)d\mu(z)\nn \\
&&\qquad   \leq \,2 \iint_{2^{-12}\ell(P) < |x-t|\leq 4\ell(R)} c(y,t,z)^2  
\,d\mu(t)d\mu(z)\nn\\ &&\qquad\quad+
2\iint_{2^{-12}\ell(P) < |x-t|\leq 4\ell(R)} \Bigl[c(x,t,z) -  
c(y,t,z)\Bigr]^2\,d\mu(t)d\mu(z)\nn \\
&&\qquad=:2I_1 + 2I_2.
\nn
\end{eqnarray}
We have
$$I_1 =  K_{\mu,J(P)+12}\chi_E(y) - K_{\mu,J(R)-2}\chi_E(y) \leq  
\ve_0\theta_\mu(R)^2.$$
To estimate $I_2$ notice that by Lemma \ref{curpert},
$$\Bigl[c(x,t,z) - c(y,t,z)\Bigr]^2 \lesssim  
\max\biggl(\frac1{\ell(P)^2},\,\frac{\ell(2^nQ_x)^2}{|x-t|^2\,|x- 
z|^2}\biggr).$$
Integrating this inequality over $\{(t,z)\in\C^2:\,2^{-12}\ell(P) <  
|x-t|\leq 4\ell(R)\}$ (dividing $\C$ into annuli, for
example), one easily gets
$$I_2 \leq \biggl(\sup_{\lambda>1}\frac{\mu(\lambda P\cap  
16R)}{\ell(\lambda
P)}\biggr)^{2} \,\lesssim\; A^2\theta_\mu(R)^2.$$
Summing the estimates for $I_1$ and $I_2$, we see that \rf{qaz2}  
follows.
\hfill\qed

\Subsec{Elimination of unnecessary squares from $\ttt_0(E)$}
\label{secelim}
Some of the bad squares generated by each square $R\in
\ttt_0(E)$ may not be contained in $R$. This fact may cause
troubles when we try to prove a packing condition such as
\rf{pack} (because of the possible superposition of squares coming from  
different $R$'s in $\ttt_0(E)$).
The class $\ttt(E)$ that we are going to construct will
be a refined version of $\ttt_0(E)$, where some unnecessary squares  
will be eliminated.

Let us introduce some notation first. We say that a square
$Q\in\ttt_0(E)$ is a descendant of another square $R\in\ttt_0(E)$ if
there is a chain $R=Q_1,\,Q_2,\dots,Q_n=Q$, with
$Q_i\in\ttt_0(E)$ such that $Q_{i+1}\in\maxbad(Q_i)$ for each $i$.
Observe that, in principle, some square $Q$ may be a descendant of
more than one square $R$. Then, in the algorithm of
elimination below, $Q$ must be counted with multiplicity (so that $Q$ is
completely eliminated if it has been eliminated $m$ times, where
$m$ is the multiplicity of $Q$, etc.)
\smallbreak 
Let us describe the algorithm for constructing $\ttt(E)$. We have
to decide when any square in $\ttt_0(E)$ belongs to $\ttt(E)$. We
follow an induction procedure of elimination. A square in
$\ttt_0(E)$ that (during the algorithm we decide  belongs to
$\ttt(E)$) is said to be ``chosen for $\ttt(E)$''. If we decide
that it will not belong to $\ttt(E)$ (i.e. we eliminate it), we
say that it is ``unnecessary''. If we have not decided already if
it is either chosen for $\ttt(E)$ or unnecessary, we say that it
is ``available''.  We start with all the squares in $\ttt_0(E)$
being available, and at each step of the algorithm, some squares
are chosen for $\ttt(E)$ and others become unnecessary.

Let $R_0$ be the $4$-dyadic square containing $E$ defined at the
beginning of Subsection \ref{secpreselec}. We start by choosing
$R_0$ for $\ttt(E)$. Let $R_1$ be (one of) the squares from
$\maxbad(R_0)$ with maximal side length. We choose $R_1$ for
$\ttt(E)$ too. Next, we choose for $\ttt(E)$ (one of) the
available square(s) $R_2\in\ttt_0(E)$ with maximal side length. At
this point, some available squares in $\ttt_0(E)$ may become
unnecessary. First, these are the squares $Q$ with
$Q\in\maxbad(R)$ for some $R\in\ttt_0(E)$ such that $Q\subset R_2$
and $\ell(R_2)\leq\ell(R)/8$ (notice that this implies that either
$R=R_0$ or $R=R_1$). Also, all the squares which are descendants
of unnecessary squares become unnecessary.

Suppose now that we have chosen $R_0,\,R_1,\dots,R_{k-1}$ for
$\ttt(E)$, with $\ell(R_0)\geq\ell(R_1)\geq
\cdots\geq\ell(R_{k-1})$, and that the available squares in
$\ttt_0(E)$ have side length $\leq \ell(R_{k-1})$. Let $R_{k}$ be
(one of) the {\em available} square(s) in $\ttt_0(E)$ with maximal
side length. We choose $R_{k}$ for $\ttt(E)$. The squares that
become unnecessary are those available squares $Q$ such that
$Q\in\maxbad(R)$ for some $R\in\ttt_0(E)$ with $Q\subset R_{k}$ and
$\ell(R_{k})\leq\ell(R)/8$ (this implies that $R$ coincides with
one of the squares $R_1,\dots,R_{k-1}$). Again, all the squares
which are descendants of unnecessary squares become unnecessary
too.

It is easily seen that following  this algorithm one will decide
if any square in $\ttt_0(E)$ is unnecessary or chosen for
$\ttt(E)$. Notice that the squares $Q\in \maxbad(R)$ which are
eliminated after choosing $R_{k}$ are the ones such that $R_k$
becomes an ``intermediate'' square (in a sense) between $Q$ and
its generator $R$ (i.e.\ the square $R\in\ttt_0(E)$ such that
$Q\in\maxbad(R)$), as well as descendants of already eliminated
squares. Moreover, if a square $Q$ has been eliminated but its
generator $R$ has been chosen for $\ttt(E)$, it means that there
is another chosen square $Q'\in\ttt(E)$ which contains $Q$, with
$\ell(Q')\leq\ell(R)/8$. Thus, if $R\in\ttt(E)$, then
\begin{equation} \label{unio22}
\bigcup_{Q\in\maxbad(R)}Q \subset \bigcup Q',
\end{equation}
  where the union on the right side is over the squares $Q'\in\ttt(E)$
  such that there exists $Q\in\maxbad(R)$ contained in $Q'$ and
$\ell(Q')\leq\ell(R)/8$.

Remember the definition of $\sss(R)$, for $R\in\ttt(E)$, given at
the beginning of Section \ref{seccorona}.
Notice that
$\maxbad(R)\cap\ttt(E)\subset\sss(R)$. Of course, the opposite
inclusion is false in general. By \rf{unio22}, we also have
  $$\bigcup_{Q\in\maxbad(R)}Q \subset \bigcup_{Q'\in\sss(R)}Q'.$$

Given $R\in\ttt(E)$ and $Q\in\sss(R)$, we write $Q\in {\rm HD}(R)$ if
there exists some $R'\in\ttt(E)$ such that $Q\in\maxbad(R')\cap
{\rm HD}(R')$  analogously with ${\rm LD}(R)$ and ${\rm HC}(R)$.

\begin{remark} \label{statclau}
Let us insist again on the following fact: for any $R\in\ttt(E)$,
if $Q\in\maxbad(R)$, then either $Q\in\sss(R)$ or otherwise there is  
some $Q'\in\sss(R)$ such that
$Q'\supset Q$.
\end{remark}

\begin{remark}  \label{remww}
Changing constants if necessary, the properties (a) and (b) of Lemma  
\ref{lemprop1} also hold if
instead of assuming $Q\in \maxbad(R)$, we suppose that $Q\in\sss(R)$.
This is due to the fact that, with the new assumption $Q\in\sss(R)$,
the squares $P$ considered in Lemma \ref{lemprop1} (a), (b) will be,  
roughly speaking, a subset
of the corresponding squares $P$ with the assumption $Q\in \maxbad(R)$,  
because of the preceding remark.

On the other hand, in principle, (b), (c) and (d) in Lemma  
\ref{lemprop01} may fail.
Nevertheless, we will see in Lemma \ref{lemprop4}
below that they still do hold in some special cases.
\end{remark}

\section{Construction of the curves $\Gamma_R$, $R\in \ttt(E)$}  
\label{secgam}

 \vglue-12pt
\Subsec{P. Jones\/{\rm '}\/ travelling salesman theorem}
To construct the curves $\Gamma_R$ for $R\in\ttt(E)$, we will
apply P. Jones' techniques. Before stating the precise result that
we will use, we need to introduce some notation. Given a set
$K\subset \C$ and a square $Q$, let $V_Q$ be an infinite strip (or
line in the degenerate case) of smallest possible width which
contains $K\cap 3Q$, and let $w(V_Q)$ denote the width of $V_Q$.
Then we set
  $$\beta_K(Q) = \frac{w(V_Q)}{\ell(Q)}.$$ We will use the
following version of Jones' travelling salesman theorem
\cite{Jones}:

\demo{\scshape Theorem C {\rm (P.\ Jones)}} 
{\it A set $K\subset \C$ is contained in an {\rm AD} regular curve if and
only if there exists some constant $C_7$ such that for every
dyadic square $Q$
\begin{equation} \label{betajones}
\sum_{P\in \DD,P\subset Q} \beta_K(P)^2\ell(P) \leq C_7\ell(Q).
\end{equation}
The {\rm AD} regularity constant of the curve depends on $C_7$.}
 \Enddemo

Let us mention that in \cite{MMV}, using an $L^2$ version of
Jones' theorem due to David and Semmes \cite{DS2}, the authors showed
that if $\mu$ is a measure such that $\mu(B(x,r))\simeq r$
and $c^2(\mu_{|B(x,r)}) \leq C\mu(B(x,r))$ for all
$x\in\supp(\mu)$, $0<r\leq\diam(\supp(\mu))$, then $\supp(\mu)$ is
contained in an {\rm AD} regular curve. In our case, the measure
$\mu_{|R}$ does not satisfy these conditions. However, in a sense,
they do hold for ``big'' balls $B(x,r)$, at scales sufficiently
above the stopping squares.

In order to apply Jones' result, we will construct a set $K$ which
approximates $\supp(\mu)\cap 3R$ at some level above the stopping
squares and then, using Theorem C, we will show that there exists a
curve $\Gamma_R$ which contains $K$. We have not been able to
extend the arguments in \cite{MMV} to our situation. Instead, our
approach is based on another idea of Jones which shows how one
can estimate the $\beta$ numbers of an {\rm AD} regular set in terms of
the curvature of a measure (see \cite[Th.~38]{Pajot}).

 \Subsec{Balanced squares}
Before constructing the appropriate set $K$ which approximates
$\supp(\mu)$ we need to show the existence of some squares that
will be called {\em balanced squares}, which will be essential for
our calculations.

Given $a,b>0$, we say that a square $Q$ is balanced with
parameters $a,\,b$ (or $(a,b)$-balanced) with respect to $\mu$ if
there exist two squares $Q_1,\,Q_2\subset Q$ such that
\begin{itemize}
\item $\dist(Q_1,Q_2)\geq a\,\ell(Q)$,
\item $\ell(Q_i) \leq \dfrac{a}{10} \ell(Q)$ for $i=1,2$, and
\item $\mu(Q_i)\geq b\mu(Q)$ for $i=1,2$.
\end{itemize}
We write $Q\in \bal_\mu(a,b)$.

The next lemmas deal with the existence of balanced squares.

\begin{lemma} \label{lembalan}
Let $a=1/40$. If $Q\not\in\bal_\mu(a,b)${\rm ,} then there exists a
square $P\subset Q$ with $\ell(P)\leq\ell(Q)/10$ such that
$\mu(P)\geq (1-2\cdot10^5b)\mu(Q).$
\end{lemma}

\Proof 
We set $N:=400$. We split $Q$ into $N^2$ squares $Q_k$,
$k=1,\dots,N^2$, of side length $\ell(Q)/N=a\ell(Q)/10$. We
put$$G = \{Q_k:\,1\leq k \leq N^2,\;\mu(Q_k)\geq b\mu(Q)\}.$$
Since $Q\not\in\bal_\mu(a,b)$, given any pair of squares
$Q_j,\,Q_k$, $1\leq j,k \leq N^2$ such that $\dist(Q_j,Q_k)\geq
a\ell(Q)$, it turns out that one of the two squares, say $Q_j$,
satisfies $\mu(Q_j)\leq b\mu(Q)$. Therefore, all the squares from
$G$ are contained in a ball $B_0$ of radius $2a\ell(Q)$, since
$$\Bigl(\frac{2a}{10} + a\Bigr)2^{1/2}\ell(Q) \leq 2a\ell(Q).$$
Thus, $$\sum_{k:Q_k\in G} \mu(Q_k) \leq \mu(B_0\cap Q).$$ Also,
$$\sum_{k:Q_k\not\in G} \mu(Q_k) \leq b\sum_{k:Q_k\not\in
G}\mu(Q)\leq bN^2\mu(Q) = \frac{100b}{a^2}\,\mu(Q).$$ Then we have
$$\mu(Q) \leq \mu(B_0\cap Q) + \frac{100b}{a^2}\,\mu(Q).$$ Thus,
$$\mu(B_0\cap Q)\geq \Bigl(1- \frac{100b}{a^2}\Bigr)\mu(Q)\geq
(1- 2\cdot10^5b)\mu(Q).$$ Since the radius of $B_0$ equals
$2a\ell(Q) = \ell(Q)/20$, there exists some square $P\subset Q$
with side length $\ell(Q)/10$ which contains $B_0\cap Q$, and we
are done.
\hfill\qed

\begin{lemma}
Let $Q$ be a square such that $2Q\not\in \bal_\mu(1/40,\,b)${\rm ,} and
suppose that $\theta_\mu(2Q) \leq C_4A\theta_\mu(R)$ and
$\theta_\mu(\frac12 Q)\geq C_4^{-1}\delta\theta_\mu(R)$ \/{\rm (}\/with
$C_4$ given by Lemma {\rm \ref{lemprop1} (a)).} If $b\ll \delta/A${\rm ,}
then $$\mu(2Q\setminus Q)\leq \frac1{10}\,\mu(Q).$$
\end{lemma}

\Proof 
By Lemma {\rm \ref{lembalan}} there exists a square $P\subset2Q$ with
$\ell(P)\leq \ell(2Q)/10$ such that $\mu(P)\geq (1-2\cdot10^5
b)\mu(2Q)$. If $P\not\subset Q${\rm ,} then $P\subset 2Q\setminus
\frac12Q${\rm ,} and so
$$\mu(2Q\setminus \tfrac12 Q) \geq \mu(P) \geq (1-2\cdot10^5b)\mu(2Q).$$
Thus{\rm ,} $\mu(\frac12 Q) \leq 2\cdot10^5b\mu(2Q)$. If $b\ll\delta/A${\rm ,}
then we derive $\theta_\mu(\frac12Q) < C_4^{-2}\delta
A^{-1}\theta_\mu(2Q)$, which is a contradiction. Therefore{\rm ,}
$P\subset Q${\rm ,} and then
$$\mu(2Q\setminus Q) \leq \mu(2Q\setminus P) \leq 2\cdot10^5b\mu(2Q)\leq
\frac1{10}\,\mu(\tfrac12 Q),$$ since $b\ll\delta/A$.
\hfill\qed

\begin{remark} \label{remab}
 From now on, we assume that $a=1/40$, and moreover that
$b=b(A,\delta)\ll\delta/A$, so that the preceding lemma holds. For
these precise values of $a$ and $b$ we write
$\bal(\mu):=\bal_\mu(a,b)$.
\end{remark}

\begin{lemma} \label{lemhq}
Let $Q$ be a square whose center lies on $3R$ and is such that
$\ell(Q)\leq C_8^{-1}A^{-2}\ell(R)$ and
  $$\delta\theta_\mu(R)\lesssim \theta_\mu(2^nQ) \lesssim
A\theta_\mu(R)$$ for all $n\geq0$ with $\ell(2^nQ)\leq8\ell(R)$.
If $C_8$ is big enough{\rm ,} then there exists some square $\wh{Q}\in
\bal(\mu)$ concentric with $Q$ such that $2\ell(Q)\leq
\ell(\wh{Q}) \leq 8\ell(R)$ and also $\ell(\wh{Q})\leq
C_9(A,\delta)\ell(Q)$.
\end{lemma}

\Proof 
If $2^nQ\not\in \bal(\mu)$ for $n=1,\dots,N$, with
$\ell(2^NQ)\leq 8\ell(R)$, then $$\mu(2^nQ\setminus 2^{n-1}Q)\leq
\frac1{10}\mu(2^{n-1}Q) \qquad\mbox{for $n=1,\dots,N$.}$$ Thus,
$\mu(2^nQ)\leq 1.1\,\mu(2^{n-1}Q)$ for $n=1,\dots,N$, and so
\begin{eqnarray*}
\mu(2^NQ)&\leq &1.1^N\mu(Q) \leq 2^{N/2}\mu(Q) \\& \leq&
CA\theta_\mu(R)2^{N/2}\ell(Q) = CA\theta_\mu(R)2^{-N/2}\ell(2^NQ).
\end{eqnarray*}
Therefore,
\begin{equation}\label{tel2}
\theta_\mu(2^NQ)\leq CA2^{-N/2}\theta_\mu(R).
\end{equation}
   Suppose that $N$ is such
that $\ell(2^NQ) = 8\ell(R)$. Then we have $\theta_\mu(R)\simeq
\theta_\mu(2^NQ)$, and by \rf{tel2} we get $CA2^{-N/2}\geq 1$, and
so $$\ell(R)=\frac18 \,\ell(2^NQ)\leq CA^2\ell(Q) =:
\frac{C_8A^2}2\, \ell(Q),$$ which contradicts the assumptions of
the lemma. \medbreak 

We infer that there is a square $2^nQ\in\bal(\mu)$ with $n\geq 1$
and $\ell(2^nQ)\leq 8\ell(R).$ Let $\wh{Q}$ be the smallest one.
Let $N$ be such that $2^NQ=\frac12\wh{Q}$. From \rf{tel2}, since
$\theta_\mu(2^{N}Q)\gtrsim\delta\theta_\mu(R)$, we deduce
$CA\delta^{-1}2^{-N/2}\geq1$, and then $\ell(\wh{Q})\leq
CA^2\delta^{-2}\ell(Q).$
\hfill\qed

\Subsec{Construction of $K$}
From Lemma \ref{lemprop1} (b) it easily follows that for any\break\vskip-12pt\noindent 
square $Q\in \maxbad(R)$ there exists another square $\wt{Q}$
concentric with $Q$ satisfying $2\ell(Q)\leq \ell(\wt{Q})\leq
8\ell(R)$ and $\delta_\mu(Q,\wt{Q})\lesssim A\theta_\mu(R)$, such that
\begin{equation} \label{ifdens}
\mbox{if $n\geq-1$ and $\ell(2^n\wt{Q})\leq 8\ell(R)$, then}\quad
\delta\theta_\mu(R)\lesssim\theta_\mu(2^n\wt{Q})\lesssim
A\theta_\mu(R)
\end{equation}
and moreover
\begin{multline} \label{cuaq11}
\iint_{\begin{subarray}{l} y,z\in 3R\\  2^{-12}\ell(\wt{Q})<
|x-y|\leq 4\ell(R)
\end{subarray}}\!\!\!\! c(x,y,z)^2\,d\mu(y)d\mu(z) \\
\;\leq\; K_{\mu,J(\wt{Q})+12}\chi_E(x) -
K_{\mu,J(R)-2}\chi_E(x)\;\lesssim\; A^2\theta_\mu(R)^2
\end{multline}
for all $x\in \wt{Q}$.

As explained in Remark \ref{remww}, the same holds if $Q\in \sss(R)$
(instead of $\maxbad(R)$). Further, we can take the squares $\wh{Q}$ to  
be $4$-dyadic (in this way, some of the calculations below will
become simpler).
That is to say, given $Q\in \sss(R)$ there exists a $4$-dyadic square  
$\wh{Q}$ such that $Q\subset
\frac12\wh{Q}$ (we cannot assume $Q$ and $\wh{Q}$ to be concentric) with
$2\ell(Q)\leq \ell(\wt{Q})\leq
8\ell(R)$ and $\delta_\mu(Q,\wt{Q})\lesssim A\theta_\mu(R)$\break\vskip-12pt\noindent which  
satisfies \rf{ifdens} and \rf{cuaq11}.
We denote by $\qsss(R)$ the family of
squares~$\wt{Q}$, with $Q\in\sss(R)$, and we say that $\wt{Q}$ is
a {\em quasi-stopping} square of $R$.

We intend to construct some set $K$ containing $G(R)$ (remember
that $G(R)$ is the set of good points of $R$) such that, besides
other properties, for each $Q\in \qsss(R)$ there exists some point
$a_Q\in K$ with $\dist(a_Q,Q)\leq C\ell(Q)$. In the next subsection,
we will show that $K$ verifies \rf{betajones}, and thus $K$ will
be contained in an {\rm AD} regular curve $\Gamma_R$. This curve will
fulfill the required properties in Main Lemma \ref{lemcorona}.

In the next lemma we deal with the details of the construction of $K$
and selection of the points $a_Q$, $Q\in\qsss(R)$. Most
difficulties are due to the fact that the squares $Q\in\qsss(R)$
are not disjoint, in general.

\begin{lemma} \label{lemk}
Let $\eta>3$ be some fixed constant to be chosen below. For each
$x\in 3R${\rm ,}
\begin{equation} \label{eqlx} \ell_x:=\inf_{Q\in\qsss(R)}
\Bigl(\ell(Q)+\frac1{40}\,\dist(x,Q),\;\frac1{40}\, \dist(x,G(R))\Bigr).
\end{equation}
  There exists a family of points
$\{a_Q\}_{Q\in\qsss(E)}$ such that if  $$K:= G(R)
\cup\{a_Q\}_{Q\in \qsss(R)},$$ the following properties hold\/{\rm :}\/
\begin{itemize}
\item[{\rm (a)}] For each $Q\in\qsss(R)${\rm ,} $\dist(a_Q,Q)\leq C\ell(Q).$
\item[{\rm (b)}] There exists $C>0$ such that for all $x\in K${\rm ,}
$K\cap \bar{B}(x,C^{-1}\ell_x)=\{x\}$.
\item[{\rm (c)}] If $x\in K$ and $\ell_x>0${\rm ,} then \end{itemize}
$$c^2_{\mu_{|B(x,\eta\ell_x)\setminus B(x,\ell_x)}}\!(x) \leq
\frac{C(A,\delta)}{\mu(B(x,\ell_x))}\! \iiint_{\begin{subarray}{l}
y,z,t\in B(x,C\eta\ell_x) \\|y-z|\geq \ell_x
\end{subarray}}\!
c(y,z,t)^2d\mu(y)d\mu(z)d\mu(t).$$
\end{lemma}

We note  that the lemma is understood better if we think
about the points in $G(R)$ as degenerate quasi-stopping squares
with zero side length. In our construction some points $a_Q$ may
coincide for different squares $Q\in\qsss(R)$.

We also remark that if (b) were not required in the lemma, then its  
proof
would be much simpler.

\Proof  First we explain the algorithm for assigning a point
$a_Q$ to each square $Q\in\qsss(R)$. Finally we will show that
(a), (b) and (c) hold for our construction.

Take a fixed square $Q\in \qsss(R)$. Since $\ell_x$ is
a continuous (and Lipschitz) function of $x$, there exists some
$z_0\in \overline{10Q}$ such that $\ell_x$ attains its minimum
over $\overline{10Q}$ at $z_0$. If $\ell_{z_0}=0$,
$a_Q:=z_0$.

Suppose now that $\ell_{z_0}>0$. Assume first that there exists a
sequence of squares $\{P_n\}_n\subset \qsss(R)$ with
$\ell(P_n)\to0$ or points $p_n\in G(R)$ such that
\begin{equation}\label{pn0}
\ell(P_n)+ \frac1{40}\,\dist(z_0,P_n) \to\ell_{z_0}
\end{equation}
(we  identify a point $p_n$ with a square $P_n$ with
$\ell(P_n)=0$). Since $\ell_{z_0}\leq \ell(Q)$ (because if $x\in
Q$, then $\ell_x\leq \ell(Q)$ and $\ell_{z_0}\leq\ell_x$), we may
assume $P_n\subset B(z_0,41\ell(Q))\subset \overline{90Q}$. Thus,
there exists some point $z_1\in\overline{90Q}$ such that a
subsequence $\{P_{n_k}\}_k$ accumulates on $z_1$. We set
$a_Q:=z_1$. Observe that $\ell_{a_Q}=0$ in this case.

Assume now that $\ell_{z_0}>0$ and that a sequence $\{P_n\}_n$ as
above does not exist. This implies that the infimum which defines
$\ell_{z_0}$ (in \rf{eqlx} replacing $x$ by $z_0$) is attained
over a subfamily of squares $P\in\qsss(R)$ with
$\ell(P)\geq\delta$, for some fixed $\delta>0$, which further
satisfy $\dist(P,Q)\leq 41\ell(Q)$ (because
$\ell_{z_0}\leq\ell(Q)$). Then it turns out that such a subfamily
of squares must be finite, because we are dealing with $4$-dyadic
squares. Thus the infimum in \rf{eqlx} (with $x$ replaced by
$z_0$) is indeed a minimum (attained by only a finite number of
squares in $\qsss(R)$). Among the squares where the minimum is
attained, let $P_Q$ be one with minimal side length.

Let us apply Vitali's covering theorem to the family of squares
$P_Q$ obtained above. Then there exists a countable (or finite)
subfamily of pairwise disjoint squares $P_i$ (which are of the
type $P_Q$) such that
$$\bigcup P_Q \subset \bigcup_{i} 5P_i,$$
where the union on the left side is over all the squares $P_Q$
which arise in the algorithm above when one considers all the
squares $Q\in \qsss(R)$.

Now, for each $i$ we choose a ``good'' point $a_i\in
\frac12P_i\cap\supp(\mu)$, so that
\begin{multline} \label{qas11}
\iint_{\begin{subarray}{l} y,z\in 2\eta P_i\\ |a_i-y|\geq
\ell(P_i)
\end{subarray}} c(a_i,y,z)^2\,d\mu(y)d\mu(z) \\ \leq\;  
\frac1{\mu(\frac12P_i)}
\iiint_{\begin{subarray}{l} x\in \frac12P_i\\ y,z\in 2\eta P_i\\
|x-y|\geq \ell(P_i)
\end{subarray}} c(x,y,z)^2\,d\mu(x)d\mu(y)d\mu(z).
\end{multline}

We claim that for each square $Q\in\qsss(R)$ for which $a_Q$ has
not been chosen yet (which means $\ell_{z_0}>0$ and there is no
sequence $\{P_n\}_n\subset \qsss(R)$ with $\ell(P_n)\to0$
satisfying \rf{pn0}), there exists some $a_i$ such that
$\dist(a_i,Q)\leq C\ell(Q)$. Then we set $a_Q:=a_i$.

Before proving our claim, we show that if $Q\in\qsss(R)$ is such
that $Q\cap 5P_i\neq\varnothing$ for some $i$, then $\ell(Q)\geq
\ell(P_i)$. Indeed, if $\ell(Q)\leq\ell(P_i)/2$ (remember that $Q$ and  
$P_i$ are $4$-dyadic squares),
then $Q\subset10P_i$, and so it easily follows that for any $y\in\C$ we  
have
\begin{equation} \label{eqclar}
\ell(P_i) + \frac1{40}\,\dist(y,P_i)\geq\ell(Q)+ \frac1{40}\,\dist(y,Q),
\end{equation}
which is not possible because of our construction (remember that
there exists some $z_0$ such that the infimum defining
$\ell_{z_0}$ is attained by $P_i$).

Let us prove the claim now. By our construction, there exists some
square $P_Q$ with $\ell(P_Q)\leq \ell(Q)$ such that
$\dist(P_Q,Q)\leq C\ell(Q)$. Let $P_i$ be such that $5P_i\cap
P_Q\neq\varnothing$. Since $\ell(P_Q)\geq\ell(P_i)$,
\begin{eqnarray*}
\dist(a_i,Q)&\leq& \dist(a_i,P_Q) + 2^{1/2}\ell(P_Q) +
\dist(P_Q,Q) \\ &\leq & C\ell(P_Q) + 2^{1/2}\ell(P_Q) + C\ell(Q)
\leq C\ell(Q).
\end{eqnarray*}

Let us consider the statement (b). If $\ell_x=0$, there is nothing
to prove. The only points such that $\ell_x>0$ are the $a_i$'s.
Notice that $\ell_{a_i}\leq \ell(P_i)$ because $a_i\in P_i$. On
the other hand, we also have $\ell_{a_i}\geq \frac1{40}\ell(P_i)$.
Otherwise, it is easily seen that there exists some $P\in\qsss(R)$
such that $P\cap 5P_i\neq\varnothing$ and $\ell(P)\leq
\ell(P_i)/2$, which is not possible as shown above (in the paragraph of  
\rf{eqclar}). Thus (b)
follows from the fact that $a_i\in \frac12P_i$, the squares $P_i$
are disjoint, $\ell(P_i)\simeq\ell_{a_i}$, etc.

Finally (c) follows easily from \rf{qas11} and the fact that
$B(a_i,\eta\ell_{a_i}) \subset 2\eta P_i \subset
B(a_i,C\eta\ell_{a_i})$, for some $C>0$.
\hfill\qed

\Subsec{Estimate of $\sum_{P\in \DD,P\subset Q}  
\beta_K(P)^2\ell(P)$}
We will need the following result.

\begin{lemma} \label{lemlam}
There exists some $\lambda>4$ depending on $A$ and $\delta$ such
that{\rm ,} given any $Q\in \qsss(R)${\rm ,} for each $n\geq 1$ with
$\lambda^{n+1}\ell(Q)\leq\ell(R)$ there exist two squares $Q^a_n$
and $Q^b_n$  fulfilling the following properties\/{\rm :}\/
\begin{itemize}
\item[{\rm (a)}] $Q^a_n,Q^b_n \subset \frac{\lambda^{n+1}}2 Q\setminus
\lambda^{n}Q${\rm ,} 
\item[{\rm (b)}] $\dist(Q^a_n,Q^b_n)
\gtrsim\lambda^{n} \ell(Q)${\rm ,}
\item[{\rm (c)}] $\lambda^{n} \ell(Q)\lesssim \ell(Q_n^i)
\lesssim\lambda^{n+1} \ell(Q)${\rm ,} for $i=a,b${\rm ,}
\item[{\rm (d)}]  
$C(A,\delta)^{-1}\theta_\mu(R)\lesssim\theta_\mu(Q^i_n)\lesssim
A\theta_\mu(R)${\rm ,} for $i=a,b$.
\end{itemize}
\end{lemma}

\Proof 
The lemma follows easily from the existence of balanced squares
(see Lemma \ref{lemhq}). Indeed, if $Q\in\bal(\mu)$, then there
are two squares $Q_1,\,Q_2\subset Q$ fulfilling the properties
stated just above Lemma \ref{lembalan}. Since $\dist(Q_1,Q_2)\geq
\frac1{40}\ell(Q)$, one of the squares $Q_i$ is contained in
$Q\setminus 2^{-7}Q$. From this fact one can easily deduce the
existence of some constant $\lambda_0>2$ (depending on
$A,\delta,C_8,C_9,\ldots$) such that for each $n\geq1$ with
$\lambda_0^{n+1}\ell(Q)\leq\ell(R)$ there exists some square
$P_n\subset \frac12\lambda_0^{n+1}Q\setminus \lambda_0^{n}Q$
satisfying $\ell(P_n)\simeq\lambda_0^{n+1}\ell(Q)$ and
$\theta_\mu(P_n)\geq C(A,\delta)^{-1}\theta_\mu(R)$. If we set
$\lambda:=\lambda_0^2$, $Q_n^a:=P_{2n}$, and $Q_n^b:=P_{2n+1}$,
the lemma follows.
\hfill\qed

\begin{remark}
The lemma above also holds for $x\in G(R)$ (interchanged with the
square $Q\in \qsss(R)$ in the lemma) and for $x$ such that
$\ell_x=0$. That is to say, increasing $\lambda$ if necessary, for
each $n\geq1$ we have  squares $Q^a_n,Q^b_n\subset
B(x,\frac12\lambda^{-n}\ell(R))\setminus
B(x,\lambda^{-n-1}\ell(R))$ which satisfy properties analogous to
(b), (c) and (d).
\end{remark}

In the following lemma we show a version of \rf{cuaq11} which
involves the curvature $c^2(\mu)$ truncated by the function
$\ell_x$.

\begin{lemma} \label{lemcurlx}
Let $C_{10}>0$ be a fixed constant. For all $x\in 3R${\rm ,}
  $$\iint_{\begin{subarray}{l} y,z\in
3R\\ |x-y|\geq C_{10}\ell_x
\end{subarray}}\!\!\!\! c(x,y,z)^2\,d\mu(y)d\mu(z) \\
\;\leq\; C_{11} A^2\theta_\mu(R)^2,$$ where $C_{11}$ depends on
$C_{10}$.
\end{lemma}

The proof of the preceding estimate follows easily from \rf{qaz2} and  
Lemma
\ref{curpert} (as in \rf{qaz33}). We will not go through
the details.

\demo{Proof of \rf{betajones}}  We follow quite closely
Jones' ideas (see \cite[pp.\ 39--44]{Pajot}).

It is enough to show that \rf{betajones}
holds for dyadic squares $Q$ with $\ell(Q)\leq
(C_{12}\lambda)^{-1}\ell(R)$ (with $\lambda$ given by Lemma
\ref{lemlam} and $C_{12}$ to be fixed below). So we assume
$\ell(Q)\leq (C_{12}\lambda)^{-1}\ell(R)$. Also, the sum
\rf{betajones} can be restricted to those squares $P\in\DD$ such
that $P\cap K\neq\varnothing$. That is, it suffices to prove that
\begin{equation} 
\pagebreak
\label{betajones2}
\sum_{P\in \DD,\, P\subset Q,\, P\cap K\neq\varnothing}
\beta_K(P)^2\ell(P) \leq C(A,\delta)\ell(Q).
\end{equation}

The main step of the proof of \rf{betajones2} consists of
estimating $\beta_K(P)$ for some $P$ as in the sum above in terms
of $c^2(\mu)$. By standard arguments, if $\beta(P)\neq0$, there
are three pairwise different points $z_0,\,z_1,\,z_2\in K\cap 3P$
such that
\begin{equation}\label{knm5}
\beta_K(P)\simeq \frac{w(z_0,z_1,z_2)}{\ell(P)},
\end{equation}  where
$w(z_0,z_1,z_2)$ stands for the width of the thinnest infinite
strip containing $z_0,z_1,z_2$. By (b) of Lemma \ref{lemk},
$\ell_{z_0}\leq C|z_0-z_1|\leq C\ell(P)$. So either $\ell_{z_0}=0$
or there is some $P_0\in\qsss(R)$ with
$\ell(P_0)+\frac1{40}\dist(z_0,P_0)\leq C\ell(P)$. In any case, if
$C_{12}$ has been chosen big enough, there are ``many'' balanced
squares $\wh{P}$ which contain $z_0,z_1,z_2$ such that
$\ell(P)\leq\ell(\wh{P})\leq C(A,\delta)\ell(P)$.
Arguing as in Lemma \ref{lemlam} , we deduce that there exist two  
squares $S_a,\,S_b$
contained in $C(A,\delta)P\setminus P$ satisfying the properties stated  
just
above Lemma~\ref{lembalan}. For any $x_0\in S_a$ and $y_0\in S_b$,
$$w(z_0,z_1,z_2)\lesssim \dist(z_0,L_{x_0,y_0}) +
\dist(z_1,L_{x_0,y_0}) +\dist(z_2,L_{x_0,y_0}),$$ as it is easy to
check. Integrating over $x_0\in S_a$ and $y_0\in S_b$ we get
  \begin{multline*} w(z_0,z_1,z_2)\lesssim
\frac{1}{\mu(S_a)\,\mu(S_b)}\; \\
\times \iint_{\begin{subarray}{l}x_0\in S_a\\
y_0\in S_b
\end{subarray}
}\bigl[\dist(z_0,L_{x_0,y_0}) + \dist(z_1,L_{x_0,y_0})
+\dist(z_2,L_{x_0,y_0})\bigr]\,d\mu(x_0)d\mu(y_0).
\end{multline*}
Thus there exists some $z_i$, say $z_0$, such that
\begin{equation}\label{knm1}
w(z_0,z_1,z_2)\lesssim \frac{1}{\mu(S_a)\,\mu(S_b)}
\iint_{\begin{subarray}{l}x_0\in S_a\\
y_0\in S_b
\end{subarray}}
\dist(z_i,L_{x_0,y_0})\,d\mu(x_0)d\mu(y_0).
\end{equation}

 From Lemma \ref{lemlam} and the subsequent remark we deduce that
there are two families of squares $\{P^a_n\}_{n\geq1}$,  
$\{P^b_n\}_{n\geq1}$
which satisfy the following properties for any $n\geq1$ such that
$\lambda^{n+1}\leq \ell(P)/\ell_{z_0}$:
\begin{itemize}
\item[(a)] $P^a_n,P^b_n\subset
B(z_0,\frac12\lambda^{-n}\ell(P))\setminus
B(z_0,\lambda^{-n-1}\ell(P))$,

\vspace{1mm}\item[(b)] $\dist(P^a_n,P^b_n) \gtrsim\lambda^{-n-1}
\ell(P)$,

\vspace{1mm}\item[(c)] $\lambda^{-n-1} \ell(P)\lesssim \ell(P_n^i)
\lesssim\lambda^{-n} \ell(P)$, for $i=a,b$,

\vspace{1mm}
\item[(d)]  
$C(A,\delta)^{-1}\theta_\mu(R)\lesssim\theta_\mu(P^i_n)\lesssim
A\theta_\mu(R)$, for $i=a,b$,
\end{itemize}
with $\lambda>4$. We also set $P_0^a:=S_a$ and $P_0^b:=S_b$.

  Let $N$ be the biggest positive integer such
that $\lambda^{N+1}\leq \ell(P)/\ell_{z_0}$ (with $N=\infty$ if
$\ell_{z_0}=0$).
  We claim that for all points $x_n\in P^a_n$ and $y_n\in P^b_n$
  we have
\begin{equation}
\pagebreak
\label{claimn}
\dist(z_0,L_{x_0,y_0}) \leq C\sum_{n=0}^N
\bigl[\dist(x_{n+1},L_{x_n,y_n}) +
  \dist(y_{n+1},L_{x_n,y_n})\bigr],
\end{equation}
where $x_{N+1}=y_{N+1}=z_0$ if $N<\infty$, and $C$ depends
on $A,\,\delta,\,\lambda$ (like all the following constants
denoted by $C$ in the rest of the proof of \rf{betajones2}).

Assuming  the claim for the moment, from \rf{claimn} we get
\begin{eqnarray*}
\dist(z_0,L_{x_0,y_0}) &\lesssim &\sum_{n=0}^{N}
\Bigl[c(x_n,y_n,x_{n+1})\,|x_n-x_{n+1}|\,|y_n-x_{n+1}| \\
&&
\mbox{}+c(x_n,y_n,y_{n+1})\,|x_n-y_{n+1}|\,|y_n-y_{n+1}|\Bigr]\\
& \lesssim & \sum_{n=0}^N \bigl[c(x_n,y_n,x_{n+1}) +
c(x_n,y_n,y_{n+1})\bigr]\lambda^{-2n}\ell(P)^2.
\end{eqnarray*}
 From \rf{knm1}, taking the $\mu$-mean of the above inequality over
$x_0\in P_0^a=S_a$ and $y_0\in P_0^b=S_b$, then over $x_1\in P_1^a$ and
$y_1\in P_1^b$, over $x_2\in P_2^a$ and $y_2\in P_2^b$, and so
on, we obtain
\begin{eqnarray} &&\label{eqliooo}\\
  w(z_0,z_1,z_2) &\lesssim & \ell(P)^2\sum_{n=0}^{N-1}  
\lambda^{-2n}\biggl[
\frac{A_n}{\mu(P_n^a)\mu(P_n^b)\mu(P_{n+1}^a)} +
\frac{B_n}{\mu(P_n^a)\mu(P_n^b)\mu(P_{n+1}^b)}
  \biggr]\nn\\
&& + \ell_{z_0}^2
\iint_{\begin{subarray}{l} x_N\in P^a_N\\y_N\in
P^b_N
\end{subarray}}
c(x_N,y_N,z_0)\,d\mu(x_N)d\mu(y_N),
\nn
\end{eqnarray}
where $$A_n:= \iiint_{\begin{subarray}{l} x_n\in P^a_n\\y_n\in
P^b_n\\x_{n+1}\in P^a_{n+1}
\end{subarray}}\!\!\!\!
c(x_n,y_n,x_{n+1})\,d\mu(x_n)d\mu(y_n)d\mu(x_{n+1})$$ and $$B_n:=
\iiint_{\begin{subarray}{l} x_n\in P^a_n\\y_n\in P^b_n\\y_{n+1}\in
P^b_{n+1}
\end{subarray}}\!\!\!\!
c(x_n,y_n,y_{n+1})\,d\mu(x_n)d\mu(y_n)d\mu(y_{n+1}).$$
Note that the last term in \rf{eqliooo} (the one involving $\ell_{z_0}$)
only appears when $N<\infty$ (i.e.\ when $\ell_{z_0}>0$).
  By H\"older inequality, the estimates for the squares
$P_n^i$ below \rf{knm1}, and Lemma \ref{lemk} (c) (with $\eta$ big  
enough), we get
\begin{eqnarray*}
  w(z_0,z_1,z_2)&\lesssim & \ell(P)^2\sum_{n=0}^{N-1}  
\lambda^{-2n}\Biggl[
\frac{c_\mu^2(P_n^a,P_n^b,P_{n+1}^a)^{1/ 
2}}{\bigl(\mu(P_n^a)\mu(P_n^b)\mu(P_{n+1}^a)\bigr)^{1/2}}\\
&&
+\frac{c_\mu^2(P_n^a,P_n^b,P_{n+1}^b)^{1/ 
2}}{\bigl(\mu(P_n^a)\mu(P_n^b)\mu(P_{n+1}^b)\bigr)^{1/2}}
  \Biggr]
+ \ell(P)^2 \lambda^{-2N}  
\frac{c_\mu^2(z_0,P_N^a,P_N^b)^{1/ 
2}}{\bigl(\mu(P_N^a)\mu(P_N^b)\bigr)^{1/2}}\\
&    \lesssim & \theta_\mu(R)^{-3/2}\ell(P)^{1/2}\Biggl[ \sum_{n=0}^{N-1}
  \lambda^{-n/2}c_\mu^2(P_n^a,P_n^b,P_{n+1}^a\cup P_{n+1}^b)^{1/2} \\
&&+ 
\lambda^{-N/2} \biggl(\iiint_{\begin{subarray}{l}
x,y,z\in B(z_0,C\eta\ell_{z_0}) \\|y-z|\geq \ell_{z_0}
\end{subarray}}\! c(x,y,z)^2\,d\mu(x)d\mu(y)d\mu(z)\biggr)^{1/2}  
\Biggr].
\end{eqnarray*}
By \rf{knm5} and Cauchy-Schwartz, we obtain
\begin{equation*}
\begin{split}
  \beta_K(P)^2 \lesssim & \,\frac{1}{\theta_\mu(R)^{3}\ell(P)} \Biggl[  
\sum_{n=0}^N
\lambda^{-n/2} c_\mu^2(P_n^a,P_n^b,P_{n+1}^a\cup P_{n+1}^b)\\
&\mbox{} + \lambda^{-N/2} \iiint_{\begin{subarray}{l}
x,y,z\in B(z_0,C\eta\ell_{z_0}) \\|y-z|\geq \ell_{z_0}
\end{subarray}}\! c(x,y,z)^2\,d\mu(x)d\mu(y)d\mu(z) \Biggr].
\end{split}
\end{equation*}
  Notice that for every $x_n\in P_n^a$, we have
  $\ell_{x_n}\lesssim\lambda^{-n}\ell(P)$,
  because $\ell_{z_0}\lesssim\lambda^{-n} \ell(P)$,
   $\dist(x_n,z_0)\simeq \lambda^{-n} \ell(P)$, and $\ell_x$ is a
   Lipschitz function of $x$. Analogously, if $y_n\in P^b_n$, then
$\ell_{y_n}\lesssim\lambda^{-n}\ell(P)$. Moreover, since for each
$n$ there exists some dyadic square $S\subset P$ such that
$P_n^a\cup P_n^b\cup P_{n+1}^a\cup P_{n+1}^b$ is contained in
$3S$, with $\ell(S)\simeq\lambda^{-n}\ell(P)$, we infer that
  \begin{multline*}
\theta_\mu(R)^{3}\beta_K(P)^2 \ell(P)\\ \lesssim
  \sum_{S\subset
  P} \biggl(\frac{\ell(S)}{\ell(P)}\biggr)^{1/2}\!\!\iiint_{(x,y,z)\in  
S^*\cap
  R^\ell} \!\!c(x,y,z)^2
d\mu(x)d\mu(y)d\mu(z),
  \end{multline*}
where $S^*$ is the set of $(\xi_1,\xi_2,\xi_3)\in (3S)^3$ such
that $|\xi_2-\xi_3|\geq C^{-1}\ell(S)$, and $R^\ell$
is the set of $(\xi_1,\xi_2,\xi_3)\in (3R)^3$ such that
$|\xi_2-\xi_3|\geq C^{-1}(\ell_{\xi_2}+\ell_{\xi_3})$.

Now, for a fixed dyadic square $Q$ with $\ell(Q)\leq
(C_{12}\lambda)^{-1}\ell(R)$, by Lemma \ref{lemcurlx}, we get
(with the sums over $P$ and $S$ only for dyadic squares)
\begin{multline*}
\theta_\mu(R)^{3}\sum_{P:P\subset Q}\beta_K(P)^2 \ell(P) \\
\begin{split}
\lesssim &\; \sum_{S:S\subset Q} \iiint_{(x,y,z)\in S^*\cap
R^\ell} \!\!c(x,y,z)^2 d\mu(x)d\mu(y)d\mu(z)
  \sum_{P:S\subset
P\subset Q}\biggl(\frac{\ell(S)}{\ell(P)}\biggr)^{1/2}\\
\lesssim &\; \sum_{S:S\subset Q} \iiint_{(x,y,z)\in S^*\cap
R^\ell} \!\!c(x,y,z)^2 d\mu(x)d\mu(y)d\mu(z)\\
\lesssim &\;\iiint_{(x,y,z)\in (3Q)^3\cap R^\ell} \!\!c(x,y,z)^2
d\mu(x)d\mu(y)d\mu(z) \\ \lesssim &\; \theta_\mu(R)^2\mu(3Q)
\;\lesssim\; \theta_\mu(R)^3\ell(Q).
\end{split}
\end{multline*}

It only remains to prove \rf{claimn}. Indeed, suppose for
simplicity that\break $N<\infty$. We set $a_{N+1}:=z_0$, and for $0\leq
n \leq N$, let $a_n$ be the orthogonal projection of $a_{n+1}$
onto $L_{x_n,y_n}$. Then
\begin{equation}\label{fbb22}
  \dist(z_0,L_{x_0,y_0}) \leq \sum_{n=0}^N \dist(a_n,a_{n+1}).
\end{equation}
  Let us check that $a_n\in B(z_0,\lambda^{-n}\ell(R))$ for $0\leq
n\leq N+1$. We argue by (backward) induction. The statement is
clearly true for $a_{N+1}$. Suppose now that $a_{n+1}\in
B(z_0,\lambda^{-n-1}\ell(R))$. By construction, we have
$x_n,y_n\in B(z_0,\frac12\lambda^{-n}\ell(R))$. On the other hand,
since we are assuming $\lambda>4$, then $a_{n+1}\in
B(z_0,\frac12\lambda^{-n}\ell(R))$ too. Then, by elementary
geometry, $$a_n\in
B(z_0,2^{-1/2}\lambda^{-n}\ell(R))\subset
B(z_0,\lambda^{-n}\ell(R)).$$

Since $a_{n+1},\,x_{n+1},\,y_{n+1}\in
B(z_0,\lambda^{-n-1}\ell(R))$, they are collinear, and
$$|x_{n+1}-y_{n+1}|\gtrsim \lambda^{-n-1}\ell(R),$$ we deduce
$$\dist(a_{n},a_{n+1}) = \dist(a_{n+1},L_{x_n,y_n}) \lesssim
\dist(x_{n+1},L_{x_n,y_n}) + \dist(y_{n+1},L_{x_n,y_n}).$$ Thus
\rf{claimn} follows from \rf{fbb22}. \hfill\qed

\section{The packing condition for the top squares} \label{secpack}
\vglue-12pt

\Subsec{The family $\sss_{\max}^{1/2}(R)$}
In order to prove the packing condition $$\sum_{R\in \ttt(E)}
\theta_\mu(R)^2 \mu(R) \leq C (\mu(E) + c^2(\mu))$$ (with $C$
depending on $A,\,\delta,\,\ve_0,\ldots$), we need to introduce
some auxiliary families of squares. Let $\sss_{\max}(R)$ be the
subfamily of those squares $Q\in\sss(R)$ such that there does not
exist another square $Q'\in\ttt(E)$, with
$Q'\cap3R\neq\varnothing$ and $\ell(Q')\leq\ell(R)/8$, such that
$Q\in \sss(Q')$. So $\sss_{\max}(R)$ is a maximal subfamily of
$\sss(R)$ in a sense. Notice, in particular, that if
$Q\in\sss_{\max}(R)$, then $Q\not\in\sss(Q')$ for any
$Q'\in\sss(R)$.
  We also denote by
$\sss_{\max}^{1/2}(R)$ the subfamily of the squares
$Q\in\sss_{\max}(R)$ such that $4Q\cap\frac12\,R\neq\varnothing$.

\begin{lemma}\label{lem12q}
For every $R\in \ttt(E)${\rm ,} there is  $\sum_{Q\in\sss(R)}
\chi_{\frac12Q}\leq C.$
\end{lemma}

\Proof 
Suppose that $\frac12Q\cap\frac12Q'\neq\varnothing$ for
$Q,Q'\in\sss(R)$. If $\ell(Q')\leq\ell(Q)/4$, then $Q'\subset Q$,
which contradicts the definition of $\sss(R)$. Thus,
$\ell(Q')\geq\ell(Q)/2$, and in an analogous way, we have
$\ell(Q)\geq\ell(Q')/2$. The lemma follows from the fact that
there is a bounded number of $4$-dyadic squares $Q'$ such that
$\ell(Q)/2\leq\ell(Q')\leq2\ell(Q)$ for $Q$ fixed.
\hfill\qed

\begin{lemma}  \label{propstop}
\begin{itemize}
\item[{\rm (a)}] Suppose that $R_1,R_2\in\ttt(E)$ and
$Q\in\sss_{\max}^{1/2}(R_1)\cap\sss_{\max}^{1/2}(R_2)$.
Then{\rm ,} \end{itemize} \vglue-15pt
\begin{equation} \label{dosco3}
\ell(R_2)/4 \leq \ell(R_1)\leq 4\ell(R_2).
\end{equation}
\begin{itemize}
\item[{\rm (b)}]
There exists an absolute constant $N_0$ such that for any $Q\in
\ttt(E)$
  $$\#\{R\in\ttt(E):\,Q\in\sss_{\max}^{1/2}(R)\}\leq N_0.$$

\item[{\rm (c)}] For $R\in\ttt(E)${\rm ,}  $$\bigcup_{P\in\sss(R)}
P \subset \bigcup_{Q\in\sss_{\max}(R)} 4Q,$$ and{\rm ,} more generally{\rm ,}
if $P,R\in\ttt(E)$ are such that $P\cap3R\neq\varnothing$ and
$\ell(P)\leq\ell(R)/8${\rm ,} then there exists some square
$Q\in\sss_{\max}(R)$ such that $P \subset 4Q$.

  \item[{\rm (d)}] For $R\in\ttt(E)${\rm ,}
  $$G^{1/2}(R):= \Bigl\{x\in \tfrac12R\setminus
  \bigcup_{Q\in\sss_{\max}^{1/2}(R)}4Q\Bigr\}.$$
  Then{\rm ,} for each $x\in E${\rm ,}
  $$\#\{R\in \ttt(E):\,x\in G^{1/2}(R)\} \leq N_1,$$
  where $N_1$ is an absolute constant.

\item[{\rm (e)}] If $R_1,R_2\in\ttt(E)${\rm ,} $Q\in\sss_{\max}^{1/2}(R_1)\cap
\maxbad(R_2)${\rm ,} then
  $$\ell(R_1)/4\leq\ell(R_2)\leq 4\ell(R_1)$$
  and $\theta_\mu(R_1)\simeq\theta_\mu(R_2)$.
\end{itemize}
\end{lemma}

\Proof 
First we show (a). If
$Q\in\sss_{\max}^{1/2}(R_1)\cap\sss_{\max}^{1/2}(R_2)$, then
\begin{equation}\label{dosco2}
4Q\cap \tfrac12R_1\neq\varnothing \qquad\mbox{and}\qquad 4Q\cap
\tfrac12R_2\neq\varnothing.
\end{equation}  It
is easily seen that if $\ell(R_2)\leq\ell(R_1)/4$, then
\rf{dosco2} and the fact that
$\ell(Q)\leq\min(\ell(R_1),\ell(R_2))/8$ imply that $R_2\subset
R_1$. As a consequence, if $\ell(R_2)\leq\ell(R_1)/8$, then
$Q\subset R_2\subset R_1$ and so, by definition,
$Q\not\in\sss(R_1)$.

The same happens if we reverse the roles of $R_1$ and $R_2$, and
so \rf{dosco3} holds.

It is easy to check that (b) follows from (a). This is left for
the reader.

The statement (c) of the lemma follows from the fact that if there is
a sequence of squares
$Q_1,\,Q_2,\dots,Q_n = P$ with
$Q_1\in\sss_{\max}(R)$ and $Q_{j+1}\in\sss(Q_j)$ for $j\geq1$, then   
$Q_{j+1}\cap 3Q_j\neq\varnothing$
and $\ell(Q_{j+1})\leq \ell(Q_j)/8$, and so
\begin{equation*}
\begin{split}
\dist_\infty(Q_1,Q_n) + \ell(Q_n) & \leq  \sum_{j=1}^{n-1}  
\bigl[\dist_\infty(Q_j,Q_{j+1})+ \ell(Q_{j+1})\bigr]\\
& \leq \sum_{j=1}^\infty\bigl[8^{1-j} + 8^{-j}\bigr]\,\ell(Q_1) =  
\frac97\,\ell(Q_1) \leq \frac32\,\ell(Q_1),
\end{split}
\end{equation*}
which implies that $Q_n\subset 4Q_1$ ($\dist_\infty$ stands for the  
distance induced by the
norm $\|\cdot\|_\infty$).

Let us show (d) now. Let $R_1,R_2\in\ttt(E)$ be such that $x\in
G^{1/2}(R_1)\cap G^{1/2}(R_2)$. If $\ell(R_2)\leq\ell(R_1)/8$,
then by (c), $R_2$ is contained in $4Q$ for some
$Q\in\sss_{\max}(R_1)$. Since $x\in \frac12R_1 \cap \frac12R_2$,
we have $4Q\cap\frac12R_1\neq\varnothing$, and so
$Q\in\sss_{\max}^{1/2}(R_1)$, which is a contradiction. Therefore,
$\ell(R_2)>\ell(R_1)/8$. The same inequality holds interchanging
$R_1$ by $R_2$. Thus, $$\ell(R_2)/4\leq\ell(R_1)\leq 4\ell(R_2).$$
That is, all the squares $R\in\ttt(E)$ such that $x\in G^{1/2}(R)$
have comparable sizes, which implies that the number of these
squares $R$ is bounded above by some absolute constant.

Finally we will prove (e). Suppose that
$\ell(R_2)\leq\ell(R_1)/8$. From
\begin{equation}
\label{dosco5}
4Q\cap\frac12R_1\neq\varnothing\qquad\mbox{and}\qquad
\ell(Q)\leq\ell(R_2)/8,
\end{equation}
we infer that $R_2\subset 3R_1$, and so $R_2\in\sss(R_1)$. This
implies that $Q\not\in\sss_{\max}(R_1)$, which is a contradiction.
Thus, $\ell(R_2)\geq\ell(R_1)/4$.

The inequality $\ell(R_2)\leq4\ell(R_1)$ also holds. Otherwise
$\ell(R_1)\leq\ell(R_2)/8$ and $Q\subset R_1$ (by \rf{dosco5})
imply that $Q\not\in\ttt(E)$ (it should have been eliminated when
constructing $\ttt(E)$ from $\ttt_0(E)$).

The comparability between $\theta_\mu(R_1)$ and $\theta_\mu(R_2)$
follows from $4R_1\cap 4R_2\neq\varnothing$ (since $Q$ is
contained in the intersection) and $1/4\leq
\ell(R_1)/\ell(R_2)\leq4$. Indeed, one easily deduces that then
$R_1\subset 16R_2$ and $R_2\subset16R_1$, and by the doubling
properties of $R_1$ and $R_2$, one gets
$\theta_\mu(R_1)\simeq\theta_\mu(R_2)$.
\Endproof\vskip4pt

The next result is a consequence of Lemma \ref{lemprop01} and the  
statement (e) in the preceding
lemma.

\begin{lemma} \label{lemprop4}
If $R\in\ttt(E)$ and $Q\in\sss_{\max}^{1/2}(R)${\rm ,} then\/{\rm :}\/
\begin{itemize}
\item[{\rm (a)}] If $Q\in {\rm HD}(R)${\rm ,} then $\theta_\mu(Q) \gtrsim
A\,\theta_\mu(R)$.

\item[{\rm (b)}] If $Q\in {\rm HC}(R)${\rm ,} then \end{itemize}\vglue-18pt
\begin{equation}
\label{fvc1} \mu\bigl\{x\!\in \tfrac12Q\!:K_{\mu,J(Q)+12}\chi_E(x)
- K_{\mu,J(R)-6}\chi_E(x)\!\geq C^{-1}\ve_0\theta_\mu(R)^2\bigr\}
\geq C_{13}^{-1}\mu(Q).
\end{equation}
\begin{itemize}
\item[{\rm (c)}] If $Q\in {\rm LD}(R)${\rm ,} then there exists some square $S_Q$
such that $Q\subset\frac1{20}S_Q${\rm ,} with $\ell(S_Q)\leq \ell(R)/2${\rm ,}
and $\theta_\mu(S_Q)\leq C_{14}\delta\theta(R)$. Also{\rm ,}\end{itemize}\vglue-18pt
\begin{equation}
\label{fvc2}\quad  \mu\bigl\{x\in \tfrac12Q:\,K_{\mu,J(Q)+12}\chi_E(x) -
K_{\mu,J(R)-2}\chi_E(x)\leq C\ve_0\theta_\mu(R)^2\bigr\} \geq
C_{15}^{-1}\mu(Q).
\end{equation}
\end{lemma}

\Proof 
Let $R_2\in\ttt(E)$ be such that $Q\in\maxbad(R_2)$.

If $Q\in {\rm HD}(R)$, then $Q\in {\rm HD}_0(R_2)$ by definition, and by (b) in  
Lemma \ref{lemprop01} and (e)
in the preceding lemma,  $\theta_\mu(Q)\gtrsim
A\theta_\mu(R_2)\simeq A \theta_\mu(R).$

If $Q\in {\rm HC}(R)$, then $Q\in {\rm HC}_0(R_2)$. Inequality \rf{fvc1}
follows from Lemma \ref{lemprop01} (c), and the fact that
$\theta_\mu(R)\simeq\theta_\mu(R_2)$ and $|J(R)-J(R_2)|\leq2$ by (e) in  
the preceding lemma.

The statement (c) also follows easily from Lemma
\ref{propstop} (e) and the definition of ${\rm LD}_0(R)$ and Lemma  
\ref{lemprop01} (d).
\Endproof\vskip4pt

Notice that, since $\frac12R$ is doubling,
\begin{equation}
  \label{des3}
  \mu(R)\,\lesssim\,\mu(\tfrac12\,R) \,\leq \,\mu\Bigl(
\bigcup_{Q\in\sss_{\max}^{1/2}(R)}\!\!\! 4Q\Bigr) +
\mu\Bigl(\tfrac12\,R\setminus\bigcup_{Q\in\sss_{\max}^{1/2}(R)}\!\!\!
4Q\Bigr).
\end{equation}

We distinguish a special kind of square $R\in\ttt(E)$. We set
$R\in VC(E)$ (and we say that $\mu$ is {\em very concentrated} on
$R$) if $$ \mu\Bigl(\bigcup_{Q\in \sss_{\max}^{1/2}(R)\cap
{\rm HD}(R)}\!\!\!\! 4Q\Bigr)
> \frac12 \,\mu(\tfrac12R).$$

\Subsec{Squares with $\mu$ very concentrated} \label{sub53}
For $R\in VC(E)$, using the doubling properties of $\frac12R$ and
$Q\in\sss_{\max}^{1/2}(R)$, we get
\begin{eqnarray*} \mu(R)&\lesssim &  \mu(\tfrac12R) \,\leq \,2
\mu\Bigl(\bigcup_{Q\in \sss_{\max}^{1/2}(R)\cap {\rm HD}(R)}
4Q\Bigr)\\
&\lesssim &
\sum_{Q\in \sss_{\max}^{1/2}(R)\cap {\rm HD}(R)} \mu(Q),
\end{eqnarray*}
and since $\theta_\mu(Q) \gtrsim A \theta_\mu(R)$ for $Q\in
\sss_{\max}^{1/2}(R)\cap {\rm HD}(R)$,
\begin{eqnarray*}
\theta_\mu(R)^2\mu(R) &\lesssim& \frac{1}{A^2} \sum_{Q\in
\sss_{\max}^{1/2}(R)\cap {\rm HD}(R)} \theta_\mu(Q)^2\mu(Q) \\
&\lesssim &\frac{1}{A^2}\,\sum_{Q\in\sss_{\max}^{1/2}(R)}
\theta_\mu(Q)^2\mu(Q).
\end{eqnarray*}
Then by (b) of Lemma \ref{propstop},
\begin{eqnarray*}
\sum_{R\in\ttt(E)\cap VC(E)} \theta_\mu(R)^2\mu(R) &\leq &
\frac{C}{A^2}\sum_{R\in\ttt(E)}\sum_{Q\in\sss_{\max}^{1/2}(R)}
\theta_\mu(Q)^2\mu(Q) \\
&\leq &\frac{C_{16}N_0}{A^2} \sum_{Q\in\ttt(E)}
\theta_\mu(Q)^2\mu(Q).
\end{eqnarray*}
If we choose $A$ such that $C_{16}N_0/A^2\leq 1/2$, we deduce (see
Remark \ref{sumfinita} below)
\begin{equation} \label{trucc}
\sum_{R\in\ttt(E)} \theta_\mu(R)^2\mu(R) \leq
2\sum_{R\in\ttt(E)\setminus VC(E)} \theta_\mu(R)^2\mu(R).
\end{equation}

\Subsec{Squares with $\mu$ not very concentrated}
If $R\not\in VC(E)$, then
  $$\mu\Bigl(\bigcup_{Q\in
\sss_{\max}^{1/2}(R)\cap {\rm HD}(R)}\!\!\!\! 4Q\Bigr) \leq \frac12
\,\mu(\tfrac12R),$$
  and by \rf{des3} we get
  $$\mu(\tfrac12\,R) \leq 2\mu\Bigl( \bigcup_{Q\in\sss_{\max}^{1/2}(R)
  \setminus {\rm HD}(R)}\!\!\! 4Q\Bigr) +
   2\mu\Bigl(\tfrac12\,R\setminus\bigcup_{Q\in\sss_{\max}^{1/2}(R)}\!\!\! 
4Q\Bigr).$$
  Therefore,
  \begin{eqnarray}\label{des4}
  \mu(R)& \leq & C_{17}\mu\Bigl( \bigcup_{Q\in\sss_{\max}^{1/2}(R)\cap
  {\rm LD}(R)}\!\!\! 4Q\Bigr) + C_{17}\mu\Bigl(
  \bigcup_{Q\in\sss_{\max}^{1/2}(R)\cap {\rm HC}(R)}\!\!\! 4Q\Bigr)
    \\ &&\mbox{}+
   C_{17}\mu\Bigl(\tfrac12\,R\setminus\bigcup_{Q\in\sss_{\max}^{1/2}(R)}\! 
\!\!4Q\Bigr).\nonumber
  \end{eqnarray}

  We will show in Section \ref{seclds} below that
  \begin{equation} \label{eqld1}
  \mu\Bigl( \bigcup_{Q\in\sss_{\max}^{1/2}(R)\cap
  {\rm LD}(R)}\!\!\! 4Q\Bigr) \leq \eta\mu(R),
  \end{equation}
  with $\eta\leq1/(2C_{17})$, and with  $\delta$ and $\ve_0$
  chosen appropriately. Thus,
  \begin{multline} \label{eqld2}
  \sum_{R\in\ttt(E)\setminus VC(E)}\theta_\mu(R)^2
\mu\Bigl( \bigcup_{Q\in\sss_{\max}^{1/2}(R)\cap
  {\rm LD}(R)}\!\!\! 4Q\Bigr)
  \\ \leq \eta
  \sum_{R\in\ttt(E)\setminus VC(E)}\theta_\mu(R)^2 \mu(R).
  \end{multline}

  Also, in Section \ref{sechcs} we will prove that
  \begin{equation} \label{eqhc}
  \sum_{R\in\ttt(E)}\,\sum_{Q\in \sss_{\max}^{1/2}(R)\cap {\rm HC}(R)}
  \theta_\mu(R)^2\mu(Q)\leq C_{18} c^2(\mu),
  \end{equation}
  with $C_{18}$ possibly depending on $A$, $\delta$, and $\ve_0$.

  Now we deal with the term
  $$\sum_{R\in\ttt(E)}  
\mu\Bigl(\tfrac12\,R\setminus\bigcup_{Q\in\sss_{\max}^{1/2}(R)}
  \!\!\!4Q\Bigr) = \sum_{R\in\ttt(E)}\mu(G^{1/2}(R)).$$
  From (d) of Lemma \ref{propstop} we deduce
  $$\sum_{R\in\ttt(E)}\mu(G^{1/2}(R))\leq N_1 \mu(E).$$
  Therefore,
  \begin{equation} \label{eqgq}
  \sum_{R\in\ttt(E)}\theta_\mu(R)^2
\mu\Bigl(\frac12\,R\setminus\bigcup_{Q\in\sss_{\max}^{1/2}(R)}
  \!\!\!4Q\Bigr)  \leq C_0^2N_1\,\mu(E).
  \end{equation}

  From \rf{des4}, \rf{eqld2}, \rf{eqhc} and \rf{eqgq}, we get
  \begin{eqnarray*}
  \sum_{R\in\ttt(E)\setminus VC(E)} \theta_\mu(R)^2\mu(R) &\leq&
C_{17}\eta\sum_{R\in\ttt(E)\setminus VC(E)}\theta_\mu(R)^2\mu(R)
\\ &&\mbox{} + C c^2(\mu) +
C_{17} C_0^2N_1\,\mu(E).
  \end{eqnarray*}
   Therefore, since $\eta<1/(2C_{17})$,
  $$\sum_{R\in\ttt(E)\setminus VC(E)} \theta_\mu(R)^2\mu(R) \,\leq\,
  C_{18} \bigl[c^2(\mu) +\mu(E)\bigr],$$
  where $C_{18}$ depends on $A$ and $\ve_0$ (see the remark below).
  So \rf{pack} follows from this estimate and \rf{trucc}.

\begin{remark} \label{sumfinita}
The arguments above work if one assumes {\it a priori\/} that
$$\sum_{R\in\ttt(E)} \theta_\mu(R)^2\mu(R)<\infty.$$
To circumvent this difficulty it is necessary to argue more
carefully. For example, we can operate with finite subsets of
$\ttt(E)$. Let $\ttt_n(E)$ be the subfamily of $\ttt(E)$ of those
squares with side length $\geq 2^{-n}$; and let $A_n(E)$ be the
family $$\bigl\{Q\in\ttt(E)\setminus \ttt_n(E):\exists R\in
\ttt_n(E) \mbox{ such that }Q\in\sss_{\max}^{1/2}(R)\bigr\}.$$
Then, it can be checked that a slight modification of the
preceding estimates yields
\begin{eqnarray*}
\sum_{R\in\ttt_n(E)}\theta_\mu(R)^2\mu(R) &\lesssim&
\Bigl(\frac{1}{A^2} + \eta\Bigr)
\sum_{R\in\ttt_n(E)}\theta_\mu(R)^2\mu(R) \\
&&\mbox{}+ \sum_{R\in A_n(E)}\theta_\mu(R)^2\mu(R) + c^2(\mu) +
\mu(E).
\end{eqnarray*}
If we choose $A$ big enough and $\eta$ sufficiently small, we
obtain
\begin{equation}\label{eqw74}
\sum_{R\in\ttt_n(E)}\theta_\mu(R)^2\mu(R) \,\lesssim \sum_{R\in
A_n(E)}\theta_\mu(R)^2\mu(R) + c^2(\mu) + \mu(E).
\end{equation}
It can be shown that $\sum_{R\in A_n(E)}\chi_{\frac12R}\leq C$
(this is left for the reader). Then $$\sum_{R\in
A_n(E)}\theta_\mu(R)^2\mu(R)\leq C_0^2 \sum_{R\in
A_n(E)}\mu(\tfrac12R) \lesssim C\mu(E),$$ and   from
\rf{eqw74} we deduce $$\sum_{R\in\ttt_n(E)}\theta_\mu(R)^2\mu(R)
\,\lesssim  c^2(\mu) + \mu(E), $$ uniformly on $n$.
\end{remark}
 
\section{Estimates for the high curvature squares} \label{sechcs}
 \vglue-12pt
\Subsec{The class $\wh{\ttt}(E)$}
To prove \rf{eqhc} it will be simpler to use dyadic squares than
$4$-dyadic squares.

\begin{lemma} \label{bonqua}
Let $R\in\ttt(E)$ and $Q\in \sss_{\max}^{1/2}(R)\cap {\rm HC}(R)$. There
exists a dyadic square $\wh{Q}\subset \frac12 Q${\rm ,} with
$\ell(\wh{Q})= \ell(Q)/4${\rm ,} such that $\mu(\wh{Q})\geq
C_{19}^{-1}\mu(Q)$ and \begin{equation} \label{jc1} \int_{\wh{Q}}
\bigl(K_{\mu,J(\wh{Q})+10} \chi_E - K_{\mu,J(R)-6} \chi_E
\bigr)d\mu \gtrsim \ve_0\theta_\mu(R)^2\mu(\wh{Q}).
\end{equation}
\end{lemma}

\Proof 
Let $P_1,\dots,P_4$ be the disjoint dyadic squares with side
length $\ell(Q)/4$ such that $\frac12Q=\bigcup_{1\leq i
\leq4}P_i$. Remember that $$\mu\bigl\{x\!\in
\tfrac12Q\!:K_{\mu,J(Q)+12}\chi_E(x) - K_{\mu,J(R)-6}\chi_E(x)\geq
C^{-1} \ve_0\theta_\mu(R)^2\bigr\} \geq C_{13}^{-1}\mu(Q).$$
  Let $\wh{Q}$ be the square $P_i$ such that
$$\mu\bigl\{x\in P_i:\,K_{\mu,J(Q)+12}\chi_E(x) -
K_{\mu,J(R)-6}\chi_E(x)\geq C^{-1}\ve_0\theta_\mu(R)^2\bigr\}$$ is
maximal. Clearly, $\wh{Q}$ satisfies \rf{jc1}, and
$\mu(\wh{Q})\geq C_{13}^{-1}\mu(Q)/4$.
\Endproof\vskip4pt

For each square $Q\in \sss_{\max}^{1/2}(R)\cap {\rm HC}(R)$, with
$R\in\ttt(E)$, we choose a dyadic subsquare $\wh{Q}$ of $Q$ as in
the lemma. In the following Subsections \ref{secdecmu} and
\ref{secproof} we denote by $\wh{\ttt}(E)$ the class made up of
all the chosen subsquares $\wh{Q}$, and $\wh{R}_0$ (which is the
dyadic subsquare of $\frac12R_0$ with side length $\ell(R_0)/4$
which contains $E$). Notice that it may happen that
$\#\wh{\ttt}(E)<\#\ttt(E)$, because not all the squares in
$\ttt(E)$ are high curvature squares.
 
\Subsec{Decomposition of $c^2(\mu)$} \label{secdecmu}
We denote the class of all dyadic squares contained in $\wh{R}_0$
by $\Delta$, and the class of dyadic squares contained in
$\wh{R}_0$ with side length $2^{-j}$, by $\Delta_j$.

Given $Q\in\wttt(E)$, let $\wh{G}(Q)$ be the set of points $x\in
Q$ which do not belong to any square $P\in \wttt(E)$, with
$P\subsetneq Q$. Let us denote by ${\rm Term}(Q)$ the family of
{\em maximal} dyadic squares $P\in \wttt(E)$, with $P\subsetneq
Q$. Finally, we let $\tree(Q)$ be the class of dyadic squares
contained in $Q$, different from $Q$, which contain either a point
$x\in \wh{G}(Q)$ or a square from ${\rm Term}(Q)$. The squares in
${\rm Term}(Q)$ are called {\em terminal} squares of the tree
$\tree(Q)$, for obvious reasons. Notice that we have
$$\Delta = \{\wh{R}_0\}\cup\bigcup_{Q\in \wttt(E)} \tree(Q),$$
and that $\tree(Q) \cap \tree(R)=\varnothing$ if $Q\neq R$.

Given $Q\in\wttt(E)$ with $Q\neq \wh{R}_0$, we denote by $\roo(Q)$
the square $R$ such that $Q$ is a terminal square of $\tree(R)$.

We split the curvature $c^2(\mu)$ as follows:
\begin{eqnarray*}
c^2(\mu) & \gtrsim& \sum_{j} \sum_{Q\in \Delta_j} \int_Q
(K_{\mu,j+10} \chi_E -
K_{\mu,j-11} \chi_E)\, d\mu \\
& \gtrsim &  \sum_{R\in \wttt(E)} \sum_{Q\in \tree(R)} \int_Q
(K_{\mu,J(Q)+10} \chi_E - K_{\mu,J(Q)-11} \chi_E)\, d\mu.
\end{eqnarray*}
Observe that if $P\in{\rm Term}(R)$ and $x\in P$, then
  \begin{multline*}
  \sum_{Q\in \tree(R)} \chi_Q(x)\bigl(K_{\mu,J(Q)+10} \chi_E(x) -  
K_{\mu,J(Q)-11}
  \chi_E(x) \bigr) \\* \geq K_{\mu,J(P)+10} \chi_E(x) - K_{\mu,J(R)-10}
  \chi_E(x).
  \end{multline*}
Therefore,
\begin{eqnarray} \label{eq**}
c^2(\mu) & \gtrsim&  \sum_{R\in \wttt(E)} \sum_{Q\in {\rm
Term}(R)} \int_Q \bigl(K_{\mu,J(Q)+10} \chi_E - K_{\mu,J(R)-10}
\chi_E\bigr)\, d\mu\\
& = & \sum_{Q\in \wttt(E),Q\neq \wh{R}_0}  \int_Q
\bigl(K_{\mu,J(Q)+10} \chi_E - K_{\mu,J(\roo(Q))-10}
\chi_E\bigr)\, d\mu. \nonumber
\end{eqnarray}

\Subsec{Proof of \rf{eqhc}} \label{secproof}
By Lemma \ref{bonqua} (we use the same notation as in the lemma),
we have
\begin{multline*}
\sum_{R\in\ttt(E)}\,\sum_{Q\in \sss_{\max}^{1/2}(R)\cap {\rm HC}(R)}
  \theta_\mu(R)^2\mu(Q)\\
  \begin{split}
  \lesssim & \; \sum_{R\in\ttt(E)}\,\sum_{Q\in \sss_{\max}^{1/2}(R)\cap  
{\rm HC}(R)}
  \theta_\mu(R)^2\mu(\wh{Q}) \\
  \lesssim &\; N_0\sum_{P\in\wttt(E),P\neq  
\wh{R}_0}\theta_\mu(R_P)^2\mu(P),
  \end{split}
\end{multline*}
where $N_0$ is the constant appearing in (a) of Lemma
\ref{propstop}, and $R_P\in\ttt(E)$ is a square such that
$P=\wh{P}_1$ and $P_1\in\sss_{\max}^{1/2}(R_P)$ for some $P_1$. The  
square
$R_P$ is not unique, but in any case remember that if
$P_1\in\sss_{\max}^{1/2}(R_P^1)\cap\sss_{\max}^{1/2}(R_P^2)$, then
$\theta_\mu(R_P^1)\simeq\theta_\mu(R_P^2)$. By \rf{jc1},
\begin{multline}
\label{jc3} \sum_{R\in\ttt(E)}\,\sum_{Q\in
\sss_{\max}^{1/2}(R)\cap {\rm HC}(R)}
  \theta_\mu(R)^2\mu(Q) \\ \lesssim \ve_0^{-1}\sum_{P\in\wttt(E),P\neq  
\wh{R}_0}
\int_{P} \bigl(K_{\mu,J(P)+10} \chi_E - K_{\mu,J(R_P)-6} \chi_E
\bigr)d\mu.
\end{multline}

For every $P\in\wttt(E)$ different from $R_0$, we have
$\ell(R_P)\leq 16\ell(\roo(P))$. This is clear if
$\roo(P)=\wh{R}_0$. For $\roo(P)\neq \wh{R}_0$, let
$P_1,R_1\in\ttt(E)$ be such that $P=\wh{P}_1$, $\roo(P)=\wh{R}_1$.
It is easily seen that $P_1\subsetneq R_1$. If
$\ell(R_P)/16>\ell(\roo(P))=\ell(R_1)/4$, then
$P_1\not\in\sss(R_P)$, by the definition of the family
$\sss(\cdot)$ (since $P_1\subset R_1$ and
$\ell(R_1)\leq\ell(R_P)/8$). Thus, $$J(R_P)\geq J(\roo(P))-4,$$
and so $$K_{\mu,J(P)+10} \chi_E - K_{\mu,J(R_P)-6}\chi_E\leq
K_{\mu,J(P)+10} \chi_E - K_{\mu,J(\roo(P))-10}\chi_E.$$
  From \rf{jc3}, \rf{eq**}, and the preceding estimate we get
\vglue12pt
\hfill $\displaystyle{\sum_{R\in\ttt(E)}\,\sum_{Q\in \sss_{\max}^{1/2}(R)\cap {\rm HC}(R)}
  \theta_\mu(R)^2\mu(Q)\lesssim \ve_0^{-1}\,c^2(\mu).}$ 
\fiproof

\section{Estimates for the low density squares} \label{seclds}

To prove the packing condition \rf{pack} it remains to show
that
  \begin{equation}
\label{eqsj} \mu\Bigl( \bigcup_{Q\in\sss_{\max}^{1/2}(R)\cap
  {\rm LD}(R)}\!\!\! 4Q\Bigr) \leq \eta\mu(R),
  \end{equation}
with $\eta= 1/(2C_{17})$ (notice that $\eta$ is an absolute
constant and it does not depend either on $A$ or $\delta$).

\Subsec{The big and small squares $S_j${\rm ,} $j\in I_{{\rm LD}(R)}$}
\label{subsub1}
 For each $x\in 4Q$,
$Q\in\sss_{\max}^{1/2}(R)\cap {\rm LD}(R)$, let $S_x$ be a square such
that $x\in \frac15S_x$, $\theta_\mu(S_x)\leq
C_{14}\delta\theta_\mu(R)$, and $\ell(S_x)=2^{-m}\ell(R)$ with $m
\geq1$ (this square exists because of (c) in Lemma
\ref{lemprop4}). Moreover, we assume that $S_x$ has maximal side
length among all the squares with these properties.

  Let $\bigcup_{j\in I_{{\rm LD}(R)}} S_j$ be a
Besicovitch covering of $\bigcup_{Q\in\sss_{\max}^{1/2}(R)\cap
  {\rm LD}(R)}\! 4Q$, with $S_j:=S_{x_j}$ as explained above.

\begin{lemma} \label{lemmu8}
There exist $n_0\geq4$ and $C_{20}>0$ such that
if $\ell(S_j)\leq C_{20}^{-1}\ell(R)$ for
$j\in I_{{\rm LD}(R)}${\rm ,} then $\ell(64S_j)\leq \ell(2^{n_0}S_j)\leq \ell(R)$  
and
$\mu(2^{n_0}S_j)\geq2\mu(S_j)$.
\end{lemma}

\Proof 
  Given $n_0$ such that $\ell(64S_j)\leq\ell(2^{n_0}S_j)\leq \ell(R)$,
  we have
  $$\mu(2^{n_0}S_j)\gtrsim \delta\theta_\mu(R)\ell(2^{n_0}S_j)
  =
2^{n_0}\delta\theta_\mu(R)\ell(S_j) \gtrsim 2^{n_0}\mu(S_j).$$
  Thus, for $n_0$ big enough, we have $\mu(2^{n_0}S_j)\geq2\mu(S_j)$. So
  the lemma follows if $C_{20}$
  is big enough too (so that $\ell(2^{n_0}S_j)\leq \ell(R)$).
\Endproof\vskip4pt

To prove \rf{eqsj} we will distinguish two types of squares $S_j$.
If $S_j$, $j\in I_{{\rm LD}(R)}$, satisfies
\begin{equation} \label{eqbig}
\ell(S_j) \geq \min\biggl(C_8^{-1}A^{-2},C_{20}^{-1}
\biggr)\,\ell(R),
\end{equation}
where $C_8$ is as  defined in Lemma \ref{lemhq}, and $C_{20}$ in
Lemma \ref{lemmu8}, then we write $j\in I_{{\rm LD}(R)}^b$, and
otherwise we set $j\in I_{{\rm LD}(R)}^s$ (the superindices ``b'' and
``s'' stand for ``big'' and ``small'' respectively).

Next we estimate the measure $\mu$ of the family of the big
squares $S_j$:

\begin{lemma} \label{quagran}
$$\mu\Bigl(\bigcup_{j\in I_{{\rm LD}(R)}^b} S_j\Bigr) \lesssim
A^{4}\delta\mu(R).$$
\end{lemma}

\Proof 
We set $C_{21}:=\max(C_8,C_{20})$. Then each square $S_j$, $j\in
I^b_{{\rm LD}(R)}$, satisfies $\ell(S_j) \geq C_{21}^{-1}
A^{-2}\ell(R).$ Since the family $\{S_j\}_{j\in I^b_{{\rm LD}(R)}}$ has
finite superposition we have $$\# I^b_{{\rm LD}(R)} \lesssim
\Bigl(\frac{\ell(R)}{\inf_{j\in I^b_{{\rm LD}(R)}} \ell(S_j)}\Bigr)^2
\lesssim A^4.$$ Therefore,
\vglue12pt\hfill $\displaystyle{\sum_{j\in I^b_{{\rm LD}(R)}}\!\!\mu(S_j) \lesssim
\!\!\sum_{j\in I^b_{{\rm LD}(R)}}\!\!\! \delta\ell(S_j)\theta_\mu(R)
\lesssim \delta\mu(R)\cdot \# I^b_{{\rm LD}(R)} \lesssim
A^{4}\delta\mu(R).}$ 
\hfill\qed
\vglue8pt

\Subsec{Estimates for the small squares $S_j${\rm ,} $j\in I_{{\rm LD}(R)}^{s}$}
Now we turn our attention to the small squares $S_j$. For each
$Q\in\sss_{\max}^{1/2}(R)\cap {\rm LD}(R)$, let $W_Q$ be the set
  $$\{x\in
\tfrac12Q\cap\supp(\mu):\,K_{\mu,J(Q)+12}\chi_E(x) -
K_{\mu,J(R)-2}\chi_E(x)\leq C\ve_0\theta_\mu(R)^2\bigr\}.$$
Remember that, by \rf{fvc2}, $\mu(W_Q)\geq C^{-1}\mu(Q)\simeq\mu(4Q)$,  
since $Q$ is $16$-doubling.

For $j\in I_{{\rm LD}(R)}$, we denote $$W_j:= \bigcup_{
Q\in\sss_{\max}^{1/2}(R)\cap {\rm LD}(R)} W_Q\cap S_j.$$ Remember that
the family $\{S_j\}_{j\in I_{{\rm LD}(R)}^s}$ was obtained by an
application of the Besicovitch covering theorem. So $\{S_j\}_{j\in
I_{{\rm LD}(R)}^s}$ can be split into $N_B$ subfamilies of pairwise
disjoint squares $S_j$. Thus there exists a subfamily
$\{S_j\}_{j\in I_{{\rm LD}(R)}^{s,0}}$, $I_{{\rm LD}(R)}^{s,0}\subset
I_{{\rm LD}(R)}^s$, of pairwise disjoint squares such that
\begin{equation*}
\mu\biggl(\bigcup_{j\in I_{{\rm LD}(R)}^{s,0}} W_j \biggr) \geq
\frac1{N_B}\, \mu\biggl(\bigcup_{j\in I_{{\rm LD}(R)}^{s}} W_j\biggr).
\end{equation*}
We set
  $$W:= \bigcup_{j\in I_{{\rm LD}^{s,0}(R)}\,}W_j.$$

By the preceding lemma and since
$\sum_{Q\in\sss_{\max}^{1/2}(R)\cap
  {\rm LD}(R)}\chi_{\frac12Q}\leq C$ (see Lemma \ref{lem12q}), we get
\begin{eqnarray}\label{qpet1}
\mu\biggl(\bigcup_{Q\in\sss_{\max}^{1/2}(R)\cap
  {\rm LD}(R)}\!\!\! 4Q\biggr)
&\leq& \; \sum_{Q\in\sss_{\max}^{1/2}(R)\cap
  {\rm LD}(R)}\!\!\! \mu(4Q)  \\
&\lesssim& \; \sum_{Q\in\sss_{\max}^{1/2}(R)\cap
  {\rm LD}(R)}\!\!\! \mu(W_Q)\nn\\
&\lesssim& \;\mu\biggl(\bigcup_{j\in I_{{\rm LD}(R)}^b}W_j\biggr) +
\mu\biggl(\bigcup_{j\in I_{{\rm LD}(R)}^s}W_j\biggr)\nn\\
&\lesssim& \;A^{4}\delta\mu(R) + N_B\!\!\! \sum_{j\in
  I_{{\rm LD}^{s,0}(R)}}\mu(W_j).
  \nn
\end{eqnarray}

\begin{remark} \label{remcostat}
Another useful property of our construction of the squares $S_j$
is the following: {\em If $Q\in\sss_{\max}^{1/2}(R)\cap {\rm LD}(R)$ is
such that $Q\cap S_j\neq \varnothing$ \/{\rm (}\/for some $j\in I_{{\rm LD}(R)}${\rm ),}
then $\ell(Q)\leq \ell(S_j)$ and $Q\subset 3S_j$.}

Indeed, suppose that $\ell(Q)> \ell(S_j)$. By Lemma
\ref{lemprop4} there exists some square $S_Q$ such that $Q\subset
\frac1{20}S_Q$ with $\theta_\mu(S_Q)\leq
C_{14}\delta\theta_\mu(R)$, and $\ell(S_Q)=2^{-m}\ell(R)$ with $m
\geq1$. Then we have $S_j\subset 3Q \subset \frac15S_Q$, which is
not possible, because of the choice of $S_j$ with maximal size
(besides other properties). The inclusion $Q\subset 3S_j$ is a
direct consequence of the inequality $\ell(Q)\leq \ell(S_j)$ and
the fact that $Q\cap S_j\neq \varnothing$.

Similar arguments show that, {\em if $Q\in\sss_{\max}^{1/2}(R)\cap
{\rm LD}(R)$, then $$\ell(Q)\leq \dist(Q,S_j)+\ell(S_j).$$}
\end{remark}

   For $x\in W$, we set
  $$\ell_x := 2^{-12}\inf\bigl\{\ell(Q):\,Q\in \sss_{\max}^{1/2}(R)\cap  
{\rm LD}(R),\,
  x\in Q\}.$$
Notice that, by the preceding remark, we have
\begin{equation}\label{remrem1}
\ell_x\leq 2^{-12}\bigl(\dist(x,S_j) + \ell(S_j)\bigr)
\end{equation}
for each $j\in I_{{\rm LD}^{s,0}(R)}$.
  As a consequence, if $x\in S_j$, then $\ell_x\leq2^{-12}\ell(S_j)$.

We consider the following truncated version of the curvature
$c_\mu(x,2R,2R)$, for $x\in W$:
  $$ c^2_{tr,\mu}(x,2R,2R):= \iint_{\begin{subarray}{l}
y,z\in 2R\;\\
|x-y|>\ell_x
\end{subarray}}c(x,y,z)^2d\mu(y)d\mu(z).$$

The next lemma follows easily from our construction.

\begin{lemma}\label{lemcurtr}
For every $x\in W${\rm ,}  $$c^2_{tr,\mu}(x,2R,2R) \lesssim
\ve_0\theta_\mu(R)^2.$$
\end{lemma}

\Proof 
Let $Q\in\sss_{\max}^{1/2}(R)\cap {\rm LD}(R)$ be such that $x\in W_Q$.
By the definition of $W_Q$, we have
\begin{eqnarray*}
c^2_{tr,\mu}(x,2R,2R)& \leq&
\iint_{2^{-12}\ell(Q)<|x-y|\leq4\ell(R)}c(x,y,z)^2\,d\mu(y)d\mu(z)
\\ &= & K_{\mu,J(Q)+12}\chi_E(x) - K_{\mu,J(R)-2}\chi_E \lesssim  
\ve_0\theta_\mu(R)^2.
\end{eqnarray*}
\vglue-22pt
\Endproof\vskip12pt 

  For each $j\in I_{{\rm LD}(R)}^{s,0}$, let $L_j$ be a segment of length
  $\HH^1(L_j)=\ell(S_j)/8$ contained in $\frac12S_j$. The exact
  position and orientation of $L_j$ will be fixed
  below. Let $\nu$ be the following measure
  $$d\nu = \sum_{j\in I_{{\rm LD}(R)}^{s,0}} \frac{\mu(W_j)
  }{\HH^1(L_j)}\,d\HH^1_{|L_j}.$$

  \begin{lemma} \label{lemcur*}
  The measure $\nu$ satisfies
\begin{equation} \label{grnu}
  \nu(B(x,r))\lesssim CA\theta_\mu(R)r\qquad\mbox{for all $x\in\C$ and
  $r>0$},
\end{equation}
  and if the position and orientation of each $L_j$ are  chosen  
appropriately{\rm ,} also
\begin{equation}\label{curnu}
  c^2(\nu) \leq C(A,\delta)\ve_0^{1/50}\theta_\mu(R)^2\mu(R).
\end{equation}
  \end{lemma}
\vskip8pt

  We defer the proof of this lemma until Subsection \ref{secprova}.
  For the moment, we only remark that it follows from the estimate
  of $c^2_{tr,\mu}(x,2R,2R)$ for $x\in W$ in
  Lemma \ref{lemcurtr}, by comparison.

Now we recall David-L\'eger's theorem \cite{Leger} (the
quantitative version in \cite[Prop.\ 1.2]{Leger}).

\demo{\scshape Theorem D}
{\it For any $c_0>0${\rm ,} there exists some $\ve_L>0$ such that if $\tau$
is a Radon measure whose support is contained in a square $R$ and
satisfies\/{\rm :}\/ \begin{itemize}
\item[{\rm (a)}] $\tau(R) \geq \ell(R)${\rm ,}
\item[{\rm (b)}] $\tau(B(x,r))\leq
c_0r$ for any $x\in\C$, $r>0${\rm ,} and
\item[{\rm (c)}]$c^2(\tau)\leq \ve_L\ell(R)${\rm ,}
\end{itemize}
then there exists a Lipschitz graph $\Gamma$ with slope $\leq
1/10$ \/{\rm (}\/with respect to the appropriate axes\/{\rm )}\/ such that
$\tau(\Gamma)\geq\dfrac{99}{100}\,\tau(R)$.}

\Subsec{Proof of \rf{eqld1}} \label{sub73}
Suppose that
\begin{equation}
\label{supos}\mu\Bigl( \bigcup_{Q\in\sss_{\max}^{1/2}(R)\cap
  {\rm LD}(R)}\!\!\! 4Q\Bigr) > \eta\mu(R).
\end{equation}
Remember that if $\delta$ is small enough, by \rf{qpet1},
\begin{equation}
\label{qg1} \mu\biggl(\bigcup_{Q\in\sss_{\max}^{1/2}(R)\cap
  {\rm LD}(R)}\!\!\! 4Q\biggr)
\leq \frac{\eta}2\,\mu(R) + CN_B\!\!\! \sum_{j\in
  I_{{\rm LD}^{s,0}(R)}}\mu(W_j).
\end{equation}
So we only have to estimate $$\sum_{j\in
  I_{{\rm LD}^{s,0}(R)}}\mu(W_j) = \nu(R).$$
 From the assumption \rf{supos} and inequality \rf{qg1} we deduce
\begin{equation}
\label{supos2}\nu(R)\geq
C^{-1}N_B^{-1}\frac\eta2\,\mu(R)=:C_{22}^{-1}\eta\mu(R).
\end{equation}

  Considering  the measure  
$\tau:=\dfrac{C_{22}}{\eta\theta_\mu(R)}\,\nu$,
  we have
  $$\tau(R) \geq
  \frac{C_{22}}{\eta\theta_\mu(R)}\,C_{22}^{-1}\eta\mu(R)=
  \ell(R).$$
  On the other hand, any ball $B(x,r)$ satisfies
  $$\tau(B(x,r))\leq
  \frac{C_{22}}{\eta}\,CA\,r,$$
  because of the estimate on the linear growth of $\nu$ in Lemma
  \ref{lemcur*}. Further, from the same lemma we also get the following  
estimate
  for $c^2(\tau)$:
  \begin{eqnarray*}
c^2(\tau) &=
&\frac{C_{22}^3}{\eta^3\theta_\mu(R)^3}\,c^2(\nu)\\
& \leq & \frac{C_{22}^3}{\eta^3\theta_\mu(R)^3}\,
C(A,\delta)\,\ve_0^{1/50} \theta_\mu(R)^2\,\mu(R) =
\frac{C_{22}^3C(A,\delta)}{\eta^3}\, \ve_0^{1/50} \,\ell(R).
  \end{eqnarray*}
Therefore, by Theorem D, if $\ve_0$ is small enough, there exists
a Lipschitz graph $\Gamma$ with slope $\leq1/10$ such that
  $\tau(\Gamma)\geq \frac{99}{100}\,\tau(R),$
which is equivalent to saying
  $$\nu(\Gamma)\geq
\frac{99}{100}\,\nu(R).$$

Let $J$ be the subset of indices $j\in I_{{\rm LD}(R)}^{s,0}$ such that
$L_j\cap\Gamma\neq\varnothing$. Notice that if $j\in J$, we have
$\HH^1(\Gamma\cap S_j)\geq\frac12\,\ell(S_j)$ because $L_j$ is
contained in $\frac12\,S_j$. Thus, since the squares $S_j$, $j\in
J$, are disjoint, we have
  $$\sum_{j\in J}\ell(S_j) \leq 2 \HH^1(\Gamma) \leq 10\ell(R)$$
(of course, ``10'' is not the best constant here). Then we obtain
\begin{eqnarray*}
\nu(R) & \leq & \frac{100}{99}\sum_{j\in J}\nu(\Gamma\cap L_j)\; \leq \;
  \frac{100}{99}\sum_{j\in J}\nu(L_j) \; \leq \;
  \frac{100}{99}\sum_{j\in J}\mu(S_j)\\ &\lesssim &
  \delta\,\theta_\mu(R) \sum_{j\in J}\ell(S_j) \;\lesssim\;
  \delta\,\theta_\mu(R)\,\ell(R)\;= \; \delta\mu(R).
\end{eqnarray*}
Thus, if $\delta$ has been chosen small enough, we get a
contradiction to  \rf{supos2}.
  \fiproof
 \Subsec{Proof of Lemma {\rm \ref{lemcur*}}} \label{secprova}
To simplify notation, in this subsection we write
$J_0:=I_{{\rm LD}(R)}^{s,0}$.

The linear growth condition \rf{grnu} follows easily from the fact
that if $x\in S_j$, $j\in J_0$, then $\mu(B(x,r)) \lesssim A
\theta_\mu(R)r$ for $r\geq\ell(S_j)$, and also
$$\frac{\mu(W_j)}{\HH^1(L_j)} \lesssim
\theta_\mu(S_j)\lesssim\theta_\mu(R).$$ The details are left for
the reader.

The proof of \rf{curnu} is more delicate. If, instead of choosing
an appropriate orientation for each segment $L_j$, we assume all
the $L_j$'s to be parallel to the $x$ axis, say, then instead of
\rf{curnu} we would get an estimate such as $$c^2(\nu) \leq
C_{23}(A,\delta)\theta_\mu(R)^2\mu(R),$$ where  $C_{23}(A,\delta)$
is  a large constant. Unfortunately this estimate is not enough
for our purposes, because for the application of L\'eger's theorem
to the measure $\tau$ in Subsection \ref{sub73}, we need
$C_{23}(A,\delta)\leq \ve_L$.

The position and orientation of each segment $L_j$, $j\in J_0$,
will be fixed with the help of a balanced square $\wh{S}_j$
concentric with $S_j$, with $S_j\subset\wh{S}_j\subset
C(A,\delta)S_j$. We will show that $W_j$ is contained in a thin
strip $V_j$ associated to $\wh{S}_j$. The segment $L_j$ will be a
segment parallel to the strip $V_j$, with length $\ell(S_j)/8$, so
that the middle point of $L_j$ coincides with some point in
$W_j\cap \frac15S_j$.
 \Subsubsec{Preliminary lemmas}
For each $j\in J_0$, by Lemmas \ref{lemhq} and \ref{lemmu8} there
exists a square $\wh{S}_j$ concentric with $S_j$ satisfying
$$64\ell(S_j)\leq
\ell(\wh{S}_j)\leq\min\bigl(8\ell(R),C(A,\delta)\bigr)\,\ell(S_j),$$
such that $\wh{S}_j\in\bal(\mu)$ and $\mu(\wh{S}_j\setminus
S_j)\geq \frac12\, \mu(\wh{S}_j)$.

\begin{lemma} \label{lemprophsj}
For each $j\in J_0$ there exist two squares $Q_j^1,\,Q_j^2\subset
\wh{S}_j$ and an infinite strip $V_j$ of width
$\leq\ve_0^{1/6}\ell(S_j)$ which contains
$10\ve_0^{-1/50}\wh{S}_j\cap W$ such that
\begin{itemize}
\item[{\rm (a)}] $\dist(Q^1_j, Q^2_j) \geq a\ell(\wh{S}_j)${\rm ,}
\item[{\rm (b)}] $\ell(Q_j^i)\leq \dfrac{a}{10}\,\ell(\wh{S}_j)$ for $i=1,2${\rm ,}  
and
\item[{\rm (c)}] $\mu(Q_j^i\cap V_j) \geq \dfrac{b}{2}\,\mu(\wh{S}_j)$ for  
$i=1,2${\rm ,}
\end{itemize}
when  $\ve_0$ is small enough.
\end{lemma}

The constants $a$ and $b$ which appear in the lemma are the ones
in Remark~\ref{remab}.

\Proof 
Since $\wh{S}_j\in\bal(\mu)$, there are squares
$Q^1_j,Q_j^2\subset \wh{S}_j$ satisfying the properties (a) and
(b) and such that $\mu(Q_j^i)\geq b\mu(\wh{S}_j)$, $i=1,2$. In
order to show the existence of the strip $V_j$,  we will first
prove that most of $\supp(\mu)\cap \wh{S}_j$ is very close to
some line.

  Let $x_0\in W\cap \frac15S_j$ ($x_0$ exists because of the
construction of $S_j$). By \rf{remrem1}, for any
$y\in\wh{S}_j\setminus S_j$, we have $|y-x_0|> \ell_{x_0}$. Thus,
by Lemma \ref{lemcurtr},
\begin{multline*}
\iint_{\begin{subarray}{l}
y\in \wh{S}_j\setminus S_j\\
z\in \wh{S}_j
\end{subarray}}c(x_0,y,z)^2\,d\mu(y)d\mu(z) \\ \leq
\iint_{\begin{subarray}{l}
y,z\in \wh{S}_j\\
|x_0-y|>\ell_{x_0}
\end{subarray}}c(x_0,y,z)^2\,d\mu(y)d\mu(z)\lesssim  
\ve_0\theta_\mu(R)^2.
\end{multline*}
Therefore, there exists some $y_0\in\wh{S}_j\setminus S_j$ such
that
  $$\int_{z\in \wh{S}_j} c(x_0,y_0,z)^2\,d\mu(z)\lesssim
  \frac{\ve_0\theta_\mu(R)^2}{\mu(\wh{S}_j\setminus S_j)}
\leq \frac{2\ve_0\theta_\mu(R)^2}{\mu(\wh{S}_j)}.$$
  By Tchebychev, we obtain
\begin{multline*}
\mu\bigl\{z\in
\wh{S}_j:\,\dist(z,L_{x_0,y_0})>\ve_0^{1/4}\ell(\wh{S}_j)\bigr\}\\
\begin{split}
\leq & \; \frac1{4\ve_0^{1/2}}\int_{\wh{S}_j}
\biggl(\frac{2\dist(z,L_{x_0,y_0})}{\ell(\wh{S}_j)}\biggr)^2\,d\mu(z)\\
  \leq & \;
\frac{\ell(\wh{S}_j)^2}{\ve_0^{1/2}}
\int_{\wh{S}_j}\biggl(\frac{2\dist(z,L_{x_0,y_0})}{|x_0-z||z- 
y_0|}\biggr)^2\,d\mu(z)\;
  = \;\frac{\ell(\wh{S}_j)^2}{\ve_0^{1/2}}\,c^2_\mu(x_0,y_0,\wh{S}_j)\\
  \lesssim & \;
2\ve_0^{1/ 
2}\,\biggl(\frac{\theta_\mu(R)}{\theta_\mu(\wh{S}_j)}\biggr)^2\mu(\wh{S} 
_j)
\; \lesssim \;\frac{\ve_0^{1/2}}{\delta^2}\,\mu(\wh{S}_j).
\end{split}
\end{multline*}
Let $\wt{V}_j$ be the infinite strip with  axis $L_{x_0,y_0}$ and
width $2\ve_0^{1/4}\ell(\wh{S}_j)$. If $\ve_0$ is small enough, we
infer that
\begin{equation}
\label{bnv1}\mu(\wh{S}_j\setminus \wt{V}_j) \leq
  \frac{C\ve_0^{1/2}}{\delta^2}\,\mu(\wh{S}_j)\leq \frac12 \mu(Q_j^i)
\end{equation}
for $i=1,2$, since $\mu(Q_j^i)\geq b\mu(\wh{S}_j)$. Therefore,
$\mu(Q_j^i\cap \wt{V}_j)\geq\frac12\mu(Q_j^i)$ for each $i$. This
will imply the statement (c) because we will construct $V_j$ so
that $V_j\supset \wt{V}_j$.

It remains to define $V_j$ and to show that
$10\ve_0^{-1/50}\wh{S}_j\cap W\subset V_j$. Take\break\vskip-12pt\noindent  $y\in Q_j^1\cap
\wt{V}_j$ and $z\in Q_j^2\cap \wt{V}_j$. Since $\dist(y,z)\geq
a\ell(\wh{S}_j)$ (with $a=1/40$), the segment $L_{y,z}\cap
30\ve_0^{-1/50}\wh{S}_j$ is contained in some strip with the same
  axis as $\wt{V}_j$ and width
$C\ve_0^{1/4}\ell(30\ve_0^{-1/50}\wh{S}_j)\leq\ve_0^{1/ 
5}\ell(\wh{S}_j)/3$
(assuming $\ve_0$ small enough).

Let $V_j$ be the strip with the same  axis as $\wt{V}_j$ and width
$\ve_0^{1/5}\ell(\wh{S}_j)$ (which is $\leq \ve_0^{1/6}\ell(S_j)$
for $\ve_0$ small). If $x\in 10\ve_0^{-1/50}\wh{S}_j \setminus
V_j$, then $\dist(x,L_{y,z})>\ve_0^{1/5}\ell(\wh{S}_j)/3$, and so
  $$c(x,y,z)\geq
   \frac{C^{-1}\ve_0^{1/5}\ell(\wh{S}_j)}{\ell\bigl(10\ve_0^{-1/ 
50}\wh{S}_j\bigr)^2}
= C^{-1}\ve_0^{6/25}\ell(\wh{S}_j)^{-1}.$$
  Thus,
  $$c^2_\mu(x,Q_j^1,Q_j^2)\geq
  C^{-1}\ve_0^{12/25}\ell(\wh{S}_j)^{-2}\,\mu(Q_j^1)\,\mu(Q_j^2) \geq
  C(A,\delta)^{-1}\ve_0^{12/25}\theta_\mu(R)^2,$$
which is larger than $C\ve_0\theta_\mu(R)^2$ as $\ve_0$ has been
taken small enough. Further, from \rf{remrem1} it easily follows
that $\ell_x\leq 2^{-12}\bigl(\dist(x,\wh{S}_j) +
\ell(\wh{S}_j)\bigr)$, and then either $\ell_x\leq \dist(x,Q_j^1)$
or $\ell_x\leq \dist(x,Q_j^2)$. As a consequence,
  $$c^2_{tr,\mu}(x,2R,2R)\geq c^2_\mu(x,Q_j^1,Q_j^2)>C
\ve_0\theta_\mu(R)^2,$$ and so $x\not\in W$.
\Endproof\vskip4pt 

The {\it orientation of the segments $L_j${\rm ,} $j\in J_0$,} which
support $\nu$ is chosen so that each $L_j$ is supported on the
  axis of $V_j$. Remember also that $L_j$ has length $\ell(S_j)/8$.
We assume that its middle point coincides with some point in
$W\cap \frac15S_j$ (for example, the point $x_0$ appearing in the
proof of the preceding lemma). Notice that, in particular, we
have $L_j\subset\frac12S_j\cap V_j$.

  We denote $\QS_j:=\ve_0^{-1/50}\wh{S}_j$. Given two lines $L$ and $M$,
  $\measuredangle (L,M)$ stands for the angle between $L$ and $M$ (it  
does not
  matter which one of the two possible angles because we will always  
deal with
  its sinus). Also, given $x,y,z\in\C$, we set
  $\measuredangle(x,y,z):=\measuredangle(L_{x,y},L_{y,z})$.

slowly, in a sense.

\begin{lemma} \label{lemangles}
Let $S_j,\,S_k${\rm ,} $j,k\in J_0${\rm ,} be such that
$3\QS_j\cap3\QS_k\neq\varnothing$. Suppose that
$\ell(\QS_j)\geq\ell(\QS_k)$. Then{\rm ,} either
$\sin\measuredangle(L_j,L_k)\leq C\ve_0^{1/6}$ or $\ell(S_k)\leq
\ve_0^{3/5}\ell(S_j)$. In any case{\rm ,} $L_k$ is contained in a strip
with the same  axis as $V_j$ and width $\ve_0^{1/8}\ell(S_j)$.
\end{lemma}

\Proof 
First we will show that either $\sin\measuredangle(L_j,L_k)\leq
C\ve^{1/6}$ or $\ell(S_k)\leq \ve_0^{3/5}\ell(S_j)$.

 From the assumptions in the lemma we deduce that
$3\QS_k\subset9\QS_j$. By construction there exists some $x\in
W\cap L_k$. Consider the squares $Q^1_j,Q^2_j\subset \wh{S}_j$
mentioned in Lemma \ref{lemprophsj}. Suppose that $Q^1_j$ is the
one which is farther from $x$, so that $\dist(x,Q^1_j)\geq
C^{-1}\ell(\wh{S}_j)$. Take also the square $Q^i_k$, $i=1$ or $2$,
which is farther from $x$. Suppose this is $Q^1_k$, and so
$\dist(x,Q^1_k)\geq C^{-1}\ell(\wh{S}_k)$. Take $y\in Q^1_j\cap
V_j$ and $z\in Q^1_k\cap V_k$. Since $x,y\in V_j$, clearly we have
  $$\sin\measuredangle(L_j,\stackrel{\leftrightarrow}{xy})\leq
  \frac{\mbox{width of $V_j$}}{|x-y|} \leq C\ve^{1/6},$$
and, analogously, since $x\in V_k$,
$\sin\measuredangle(L_k,\stackrel{\leftrightarrow}{xz})\leq
C\ve_0^{1/6}$. So we infer that $$\sin\measuredangle(L_j,L_k)\leq
C\sin\measuredangle(y,x,z) + C\ve_0^{1/6}.$$
  Therefore,
$$c(x,y,z) = \frac{2\sin\measuredangle(y,x,z)}{|y-z|} \geq
C^{-1}\frac{\sin\measuredangle(L_j,L_k)-C\ve_0^{1/6}}{\ell(\QS_j)}.$$
Thus,
\begin{multline*}
\iint_{\begin{subarray}{l}y\in Q^1_j\cap V_j\\z\in Q^1_k\cap
V_k\end{subarray}} c(x,y,z)^2\,d\mu(y)d\mu(z)\\
\begin{split}
  \geq\;&
C^{-1}\frac{\bigl(\sin\measuredangle(L_j,L_k)-C\ve_0^{1/ 
6}\bigr)^2}{\ell(\QS_j)^2}\,
\mu(Q^1_j\cap V_j)\,\mu(Q^1_k\cap V_k)\\
\geq\;&
C(A,\delta)^{-1}\frac{\bigl(\sin\measuredangle(L_j,L_k)
-C\ve_0^{1/6}\bigr)^2}{\ve_0^{-1/25}\ell(\wh{S}_j)^2}\,
\mu(\wh{S}_j)\,\mu(\wh{S}_k)\\\geq\;&
C(A,\delta)^{-1}\frac{\bigl
(\sin\measuredangle(L_j,L_k)-C\ve_0^{1/6}\bigr)^2}{\ve_0^{-1/ 
25}\ell(\wh{S}_j)}\,
\theta_\mu(R)^2\ell(\wh{S}_k).
\end{split}
\end{multline*}
On the other hand, it is easily seen that
  $$c^2_{tr,\mu}(x,2R,2R) \geq \iint_{\begin{subarray}{l}y\in Q^1_j\cap  
V_j\\z\in Q^1_k\cap
V_k\end{subarray}} c(x,y,z)^2\,d\mu(y)d\mu(z).$$
  Since $x\in W$, we have
  $c^2_{tr,\mu}(x,2R,2R)\leq\ve_0\theta_\mu(R)^2$, and then we
  get
$$\ell(\wh{S}_k)\bigl(\sin\measuredangle(L_j,L_k)-C_{24}\ve_0^{1/ 
6}\bigr)^2
\leq C(A,\delta) \ve_0^{24/25}\ell(\wh{S}_j).$$ So we deduce that
either $\sin\measuredangle(L_j,L_k)\leq 2C_{24}\ve_0^{1/6}$, or
otherwise,
  $$\ell(\wh{S}_k) \leq  
C(A,\delta)\ve_0^{24/25}\ve_0^{-1/3}\ell(\wh{S}_j).$$
  Thus,
$$\ell(S_k) \leq C(A,\delta)\ve_0^{47/75}\ell(\wh{S}_j)\leq
\ve_0^{3/5}\ell(S_j),$$ assuming $\ve_0$ small enough.

It remains to show that, in any case, $L_k$ lies in a thin strip
with the same  axis as $V_j$. Remember that $x\in L_k\cap V_j$. If
$\ell(S_k)\leq\ve_0^{3/5}\ell(S_j)$, then $S_k$ (and thus $L_k$)
is contained in a strip with the same  axis as $V_j$ and width
$$\ve_0^{1/6}\ell(S_j)+2\ve_0^{3/5}\ell(S_j)\leq
\ve_0^{1/8}\ell(S_j)$$
  (for $\ve_0$ small).

  Supposing  now that $\sin\measuredangle(L_j,L_k)\leq C\ve_0^{1/6}$, we
  have
$$\ell(S_k)\leq \ell(\wh{S}_k) = \ve_0^{1/50}\ell(\QS_k) \leq
\ve_0^{1/50}\ell(\QS_j) = \ell(\wh{S}_j) \leq
C(A,\delta)\ell(S_j).
$$
We deduce that $L_k$ is also contained in a strip with the same
  axis as $V_j$ and width
  $$\ve_0^{1/6} \ell(S_j) + 2\ell(S_k)\sin\measuredangle(L_j,L_k)
\leq \ve_0^{1/6} \ell(S_j) + C(A,\delta)\ve_0^{1/6}\ell(S_j) \leq
\ve_0^{1/8}\ell(S_j),$$ for $\ve_0$ sufficiently small again.
\hfill\qed

\begin{lemma} \label{lemcurfina1}
Given $j\in J_0${\rm ,} let $x\in L_j${\rm ,} $y,z\not \in S_j${\rm ,} $x_1\in
Q^1_j\cap V_j${\rm ,} and $x_2\in Q_j^2\cap V_j$. Then{\rm ,}
  $$c(x,y,z)\leq C(A,\delta)\,\bigl[c(x_1,y,z) + c(x_2,y,z)\bigr] +
  \frac{C\ve_0^{1/6}\ell(S_j)}{|x-y||x-z|}.$$
\end{lemma}

\Proof 
Let $x'$ be the orthogonal projection of $x$ onto the line
$L_{x_1,x_2}$. Since $x_1,x_2\in V_j$ and $$|x_1-x_2|\geq
\ell(\wh{S}_j)/40\gg \mbox{width of $V_j$},$$ the segment
$L_{x_1,x_2}\cap S_j$ is contained in $CV_j$, where $CV_j$ stands
for the strip with the same  axis as $V_j$ and width $C$ times the
one of $V_j$. Remember also that $L_j$ is a segment supported on
the  axis of $V_j$. As a consequence,
  $$|x-x'| = \dist(x,L_{x_1,x_2}) \leq C\mbox{ width of
$V_j$} \leq C\ve_0^{1/6}\ell(S_j).$$

By Lemma \ref{curpert},
\begin{equation}\label{ccc2}
c(x,y,z) \leq  c(x',y,z) + \frac{C|x-x'|}{|x-y||x-z|},
\end{equation}
because $x,x'$ are in $S_j$ and far from $\partial S_j$, while
$y,z\not\in S_j$.

   It can be shown that there exists some absolute constant $C$ such
  that
  $$\dist(x,L_{y,z}) \leq C\bigl(\dist(x_1,L_{y,z}) +
  \dist(x_2,L_{y,z})\bigr).$$
  This follows easily from the fact that $x',x_1,x_2$ are collinear
  and $|x_1-x_2|\geq C^{-1}|x'-x_1|$.
  Notice also that, for $i=1,2$,
  $$|x_i-y|\leq C(A,\delta)|x'-y|\qquad \mbox{and}\qquad |x_i-z|\leq
  C(A,\delta)|x'-z|.$$ In fact, the constants $C(A,\delta)$ above
  depend on the ratio $\ell(\wh{S}_j)/\ell(S_j)$.
  We get
  \begin{eqnarray*}
  c(x',y,z)\! &= &\!\frac{2\dist(x,L_{y,z})}{|x-y||x-z|} \leq
  \frac{C(A,\delta)\dist(x_1,L_{y,z})}{|x_1-y||x_1-z|} +
\frac{C(A,\delta)\dist(x_2,L_{y,z})}{|x_2-y||x_2-z|}\\
  & = & C(A,\delta)\,\bigl[c(x_1,y,z) + c(x_2,y,z)\bigr].
  \end{eqnarray*}
  From this estimate, \rf{ccc2}, and the fact that $|x-x'|\lesssim
  \ve_0^{1/6}\ell(S_j)$, the lemma follows.
\hfill\qed

\begin{lemma} \label{maximaltr}
Let $M_{tr,\mu}$ be the following \/{\rm (}\/truncated\/{\rm )}\/ maximal operator
  $$M_{tr,\mu}f(x) = \sup_{\begin{subarray}{l}
Q:x\in\frac12Q\\ \ell(Q)>8\ell(S_x)
  \end{subarray}}
  \frac{1}{\mu(Q)}\int_{Q\cap2R}|f|\,d\mu,\qquad \mbox{for $
x\in \bigcup_{j\in J_0} S_j,$}$$
  where $S_x$ is the square $S_j${\rm ,}
$j\in J_0${\rm ,} which contains $x$. Then{\rm ,} $M_{tr,\mu}$ is bounded from
$L^2(\mu)$ into $L^2(\nu)${\rm ,} with norm depending on $A$ and
$\delta$.
\end{lemma}

Notice that the notation ``$S_x$'' was also used at the
beginning of Subsection \ref{subsub1}, but with a different
meaning.

\Proof 
We immediately  check that $M_{tr,\mu}$ is bounded from
$L^\infty(\mu)$ into $L^\infty(\nu)$. So, by interpolation it is
enough to show that it is also bounded from $L^1(\mu)$ into
$L^{1,\infty}(\nu)$. Take a fixed $\lambda>0$. If
$M_{tr,\mu}f(x)>\lambda$ for some $x\in\supp(\nu)$, there is a
square $Q_x$ such that $x\in\frac12Q_x$, $\ell(Q_x)>8\ell(S_x)$,
and $\int_{Q_x}|f|d\mu/\mu(Q_x)>\lambda$. By the Besicovitch covering
theorem, there exists a family of points
$\{x_i\}_i\subset\supp(\nu)$ so that the family of squares
$\{Q_{x_i}\}_i$ has finite overlap and
$\{x:M_{tr,\mu}f(x)>\lambda\}\subset\bigcup_i Q_{x_i}$. Since
$\ell(Q_{x_i})>8\ell(S_{x_i})$ and $x_i\in\frac12Q_{x_i}$, it is
easy to check that there exists a square $P$ concentric with
$S_{x_i}$ and with side length $\ell(P)=\ell(Q_{x_i})/2$ such that
$S_{x_i}\subset P \subset Q_{x_i}$. Then,
  $$\nu(Q_{x_i})\leq CA\theta_\mu(R)\ell(Q_{x_i})
= CA\theta_\mu(R)\ell(P) \leq CA\delta^{-1}\mu(P) \leq
C(A,\delta)\mu(Q_{x_i}).$$ Thus,
\begin{eqnarray*}
\nu\{x:M_{tr,\mu}f(x)>\lambda\}& \leq &\sum_i\nu(Q_{x_i}) \leq
C(A,\delta)\sum_i\mu(Q_{x_i}) \\*& \leq&
\frac{C(A,\delta)}{\lambda}\sum_i \int_{Q_{x_i}}|f|\,d\mu \leq
\frac{C(A,\delta)}{\lambda} \int|f|\,d\mu.
\end{eqnarray*}
\vglue-28pt
\Endproof\vskip4pt

\Subsubsec{Proof of \rf{curnu}}
As in the preceding lemma, for $x\in \bigcup_{j\in J_0} S_j$, we
denote by $S_x$ be the square $S_j$, $j\in J_0$, which contains
$x$. Analogously, $\wh{S}_x$, $\QS_x$, $Q_x^1$, $Q^2_x$, and $V_x$
stand for $\wh{S}_j$, $\QS_j$, $Q_j^1$, $Q^2_j$, and $V_j$
respectively.

We denote
\begin{eqnarray*}
F_1  &:= & \bigl\{(x,y,z)\in (\textstyle \bigcup_{j\in J_0}
S_j)^3:\,S_x=S_y \neq
S_z\bigr\},\\
F_2& := & \bigl\{(x,y,z)\in (\textstyle\bigcup_{j\in J_0} S_j)^3:\,
S_x\neq S_y \neq
S_z\neq S_x\bigr\},\\
F_3& := & \bigl\{(x,y,z)\in F_2:\,
3\QS_y\cap3\QS_z=\varnothing\bigr\},\\
F_4& := & \bigl\{(x,y,z)\in F_2:\,
3\QS_x\cap3\QS_y\neq\varnothing,\,
3\QS_x\cap3\QS_z\neq\varnothing,\,3\QS_y\cap3\QS_z\neq\varnothing\bigr\} 
.
\end{eqnarray*}

Since $c^2(\nu_{|S_j})=0$ for all $j\in J_0$,
\begin{equation}\label{spl1}
\begin{split}
c^2(\nu) & \;= \;
  \iiint_{\bigl(\bigcup_{j\in J_0}
  S_j\bigr)^3}c(x,y,z)^2\,d\nu(x)d\nu(y)d\nu(z)\\ &\;=\;
  3\iiint_{F_1}\cdots
  +\iiint_{F_2}\cdots
   \;\leq \;
  3\iiint_{F_1}\cdots + 3\iiint_{F_3}\cdots +\iiint_{F_4}\cdots\\
   & \;=: \; 3I_1 + 3I_3 +I_4.
\end{split}
\end{equation}
 
\vglue12pt $\bullet$ {\it Estimates for} $I_3$.  If
$y'\in S_y\cap W$ and $z'\in S_z\cap W$, by Lemma \ref{curpert}
  \begin{eqnarray}\label{ldf10} 
c(x,y,z)&\leq& c(x,y',z')+ \frac{C\ell(S_y)}{|y-x||y-z|} +
\frac{C\ell(S_z)}{|z-x||z-y|} \\
  &=:&  c(x,y',z')+ C\,\bigl[T_y(x,y,z) +
T_z(x,y,z)\bigr].\nonumber
\end{eqnarray}
Then it easily follows that
\begin{eqnarray}\label{ldf11}\qquad
I_3 &\leq&\; 2\iiint_{\begin{subarray}{l}(x,y,z)\in F_2\\y,z\in W\\
3\QS_y\cap 3\QS_z=\varnothing
\end{subarray}}
c(x,y,z)^2\,d\nu(x)d\mu(y)d\mu(z) \\ && + C
\iiint_{\begin{subarray}{l}(x,y,z)\in F_2\\
3\QS_y\cap 3\QS_z=\varnothing
\end{subarray}}
\bigl[T_y(x,y,z)^2 + T_z(x,y,z)^2\bigr]\,d\nu(x)d\nu(y)d\nu(z)\nn\\
&=: & 2 I_{3,1} + C\,I_{3,2}.
\nn
\end{eqnarray}
Although it is not written explicitly, all the integrals above
are restricted to $(2R)^3$ (and the same for the rest of the proof
of \rf{curnu}).

First, dealing  with the term $I_{3,2}$, we have
\begin{multline*}
\iiint_{\begin{subarray}{l}(x,y,z)\in F_2\\
3\QS_y\cap 3\QS_z=\varnothing
\end{subarray}}T_y(x,y,z)^2\,d\nu(x)d\nu(y)d\nu(z)
\;\leq\; \iiint_{\begin{subarray}{l}|y-x|>\ell(S_y)/2\\
|y-z|>\ell(\QS_y)
\end{subarray}}\cdots\\
\begin{split}
  \leq &\;
\int \biggl(\int_{|y-x|>\ell(S_y)/2}\frac{\ell(S_y)}{|y-x|^2}
\,d\nu(x)\biggr)\biggl(\int_{|y-z|>\ell(\QS_y)}\frac{\ell(S_y)}{|y-z|^2}
\,d\nu(z)\biggr)d\nu(y)\\
  \leq &\;
C\Bigl(A\theta_\mu(R)\Bigr)\,\biggl(A\theta_\mu(R)\frac{\ell(S_y)}{\ell( 
\QS_y)}\biggr)\,\nu(\C)
\leq CA^2\theta_\mu(R)^2\,\ve_0^{1/50}\,\mu(R).
\end{split}
\end{multline*}
We have analogous estimates for the integral of $T_z(\cdots)^2$.
Thus,
  $$I_{3,2}\leq
CA^2\theta_\mu(R)^2\,\ve_0^{1/50}\,\mu(R).$$

Now we consider the term $I_{3,1}$ in \rf{ldf11}. By Lemma
\ref{lemcurfina1}, for all $x_1\in Q^1_x\cap V_x$ and $x_2\in
Q_x^2\cap V_x$ we have
  $$c(x,y,z)\leq C(A,\delta)\,\bigl[c(x_1,y,z) + c(x_2,y,z)\bigr] +
  \frac{C\ve_0^{1/6}\ell(S_x)}{|x-y||x-z|}.$$
Integrating over $x_1\in Q_x^1\cap V_x$ and over $x_2\in Q_x^2\cap
V_x$ with respect to $\mu$, we obtain
\begin{eqnarray*}
c(x,y,z)&\leq
&\frac{C(A,\delta)}{\mu(\wh{S}_x)}\biggl(\int_{x_1\in V_x\cap
Q_x^1}
  \!\!\!\! c(x_1,y,z)d\mu(x_1)
+ \int_{x_2\in V_x\cap Q_x^2}\!\!\!\! c(x_2,y,z)d\mu(x_2)\biggr)\\
& & +C\ve_0^{1/6} T_x(x,y,z),
\end{eqnarray*}
where $T_x(x,y,z):=\ell(S_x)/\bigl(|x-y||x-z|\bigr)$.
  Therefore,
\begin{eqnarray*}
c(x,y,z)&\leq& \frac{C(A,\delta)}{\mu(\wh{S}_x)} \int_{\wh{S}_x}
c(w,y,z)\,d\mu(w) +C\ve_0^{1/6} T_x(x,y,z)\\ &\leq&
C(A,\delta)\,M_{tr,\mu}\bigl[c(\cdot,y,z)\bigr](x)+C\ve_0^{1/6}
T_x(x,y,z).
\end{eqnarray*}
Thus,
\begin{eqnarray*}
I_{3,1} &\leq& C(A,\delta)\iiint_{3\QS_y\cap 3\QS_z=\varnothing}
M_{tr,\mu}\bigl[c(\cdot,y,z)\bigr](x)^2
\,d\nu(x)d\mu(y)d\mu(z) \\
&&\mbox{} + C\ve_0^{1/3}\iiint_{|x-y|,|x-z|\geq
\ell(S_x)/2}T_x(x,y,z)^2\,d\nu(x)d\mu(y)d\mu(z).
\end{eqnarray*}
The last integral on the right side is estimated as follows:
\begin{multline*}
\iiint_{|x-y|,|x-z|\geq
\ell(S_x)/2}T_x(x,y,z)^2\,d\nu(x)d\mu(y)d\mu(z) \\ \leq
\int\biggl(
\int_{|x-y|>\ell(S_x)/2}\frac1{|x-y|^2}\,d\mu(y)\biggr)^2d\nu(x)\leq
CA^2\theta_\mu(R)^2\mu(R).
\end{multline*}

On the other hand, by Lemma \ref{maximaltr} we know that
$M_{tr,\mu}$ is bounded from $L^2(\mu)$ into $L^2(\nu)$. So we
have
  \begin{eqnarray*}
I_{3,1} &\leq &C(A,\delta)\iiint_{3\QS_y\cap
3\QS_z=\varnothing} c(x,y,z)^2\,d\mu(x)d\mu(y)d\mu(z) \\
\mbox{} &&+ C\ve_0^{1/3} A^2\theta_\mu(R)^2\mu(R).
  \end{eqnarray*}
It is easy to check that
  \begin{eqnarray*}
\iiint_{3\QS_y\cap 3\QS_z=\varnothing}
c(x,y,z)^2\,d\mu(x)d\mu(y)d\mu(z)&\leq &\int
c^2_{tr,\mu}(y,2R,2R)\,d\mu(y) \\ &\leq&
C\ve_0\theta_\mu(R)^2\mu(2R).
  \end{eqnarray*}

 From the preceding estimates for $I_{3,1}$ and $I_{3,2}$, we get
$$I_3 \leq C(A,\delta)\ve_0^{1/50}\theta_\mu(R)^2\mu(R).$$

\demo{$\bullet$ Estimates for $I_4$}  By
Fubini, we have
  $$I_4 \leq 3\iiint_{\begin{subarray}{l}(x,y,z)\in
  F_4\\ \ell(\QS_x)\geq \ell(\QS_y)\geq\ell(\QS_z)\end{subarray}}
c(x,y,z)^2\,d\nu(x)d\nu(y)d\nu(z) =:3I_4'.$$ Now we split $I_4'$
as follows:
\begin{equation*}
\begin{split}
I_4' \,&= \; \biggl(\iiint_{\begin{subarray}{l}(x,y,z)\in
  F_4\\ \ell(\QS_x)\geq \ell(\QS_y)\geq\ell(\QS_z)\\
  \ell(S_x)\geq \ell(\QS_y)\end{subarray}} \!\!\!
+ \iiint_{\begin{subarray}{l}(x,y,z)\in
  F_4\\ \ell(\QS_x)\geq \ell(\QS_y)\geq\ell(\QS_z)\\
  \ell(S_x)< \ell(\QS_y)\end{subarray}}\biggr)
c(x,y,z)^2\,d\nu(x)d\nu(y)d\nu(z)\\
  &=: \; I_{4,1} + I_{4,2}.
\end{split}
\end{equation*}

First we will study $I_{4,1}$. For $(x,y,z)$ in the domain of
integration of $I_{4,1}$ we have
  $$|x-y| \geq \frac{\ell(S_x)}4\geq \frac{\ell(\QS_y)}4\qquad\mbox{and}
\qquad |x-z| \geq \frac{\ell(S_x)}4\geq \frac{\ell(\QS_y)}4 \geq
\frac{\ell(\QS_z)}4.$$
  The estimates for $I_{4,1}$ are similar to the ones
  for $I_3$. Indeed, consider $y'\in S_y\cap W$ and
  $z'\in S_z\cap W$, so that \rf{ldf10} also holds in this case.
  Instead of \rf{ldf11} now we get
\begin{equation}
\pagebreak
\label{ldf12}
\begin{split}
I_{4,1} \leq&\; 2\iiint_{\begin{subarray}{l}(x,y,z)\in F_4\\y,z\in W\\
\end{subarray}}
c(x,y,z)^2\,d\nu(x)d\mu(y)d\mu(z)\\ &\mbox{} + C
\iiint_{\begin{subarray}{l}(x,y,z)\in F_4\\
|x-y|\geq\ell(\QS_y)/4\\ |x-z| \geq\ell(\QS_z)/4
\end{subarray}}
\bigl[T_y(x,y,z)^2 + T_z(x,y,z)^2\bigr]\,d\nu(x)d\mu(y)d\mu(z).
\end{split}
\end{equation}
We have
\begin{multline*}
\iiint_{\begin{subarray}{l}(x,y,z)\in F_4\\
|x-y|\geq\ell(\QS_y)/4\\ |x-z| \geq\ell(\QS_z)/4
\end{subarray}}
T_y(x,y,z)^2 \,d\nu(x)d\mu(y)d\mu(z)
\;\leq\; \iiint_{\begin{subarray}{l}|y-x|>\ell(\QS_y)/4\\
|y-z|>\ell(S_y)/4
\end{subarray}}\cdots\\
\begin{split}
  \leq &\;
\int \biggl(\int_{|y-x|>\ell(S_y)/4}\frac{\ell(S_y)}{|y-x|^2}
\,d\nu(x)\biggr)\biggl(\int_{|y-z|>\ell(\QS_y)/4}\frac{\ell(S_y)}{|y- 
z|^2}
\,d\mu(z)\biggr)d\mu(y)\\
  \leq &\;
C\Bigl(A\theta_\mu(R)\Bigr)\,\biggl(A\theta_\mu(R)\frac{\ell(S_y)}{\ell( 
\QS_y)}\biggr)\,\mu(R)
\leq CA^2\theta_\mu(R)^2\,\ve_0^{1/50}\,\mu(R).
\end{split}
\end{multline*}
Analogous estimates hold for the integral of $T_z(\cdots)^2$
since in the domain of integration we have
$|x-z|\geq\ell(\QS_z)/4$.

To estimate the integral $$\iiint_{\begin{subarray}{l}(x,y,z)\in  
F_4\\y,z\in W\\
\end{subarray}}
c(x,y,z)^2\,d\nu(x)d\mu(y)d\mu(z)$$ in \rf{ldf12}, the same
arguments used for $I_{3,1}$ work in this case. Only some minor
changes which are left for the reader are required.

Now we deal with $I_{4,2}$. Take $(x,y,z)$ in the domain of
integration of $I_{4,2}$. Since $\ell(\QS_x)\geq
\ell(\QS_y)\geq\ell(\QS_z)$, $3\QS_x\cap 3\QS_y\neq\varnothing$,
and $3\QS_x\cap 3\QS_z\neq\varnothing$, we have
$3\QS_y,3\QS_z\subset 9\QS_x$. Remember that, by Lemma
\ref{lemangles}, the segments $L_y,L_z$ (and thus $y$ and $z$) are
contained in a strip with the same  axis as $V_x$ and width
$\ve_0^{1/8}\ell(S_x)$. Therefore,
$$\sin\measuredangle(L_x,\stackrel{\leftrightarrow}{xy})\leq
\frac{\ve_0^{1/8}\ell(S_x)}{|x-y|}\leq
\frac{\ve_0^{1/8}\ell(S_x)}{C^{-1}\ell(S_x)} \leq C\ve_0^{1/8},$$
and, in the same way,
$\sin\measuredangle(L_x,\stackrel{\leftrightarrow}{xz})\leq
C\ve_0^{1/8}$. Thus, we get
  $$\sin\measuredangle(y,x,z) \leq C
  \bigl(\sin\measuredangle(L_x,\stackrel{\leftrightarrow}{xy}) +
  \sin\measuredangle(L_x,\stackrel{\leftrightarrow}{xz})\bigr)
  \leq C\ve_0^{1/8}.$$
Therefore, since $\ell(S_x)<\ell(\QS_y)$,
  $$c(x,y,z) = \frac{2\sin\measuredangle(y,x,z)}{|y-z|} \leq
  \frac{C\ve_0^{1/8}}{\ell(S_y)}\leq
  \frac{C(A,\delta)\,\ve_0^{1/8}}{\ve_0^{1/50}\ell(\QS_y)}
  \leq \frac{C(A,\delta)\,\ve_0^{21/200}}{\ell(S_x)},$$
and so
  \begin{eqnarray*}
I_{4,2}& \leq &\int_{x\in 2R}\biggl(\iint_{y,z\in 9\QS_x}
   \frac{C(A,\delta)\,\ve_0^{21/ 
100}}{\ell(S_x)^2}\,d\mu(y)d\mu(z)\biggr)d\mu(x)\\
&\leq& C(A,\delta)\ve_0^{21/100}\ve^{-1/25}\theta_\mu(R)^2 \mu(R) =
C(A,\delta)\ve_0^{17/100}\theta_\mu(R)^2 \mu(R).
  \end{eqnarray*}

Gathering the estimates for $I_{4,1}$ and $I_{4,2}$, we obtain
$$I_4 \leq
C(A,\delta)\,(\ve_0^{1/50}+\ve_0^{17/100})\,\theta_\mu(R)^2 \mu(R)
\leq C(A,\delta)\ve_0^{1/50}\theta_\mu(R)^2 \mu(R).$$

\demo{$\bullet$ Estimates for $I_1$}
We have
\begin{equation*}
\begin{split}
I_1 =\;& \iiint_{\begin{subarray}{l} x\in\supp(\nu)\\ y\in
  L_x\\z\not\in\QS_x
  \end{subarray}} c(x,y,z)^2\,d\nu(x)d\nu(y)d\nu(z) \\ &\mbox{} +
\iiint_{\begin{subarray}{l} x\in\supp(\nu)\\ y\in
  L_x\\z\in\QS_x\setminus S_x
  \end{subarray}} c(x,y,z)^2\,d\nu(x)d\nu(y)d\nu(z) \;=:\; I_{1,1} +
  I_{1,2}.
\end{split}
\end{equation*}
The term $I_{1,1}$ is estimated as follows:
\begin{eqnarray*}
I_{1,1} & \leq & \iiint_{\begin{subarray}{l} x\in\supp(\nu)\\ y\in
  L_x\\ |x-z|>\ell(\QS_x)/4
  \end{subarray}} \frac1{|x-z|^2}\,d\nu(x)d\nu(y)d\nu(z)\\
  & \leq &
   \int_{x\in\supp(\nu)}\frac{CA\theta_\mu(R)\nu(L_x)}{\ell(\QS_x)}\,d\nu(x 
)
  \;\leq\;  CA\theta_\mu(R)^2\ve_0^{1/50}\mu(R),
\end{eqnarray*}
since $\nu(L_x)/\ell(\QS_x)\leq
C\theta_\mu(R)\ell(S_x)/\ell(\QS_x) \leq C\theta_\mu(R)\ve^{1/50}$
and $\nu(\C)\leq C\mu(R)$.

Finally we turn our attention to $I_{1,2}$. Consider $(x,y,z)$ in
the domain of integration of $I_{1,2}$. Clearly, in this case we
have $\QS_x\cap\QS_z\neq\varnothing$. If $\ell(\QS_z)\leq
\ell(\QS_x)$, by Lemma \ref{lemangles}, $z$ is contained in a
strip with the same  axis as $V_x$ and width
$\ve_0^{1/8}\ell(S_x)$.

Suppose now that $\ell(\QS_z)> \ell(\QS_x)$. Then, again by Lemma
\ref{lemangles}, either $\sin\measuredangle(L_x,L_z)\leq
C\ve_0^{1/6}$ or $\ell(S_x)\leq\ve_0^{3/5}\ell(S_z)$. However, the
latter inequality cannot hold because it implies
  $$|x-z| \geq\ell(S_z)/4 \geq \ve_0^{-3/5}\ell(S_x)/4 \gg
  \ell(\QS_x),$$
  and so $z\not\in\QS_x$. Then the condition  
$\sin\measuredangle(L_x,L_z)\leq
C\ve_0^{1/6}$ holds. Further, we have $\ell(S_z)\leq 2\ell(\QS_x)$
because $z\in\QS_x$. As a consequence, we easily infer that $z$
lies in a thin strip with the same  axis as $V_x$ and width
$$\ve_0^{1/6}\ell(S_x) + C\ve_0^{1/6}\ell(\QS_x) \leq
C(A,\delta)\ve_0^{1/6-1/50}\ell(S_x) \leq
\ve_0^{2/15}\ell(S_x),$$ for $\ve_0$ small enough.

So in any case $z$ is contained in the strip with the same  axis
as $V_x$ and width $\ve_0^{2/15}\ell(S_x)$. As a consequence, we
deduce
  $$\sin\measuredangle(x,y,z)\lesssim
  \frac{\ve_0^{2/15}\ell(S_x)}{|y-z|}\lesssim \ve_0^{2/15}.$$
  Thus, $c(x,y,z)\lesssim \ve_0^{2/15}/|x-z|\lesssim
  \ve_0^{2/15}/\ell(S_x)$, and so
  \begin{eqnarray*}
I_{1,2} &\lesssim &
\int_{x\in\supp(\nu)}\frac{A\ve_0^{4/ 
15}}{\ell(S_x)^2}\,\theta_\mu(R)^2\ell(\QS_x)\ell(S_x)
\,d\nu(x)\\*
& \leq & C(A,\delta)\ve_0^{4/15-1/50}\,\theta_\mu(R)^2\mu(R) \;=\;
C(A,\delta)\ve_0^{37/150}\,\theta_\mu(R)^2\mu(R).
  \end{eqnarray*}

\demo{$\bullet$ End of the proof}
By the estimates obtained for $I_1$, $I_3$ and $I_4$,
we get 
\vglue12pt
\hfill $\displaystyle{c^2(\nu) \leq
C(A,\delta)\theta_\mu(R)^2\ve_0^{1/50}\mu(R).}$  \Endproof

\section{The curvature of $\vphi_\sharp\mu$} \label{seccurv}

In this section we denote $\sigma:=\vphi_\sharp\mu$ and $F:=\vphi(E)$.

Given a square $Q$, we say that $\vphi(Q)$ is a $\vphi$-square. If
$x_Q$ is the center of $Q$, then we call $\vphi(x_Q)$ the center
of $\vphi(Q)$. We also set $\ell(\vphi(Q)) := \ell(Q)$. Since
$\vphi$ is bilipschitz, we have $\ell(\vphi(Q)) \simeq
\diam(\vphi(Q))$. We will often use the letters $P,Q,R$ to denote
$\vphi$-squares too. If $Q$ is a dyadic (or $4$-dyadic) square, we
say that $\vphi(Q)$ is a dyadic (or $4$-dyadic) $\vphi$-square.

If $Q= \vphi(Q_0)$ is a $\vphi$-square, we let  $\lambda Q =
\vphi(\lambda Q_0)$, for $\lambda>0$. Then  $Q$ is
$\lambda$-doubling if $\sigma(\lambda Q)\leq C\sigma(Q)$ for some
$C\geq1$. We also set $$\theta_\sigma(Q) : =
\frac{2^{1/2}\sigma(Q)}{\diam(Q)}$$ (the number $2^{1/2}$ is due
to aesthetic reasons; if $\vphi$ is the identity, then the
definition coincides with \rf{eqdens}) and if $R$ is another
$\vphi$-square which contains $Q$, we put
$$\delta_\sigma(Q,R) := \int_{R_Q\setminus Q}
\frac1{|y-x_{Q}|}\,d\sigma(y),$$ where $x_{Q}$ stands for the
center of $Q$  and $R_Q$ is the smallest $\vphi$-square concentric
with $Q$ that contains $R$.

Given a family $\ttt(F)$ of $4$-dyadic $\vphi$-squares and a fixed
$Q\in \ttt(F)$, we denote by $\sss(Q)$ the subfamily of
$\vphi$-squares which satisfy the properties (a), (b), (c) stated
at the beginning of Section \ref{seccorona} (with squares
replaced by $\vphi$-squares). The set $G(Q)$ is also defined as in
Section \ref{seccorona}, with $\vphi$-squares instead of squares.

\begin{mlemma} \label{lemafi}
Let $\sigma$ be a Radon measure supported on a compact
$F\subset\C$. Suppose that $\sigma(B(x,r))\leq C_0r$ for all
$x\in\C,\,r>0$. Let $\ttt(F)$ be a family of $4$-dyadic
$16$-doubling $\vphi$-squares \/{\rm (}\/called top $\vphi$-squares\/{\rm )}\/ which
contains some $4$-dyadic $\vphi$-square $R_0$ such that $F\subset
R_0${\rm ,} and such that for each $Q\in \ttt(F)$ there exists a
$C_{25}$-{\rm AD} regular curve $\Gamma_Q$ satisfying\/{\rm :}\/
\begin{itemize}
\item[{\rm (a)}] $\sigma$-almost every point in $G(Q)$ belongs to
$\Gamma_Q$.

\item[{\rm (b)}] For each $P\in \sss(Q)$ there exists some $\vphi$-square
$\wt{P}$ containing $P$ such that $\delta_\sigma(P,\wt{P})\leq
C\theta_\sigma(Q)$ and $\wt{P}\cap \Gamma_Q\neq \varnothing$.

\item[{\rm (c)}] If $P$ is a $\vphi$-square with $\ell(P)\leq \ell(Q)$ such  
that either
$P\cap G(Q)\neq\varnothing$ or there is another $\vphi$-square
$P'\in\sss(Q)$ such that $P\cap P'\neq\varnothing$ and
$\ell(P')\leq\ell(P)${\rm ,} then $\sigma(P)\leq
C\,\theta_\sigma(Q)\,\ell(P).$
\end{itemize}
Then{\rm ,}
$$c^2(\sigma) \leq C \sum_{Q\in \ttt(F)} \theta_\sigma(Q)^2
\sigma(Q).$$
\end{mlemma}
\vskip8pt

We will prove this lemma in Subsections
\ref{subdecomp}--\ref{subalt}.

\Subsec{Proof of Theorem {\rm \ref{teocurv}}}
This is an easy consequence of Main Lemmas \ref{lemcorona} and
\ref{lemafi}. Indeed, if $c^2(\mu)<\infty$, then we have the
corona decomposition given by Main Lemma \ref{lemcorona}. Applying
the bilipschitz map $\vphi$, we obtain another corona decomposition
for $F=\vphi(E)$ like the one required in Main Lemma \ref{lemafi}.
In particular, notice that $\vphi$ sends {\rm AD} regular curves to {\rm AD}
regular curves, and also if $Q,R$ are squares such that $Q\subset
R$, then
\begin{eqnarray*}
\delta_\sigma(\vphi(Q),\,\vphi(R)) & = &
\int_{\vphi(R_Q)\setminus \vphi(Q)}
\frac1{|y-x_{\vphi(Q)}|}\,d\sigma(y) \\ &=& \int_{R_Q\setminus
Q}\frac1{|\vphi(y)-\vphi(x_Q)|}\,d\mu(y)\\
& \leq&
 \int_{R_Q\setminus Q}\frac{C}{|y-x_Q|}\,d\mu(y) =
C\delta_\mu(Q,R),
\end{eqnarray*}
with $C$ depending on $\vphi$. So, by Main Lemma \ref{lemafi},
$c^2(\sigma)\lesssim\bigl(\mu(E)+c^2(\mu)\bigr)$.\fiproof
 
\Subsec{Decomposition of $c^2(\sigma)$} \label{subdecomp}
We start the proof of Main Lemma \ref{lemafi}. Observe that
$$c^2(\sigma) \leq 3 \iiint_{|x-y|\geq |x-z|,|y-z|}
c(x,y,z)^2\,d\sigma(x) d\sigma(y)d\sigma(z).$$ We now introduce a
variant of the curvature operator $K_\sigma$. Consider the kernel
$$\wh{k}_\sigma(x,y) = \int_{z:|x-y|\geq |x-z|,|y-z|}
c(x,y,z)^2\,d\sigma(z),$$
  and set
  $$\wh{K}_\sigma f(x) = \int
\wh{k}_\sigma(x,y)\,f(y)\,d\sigma(y),\hspace{8mm} x\in\C,f\in
L^1_{\rm loc}(\sigma)$$ (compare with the definition of $k_\mu(x,y)$
and $K_\mu$ in Section \ref{secprelim}).  We have $$\int_F
\wh{K}_\sigma\chi_F(x)\,d\sigma(x) \leq c^2(\sigma) \leq 3 \int_F
\wh{K}_\sigma\chi_F(x)\,d\sigma(x).$$
  The truncated operator
$\wh{K}_{\sigma,j}$, $j\in\Z$, is
  $$\wh{K}_{\sigma,j} f(x) =
\int_{|x-y|>L^{-1}2^{-j}}
\wh{k}_\sigma(x,y)\,f(y)\,d\sigma(y),\hspace{8mm} x\in\C,f\in
L^1_{\rm loc}(\sigma),$$
  where $L$ is the bilipschitz constant of $\vphi$.
  We say that a $\vphi$-square
$Q\in\ttt(F)$ is a descendant of another $\vphi$-square
$R\in\ttt(F)$ if there is a chain $R=Q_1,\,Q_2,\dots,Q_n=Q$, with
$Q_i\in\ttt(F)$ such that $Q_{i+1}\in\sss(Q_i)$ for each $i$. Only
the $\vphi$-squares from $\ttt(F)$ which are descendants of $R_0$
will be relevant to estimate $c^2(\sigma)$. So we assume that all
the $\vphi$-squares in $\ttt(F)$ are of this type.

To decompose $c^2(\sigma)$, we prefer to use dyadic
$\vphi$-squares instead of $4$-dyadic $\vphi$-squares. A
$\vphi$-square $Q$ belongs to the family $\ttt_{\rm dy}(F)$ if
there exists some $R\in \ttt(F)$ such that $Q$ is one of the $16$
dyadic $\vphi$-squares contained in $R$ with side length
$\ell(R)/4$.

Note  that if $Q\in \ttt_{\rm dy}(F)$, then $Q$ is
contained in a $4$-dyadic $\vphi$-square $R$ such that $Q\subset
R\subset 7Q\subset 3R$. Moreover, since each $4$-dyadic
$\vphi$-square $R\in \ttt(F)$ is made up of $16$ dyadic
$\vphi$-squares $Q\in\ttt_{\rm dy}(F)$, we get (using the doubling
properties of the $\vphi$-squares in $\ttt(F)$)
\begin{equation} \label{erw1}
\sum_{Q\in \ttt_{\rm dy}(F)} \theta_\sigma(7Q)^2 \sigma(7Q)
\lesssim \sum_{R\in \ttt(F)} \theta_\sigma(R)^2 \sigma(R).
\end{equation}

Given $Q\in \ttt(F)$ or $Q\in\ttt_{\rm dy}(F)$, we denote by ${\rm
Term}(Q)$ the family of {\em maximal} dyadic (and thus disjoint)
$\vphi$-squares $P\in \ttt_{\rm dy}(F)$, with $P\subsetneq Q$.
Finally, we let $\tree(Q)$ be the class of dyadic $\vphi$-squares
contained in $Q$, different from $Q$, which are not proper
$\vphi$-subsquares of any $P\in{\rm Term}(Q)$.

We denote by $\vphi\Delta$ the class of dyadic $\vphi$-squares
contained in $R_0$, and by $\vphi\Delta_j$ those $\vphi$-squares
in $\vphi\Delta$ with side length $2^{-j}$. We have $$\vphi\Delta
= \bigl\{Q\in\vphi\Delta:\ell(Q)\geq\ell(R_0)/4\bigr\} \cup
\bigcup_{Q\in \ttt_{\rm dy}(F)} \tree(Q).$$
  Observe that the
$\vphi$-squares $Q\in\vphi\Delta$ such that
$\ell(Q)\geq\ell(R_0)/4$ are the only $\vphi$-squares in
$\ttt_{\rm dy}(F)$ which may not belong to any $\tree(R)$,
$R\in\ttt_{\rm dy}(F)$. Notice also that $\tree(Q) \cap
\tree(R)=\varnothing$ if $Q\neq R$.

We split the curvature $c^2(\sigma)$ as follows:
\begin{eqnarray*}
c^2(\sigma) & \simeq & \sum_{j} \sum_{Q\in \vphi\Delta_j} \int_Q
(\wh{K}_{\sigma,j+1} \chi_F -
\wh{K}_{\sigma,j} \chi_F)\, d\sigma + \int_F \wh{K}_{\sigma,J(R_0)+2}  
\chi_F\,d\sigma\\[4pt]
& = & \sum_{R\in \ttt_{\rm dy}(F)} \sum_{Q\in \tree(R)} \int_Q
(\wh{K}_{\sigma,J(Q)+1} \chi_F - \wh{K}_{\sigma,J(Q)} \chi_F)\, d\sigma  
\\[4pt]
\mbox{}&&+ \int_F \wh{K}_{\sigma,J(R_0)+2} \chi_F\,d\sigma,
\end{eqnarray*}
where $J(Q)$ stands for the integer $j$ such that $Q\in
\vphi\Delta_j$. Since
$$\int_F \wh{K}_{\sigma,J(R_0)+2} \chi_F\,d\sigma \lesssim
\theta_\sigma(R_0)^2\sigma(F),
\pagebreak
$$ to prove Main Lemma \ref{lemafi}
it is enough to show that
\begin{equation} \label{eqclau}
\sum_{Q\in \tree(R)}\int_Q (\wh{K}_{\sigma,J(Q)+1} \chi_F -
\wh{K}_{\sigma,J(Q)} \chi_F)\, d\sigma \lesssim
\theta_\sigma(7R)^2 \sigma(7R),
\end{equation}
for every $R\in \ttd(F)$, by \rf{erw1}.

\Subsec{Regularization of the stopping $\vphi$-squares}  
\label{subregul}
Given a fixed $R\in \ttd(F)$, let $R_1$ be a $4$-dyadic
$\vphi$-square $R_1\in \ttt(F)$ such that $R\subset R_1\subset
7R$ (it does not matter which $R_1$ if it is not unique). Let  
$\Gamma_R:=\Gamma_{R_1}$ be the {\rm AD}
regular curve satisfying (a) and (b) in Main Lemma \ref{lemafi}.

It seems that after defining $\ttt_{\rm dy}(F)$ we should
introduce the family $\sss_{\rm dy}(R)$ analogously. However, for
technical reasons, it is better to introduce a regularized version
of ${\rm Stop_{\rm dy}}(R)$ (it does not matter what ${\rm
Stop_{\rm dy}}(R)$ means precisely), that we will denote by
$\reg(R)$. First we set
  $$d_R(x) := \inf_{Q\in \sss(R_1)}\bigl\{
\dist(x,Q) + \ell(Q),\,\,\dist(x,G(R_1))\bigr\}.$$ For each $x\in
3R\cap\supp(\sigma)\setminus [G(R_1)\cup Z(\sigma)]$ (recall that  
$Z(\sigma)$ is a set of
zero $\sigma$-measure, defined similarly to $Z(\mu)$ at the beginning  
of Section \ref{seccorona}),
let $Q_x$ be a dyadic $\vphi$-square
containing $x$ such that
\begin{equation} \label{eqqx}
\frac{d_R(x)}{20L} < \ell(Q_x) \leq \frac{d_R(x)}{10L}.
\end{equation}
  Remember that $L$ is the bilipschitz constant of
$\vphi$. Then, $\reg(R)$ is a maximal (and thus disjoint)
subfamily of $\{Q_x\}_{x\in 3R\cap\supp(\sigma)\setminus [G(R_1)\cup  
Z(\sigma)]}$.

\begin{lemma} \label{regul}
\begin{itemize}
\item[{\rm (a)}] If $P,Q\in \reg(R)$ and $2P\cap 2Q\neq
\varnothing${\rm ,} then $\ell(Q)/2 \leq \ell(P) \leq 2 \ell(Q)$.

\item[{\rm (b)}] If $Q\in
\reg(R)$ and $x\in Q${\rm ,} $r\geq \ell(Q)${\rm ,} then $\sigma(B(x,r)\cap
4R) \leq C\theta_\sigma(R_1)r.$

\item[{\rm (c)}] For each $Q\in \reg(R)${\rm ,} there exist some $\vphi$-square  
$\wt{Q}$
which contains $Q$ such that $\delta_\sigma(Q,\wt{Q})\leq
C\theta_{\sigma}(R_1)$ and $\frac12\wt{Q}\cap\Gamma_R\neq
\varnothing$.
\end{itemize}
\end{lemma}

\Proof 
(a) Consider $P,Q \in \reg(R)$ such that $2P\cap 2Q\neq
\varnothing$. By construction, there exist  some $x\in P$ and some
$\vphi$-square $P_0\in \sss(R_1)$ or point $P_0\in G(R_1)$ (for
convenience, in this proof we identify points in $G(R_1)$ with
stopping squares in $\sss(R_1)$ with zero side length) such that
$\ell(P)\geq d_R(x)/20L$ and $$\dist(x,P_0) + \ell(P_0)\leq 1.1
d_R(x) \leq 22L\ell(P).$$ Thus, for any $y\in Q$,
\begin{eqnarray*}
\dist(y,P_0) + \ell(P_0) &\leq& \diam(2Q) + \diam(2P) +
\dist(x,P_0) + \ell(P_0) \\
& \leq & 3L\ell(Q) + 3L\ell(P) + 22L\ell(P),
\end{eqnarray*}
since $\diam(2Q) \leq L \,\diam(\vphi^{-1}(2Q)) =
L8^{1/2}\ell(Q)\leq 3L\ell(Q)$. So $d_R(y) \leq 3L\ell(Q) +
25L\ell(P)$ for all $y\in Q$. Therefore,
$$\ell(Q) \leq \frac1{10L}(3L\ell(Q) +
25L\ell(P)),$$ which yields $\ell(Q)\leq \frac{25}{7}\,\ell(P)
<4\ell(P)$. This implies $\ell(Q) \leq 2\ell(P)$, because $P$ and
$Q$ are $\vphi$-dyadic squares.

The inequality $\ell(P) \leq 2\ell(Q)$ is proved in an analogous
way.

\vspace{3mm}  (b) Take now $Q\in \reg(R)$ and $x\in Q$,
$r\geq \ell(Q)$. There exists some $y\in Q$ and some
$\vphi$-square $P_0\in \sss(R_1)$ such that
$d_R(y)/20L<\ell(Q)\leq d_R(y)/10$ and $$\dist(y,P_0) +
\ell(P_0)\leq 1.1 d_R(y) \leq 22L\ell(Q).$$ Thus $B(x,r)$ is
contained in some $\vphi$-square of the form
$\frac{Cr}{\ell(P_0)}P_0$, with $\frac{Cr}{\ell(P_0)}\geq1$ and
$C$ depending on $L$. Then, $$\sigma(B(x,r) \cap 4R) \leq
\sigma\bigl(\tfrac{Cr}{\ell(P_0)}P_0\cap 4R\bigr)
  \leq C\theta_\sigma(R_1)r.$$

\vspace{3mm} (c) We continue with the same notation as in
(b). Let $\wt{P}_0$ be a $\vphi$-square containing $P_0$ such that
$\delta_\sigma(P_0,\wt{P}_0)\leq C\theta_\sigma(R_1)$ and $\frac12
\wt{P}_0\cap \Gamma_R\neq \varnothing$ (given by (c) of Main Lemma
\ref{lemafi}). It is easily checked that there exists some
absolute constant $C_{26}\geq1$ such that $C_{26}\wt{P}_0$
contains $Q$. We set $\wt{Q}:= C_{26}\wt{P}_0$.
\Endproof\vskip4pt

Given $R \in \ttd(F)$, we denote by $\treeg(R)$ the tree of dyadic  
$\vphi$-squares whose
top $\vphi$-square is $R$ and whose terminal $\vphi$-squares are
the $\vphi$-squares $Q_i\in\reg(R)$ which are contained in $R$
(this is the same definition as the one for $\tree(R)$ in
Subsection \ref{subdecomp}, but with ${\rm Term}(R)$ replaced by
$\reg(R)$).

\begin{lemma} \label{treetreeg}
Given any $R \in \ttd(F)${\rm ,} if $Q\in\tree(R)$ and $\sigma(Q)>0${\rm ,} then  
$Q\in\treeg(R)$.
\end{lemma}

Roughly speaking, the lemma asserts that if we do not care about squares  
with vanishing $\sigma$-measure, then $\tree(R)\subset\treeg(R)$,
and so we always stop later in $\treeg(R)$ than in $\tree(R)$.

\Proof 
Let $R_1$ be the $4$-dyadic $\vphi$-square $R_1\in \ttt(F)$ such that  
$R\subset R_1\subset
7R$ is as in the definition of $\reg(R)$.

Let $Q_0\in\tree(R)$ be such that $\sigma(Q_0)>0$. To see that  
$Q_0\in\treeg(R)$, it is enough to show that $Q_0$ is not contained
in any square $Q_x$ like the ones appearing in \rf{eqqx}, with  
$x\in3R\cap\supp(\sigma)\setminus [G(R_1)\cup Z(\sigma)]$.
Suppose that this is not the case, so that $Q_0\subset Q_x$ for some  
$Q_x$ as above.
Since
$$R\cap \supp(\sigma)\subset \bigcup_{P\in\sss(R_1)} P \cup G(R_1) \cup  
Z(\sigma)$$
and $\sigma(Q_0) \neq 0$,
by the definition of $\tree(R)$, either there exists some square  
$P\in\sss(R_1)$
such that one of the $16$ dyadic squares which form $P$ (which is  
$4$-dyadic) is contained in $P$, or
there exists some $y_0 \in G(R_1)\cap Q_0$.
In any case, we deduce (identifying $y_0$ with a square $P$
with $\ell(P)=0$ in the latter case) that for any $y\in Q_x$
$$d_R(y) \leq \ell(P) + 2^{1/2} L \ell(Q_x) \leq (4 +  2^{1/2} L)\,  
\ell(Q_x).$$
In particular, this holds for $x=y$, and so
$$\ell(Q_x)\leq\frac1{10L}\,d_R(x) \leq \frac{4 +   
2^{1/2}}{10}\ell(Q_x)\leq \frac35 \ell(Q_x).$$
Thus $\ell(Q_x)=0$, which is a contradiction.
\hfill\qed

\Subsec{Construction of the approximating measure on
$\Gamma_R$} \label{subcons}
In this subsection we denote $\reg(R)=:\{Q_i\}_{i\geq1}$. For each
$i$, let $\wt{Q}_i$ be a $\vphi$-square containing $P_i$ such that
$\delta_\sigma(Q_i,\wt{Q}_i)\leq C\theta_\sigma(R_1)\simeq
\theta_\sigma(7R)$ and $\frac12 \wt{Q}_i\cap \Gamma_{R}\neq
\varnothing$. We may also suppose that $\diam(\Gamma_R)\geq
10\ell(R)$, since we can always extend $\Gamma_R$ if necessary.

\begin{lemma} \label{repart}
For each $i\geq1$ there exists some function $g_i\geq0$ supported
on $\Gamma_R\cap \wt{Q}_i$ such that
\begin{equation} \label{co1}
\int_{\Gamma_R} g_i\,d\HH^1 = \sigma(Q_i),
\end{equation}
\begin{equation} \label{co2}
\sum_i g_i \lesssim \theta_\sigma(R_1),
\end{equation}
and
\begin{equation} \label{co3}
\|g_i\|_\infty \ell(\wt{Q}_i) \lesssim\sigma(Q_i).
\end{equation}
\end{lemma}
\vskip8pt

\Proof 
The arguments are inspired by the Calder\'on-Zygmund decomposition
of \cite[Lemma 7.3]{Tolsa-bmo}.

We assume first that the family $\reg(R)=\{Q_i\}_i$ is finite. We
also suppose that $\ell(\wt{Q}_i)\leq\ell(\wt{Q}_{i+1})$ for all
$i$. The functions $g_i$ that we will construct will be of the
form $g_i=\alpha_i\chi_{A_i}$, with $\alpha_i\geq0$ and
$A_i\subset \wt{Q}_i$. We set
$\alpha_1:=\sigma(Q_1)/\HH^1(\wt{Q}_1\cap \Gamma_R)$ and
$A_1:=\wt{Q}_1\cap \Gamma_R$, so that $\int_{\Gamma_R} g_1d\HH^1 =
\sigma(Q_1).$ Notice by the way that $\|g_1\|_\infty \leq
C\sigma(Q_1)/\ell(\wt{Q}_1) \leq C\theta_\sigma(R_1)$.

To define $g_k$, $k\geq2$, we argue by induction. Suppose that
$g_1,\dots,g_{k-1}$ have been constructed, satisfy \rf{co1} and
$\sum_{i=1}^{k-1}g_i\leq B\theta_\sigma(R_1),$ where $B$ is
some\break\vskip-11pt\noindent  constant which will be chosen below. Let
$\wt{Q}_{s_1},\dots,\wt{Q}_{s_m}$ be the subfamily of
$\wt{Q}_1,\dots,\wt{Q}_{k-1}$ such that $\wt{Q}_{s_j}\cap
\wt{Q}_k\neq\varnothing$. Since $\ell(\wt{Q}_{s_j})\leq
\ell(\wt{Q}_k)$ (because of the nondecreasing sizes of the
$\wt{Q}_i$'s), we have $\wt{Q}_{s_j}\subset 3\wt{Q}_k$. Using
\rf{co1} for $i=s_j$, we get
\begin{eqnarray*}
\sum_j \int_{\Gamma_R} g_{s_j}\,d\HH^1 & \leq & \sum_j
\sigma(Q_{s_j}) \\ & \leq & \sigma(3\wt{Q}_k) \leq
C\theta_\sigma(R_1) \ell(\wt{Q}_k) \leq C_{27} \theta_\sigma(R_1)
\HH^1(\Gamma_R\cap \wt{Q}_k).
\end{eqnarray*}
Therefore, $$\HH^1\Bigl(\Gamma_R\cap \Bigl\{\sum_j g_{s_j} >
2C_{27}\theta_\sigma(R_1) \Bigr\}\Bigr) \leq \frac12
\HH^1(\Gamma_R\cap \wt{Q}_k).$$ So we set $$A_k := \Gamma_R
\cap\wt{Q}_k\cap \Bigl\{\sum_j g_{s_j} \leq
2C_{27}\theta_\sigma(R_1)\Bigr\},$$ and then $\HH^1(A_k) \geq
\HH^1(\Gamma_R\cap \wt{Q}_k)/2$. Also, we put $\alpha_k:=
\dfrac{\sigma(Q_k)}{\HH^1(A_k)},$ so that $\int_{\Gamma_R}
g_k\,d\HH^1 = \sigma(Q_k).$ Then,
\begin{equation} \label{www2}
\alpha_k \leq \frac{2\sigma(Q_k)}{\HH^1(\Gamma_R\cap \wt{Q}_k)}
\leq \frac{C\sigma(Q_k)}{\ell(\wt{Q}_k)} \leq
C_{28}\theta_\sigma(R_1).
\end{equation}
Thus, $$g_k + \sum_j g_{s_j} \leq (2C_{27} +
C_{28})\theta_\sigma(R_1).$$ We choose $B:= 2C_{27}+C_{28}$ and
\rf{co2} follows. Notice that \rf{co3} is proved in \rf{www2}.

Suppose now that $\{Q_i\}_i$ is not finite. For each fixed $N$ we
consider a family of squares $\{Q_i\}_{1\leq i\leq N}$. As above,
we construct functions $g_1^N,\dots,g_N^N$ with
$\supp(g_i^N)\subset \wt{Q}_i\cap\Gamma_R$ satisfying
$$\int_{\Gamma_R} g^N_id\HH^1 = \sigma(Q_i),\;\;\quad \sum_{i=1}^N
g_i^N \leq B\theta_\sigma(R_1),\quad \mbox{and}\quad
\|g_i^N\|_\infty \ell(\wt{Q}_i) \leq C\sigma(Q_i).$$ Then there
is a subsequence $\{g_1^k\}_{k\in I_1}$ which is convergent in
the weak $*$ topology of $L^\infty(\hhh)$ to some function $g_1\in
L^\infty(\hhh)$. Now we  take another subsequence $\{g_2^k\}_{k\in
I_2}$, $I_2\subset I_1$, convergent in the weak $*$ topology of
$L^\infty(\hhh)$ to another function $g_2\in L^\infty(\hhh)$, etc.
We have $\supp(g_i)\in \wt{Q}_i$. Further, \rf{co1}, \rf{co2} and
\rf{co3} also hold, because of the weak $*$ convergence.
\hfill\qed

\Subsec{A symmetrization lemma} \label{subsym}
Recall that by Lemma \ref{treetreeg}, $\tree(R)\subset \treeg(R)$. As a  
consequence,
$$\sum_{Q\in
\tree(R)}\int_Q (\wh{K}_{\sigma,J(Q)+1} \chi_F -
\wh{K}_{\sigma,J(Q)} \chi_F)\, d\sigma \leq \sum_{Q\in
\treeg(R)}\int_Q \cdots\,\,\, d\sigma.$$
Observe also that if $x\in Q_i$, then
\begin{multline} \label{we31}
\sum_{Q\in \treeg(R)} \chi_Q(x) \bigl(\wh{K}_{\sigma,J(Q)+1}
\chi_F(x) -
\wh{K}_{\sigma,J(Q)} \chi_F(x)\bigr)\\
\begin{split}
  & = \wh{K}_{\sigma,J(Q_i)+1} \chi_F(x) -
\wh{K}_{\sigma,J(R)+1} \chi_F(x)\\
& =  \iint_{\begin{subarray}{l}
\frac1{2L}\ell(Q_i) <|x-y|\leq \frac1{2L}\ell(R)\\
|x-z|,|y-z|\leq|x-y|
\end{subarray}}
c(x,y,z)^2d\sigma(y)d\sigma(z)\\
& \leq \iint_{\begin{subarray}{l}
y,z\in 2R\\
|x-y|>\frac1{2L}\ell(Q_i),\\
|x-z|,|y-z|\leq|x-y|\end{subarray}}c(x,y,z)^2d\sigma(y)d\sigma(z).
\end{split}
\end{multline}
In the last inequality we took into account that if $x\in R$ and
$|x-y|,|x-z|\leq \ell(R)/(2L)$, then $y,z\in 2R$.
  Analogously, if $x\in R\setminus \bigcup_i
Q_i$, we get
\begin{equation} \label{we32}
\sum_{Q\in \treeg(R)} \chi_Q(x) \bigl(\wh{K}_{\sigma,J(Q)+1}
\chi_F(x) - \wh{K}_{\sigma,J(Q)} \chi_F(x)\bigr) \leq
c^2_{\sigma|2R}(x).
\end{equation}

The lack of symmetry with respect to $x,y,z$ in the truncation of
the integrals that appear in \rf{we31} might cause some difficulties
in our estimates. This question is solved in the next lemma.

For any $y\in 3R$, we denote $\ell_y:=\ell(Q_i)$ if $y\in Q_i$,
and $\ell_y:=0$ if $y\in 3R\setminus \bigcup_i Q_i$.

\begin{lemma} \label{simpler}
\begin{itemize}
\item[{\rm (a)}] If $|x-y|\geq C_{29}^{-1} \ell_x${\rm ,} then $|x-y| \geq
C_{30}^{-1}\ell_y${\rm ,} with $C_{30}$ depending only on $C_{29}$ and $L$.
\item[{\rm (b)}] There exists a sufficiently small constant $\ve>0$ such that
\end{itemize}
\vglue-15pt
\begin{multline} \label{mult44}
\iiint_{\begin{subarray}{l}
x,y,z\in 2R\\
|x-y|>\frac1{2L}\ell_x
\\|x-z|,|y-z|\leq|x- 
y|\end{subarray}}c(x,y,z)^2d\sigma(x)d\sigma(y)d\sigma(z)
\\ \leq C\theta_\sigma(R_1)^2\sigma(R_1) +
\iiint_{\begin{subarray}{l}
x,y,z\in 2R\\
|x-y|\geq\ve(\ell_x+\ell_y)\\
|x-z|\geq\ve(\ell_x+\ell_z)\\
|y-z|\geq\ve(\ell_y+\ell_z)
\end{subarray}}c(x,y,z)^2d\sigma(x)d\sigma(y)d\sigma(z).
\end{multline}
\end{lemma}

\Proof 
First we show (a). Suppose  $\ell_x\neq0$, $\ell_y\neq0$.
Take $Q_i,Q_j\in\reg(R)$ such that $x\in Q_i$, and $y\in Q_j$. So
$\ell_x=\ell(Q_i)$ and $\ell_y=\ell(Q_j)$. If $|x-y|\leq
\ell(Q_j)/(2L)$, then $x\in 2Q_j$. Thus
$Q_i\cap2Q_j\neq\varnothing$, and then $\ell(Q_i)\geq\ell(Q_j)/2$,
which yields
$$|x-y|\geq C_{29}^{-1}\ell(Q_i)\geq \frac{C_{29}^{-1}}{2}\ell(Q_j).$$
So in any case we have
$$|x-y|\geq  
\min\Bigl(\frac1{2L},\,\frac{C_{29}^{-1}}{2}\Bigr)\ell(Q_j).$$

If $\ell_x=0$ or $\ell_y=0$ the arguments above also work, with
the convention $Q_i\equiv\{x\}$ or $Q_j\equiv\{y\}$.

Let us prove (b) now. We put
\begin{multline*}
\iiint_{\begin{subarray}{l}
x,y,z\in 2R\\
|x-y|>\frac1{2L}\ell_x
\\|x-z|,|y-z|\leq|x- 
y|\end{subarray}}c(x,y,z)^2d\sigma(x)d\sigma(y)d\sigma(z)
\\ =
\iiint_{\begin{subarray}{l}
x,y,z\in 2R\\
|x-y|>\frac1{2L}\ell_x\\
|x-y|\geq |x-z|\geq |y-z|
\end{subarray}}\!\!\!\!\!\!\!\!\cdots \;\;+
\iiint_{\begin{subarray}{l}
x,y,z\in 2R\\
|x-y|>\frac1{2L}\ell_x\\ 
|x-y|\geq |y-z|>|x-z|
\end{subarray}}\!\!\!\!\!\!\!\!\cdots \;\; =: A + B.
\end{multline*}
First we deal with the term $A$. By (a) we deduce that if
$|x-y|\geq\frac1{2L}\ell_x$, then $|x-y|\geq C^{-1}\ell_y$, and so
$|x-y|\geq\ve(\ell_x+\ell_y).$
If moreover $|x-y|\geq |x-z|\geq |y-z|$, then
$|x-z|\geq\frac12|x-y|\geq\frac1{4L}\ell_x$. Thus, $|x-z|\geq
C^{-1}\ell_z$ by (a), and so
$|x-z|\geq\ve(\ell_x+\ell_z).$
We obtain
\begin{eqnarray*}
A & \leq & \iiint_{\begin{subarray}{l}
x,y,z\in 2R\\
|x-y|\geq \ve(\ell_x+\ell_y)\\
|x-z|\geq \ve(\ell_x+\ell_z)
\end{subarray}}c(x,y,z)^2d\sigma(x)d\sigma(y)d\sigma(z)\\
& = & \iiint_{\begin{subarray}{l}
x,y,z\in 2R\\
|x-y|\geq \ve(\ell_x+\ell_y)\\
|x-z|\geq \ve(\ell_x+\ell_z)\\
|y-z|> \ell_y
\end{subarray}}\!\!\!\!\!\!\!\!\cdots \;\;
+ \iiint_{\begin{subarray}{l}
x,y,z\in 2R\\
|x-y|\geq \ve(\ell_x+\ell_y)\\
|x-z|\geq \ve(\ell_x+\ell_z)\\
|y-z|\leq \ell_y
\end{subarray}}\!\!\!\!\!\!\!\!\cdots \;\;
=:A_1 + A_2.
\end{eqnarray*}
To estimate $A_1$ we apply (a) again. Indeed, if $|y-z|> \ell_y$,
then $|y-z|\geq C^{-1}\ell_z$, and we get
$|y-z|\geq\ve(\ell_y+\ell_z).$ Therefore, $$A_1\leq
\iiint_{\begin{subarray}{l}
x,y,z\in 2R\\
|x-y|\geq \ve(\ell_x+\ell_y)\\
|x-z|\geq \ve(\ell_x+\ell_z)\\
|y-z|\geq \ve(\ell_y+\ell_z)
\end{subarray}}c(x,y,z)^2d\sigma(x)d\sigma(y)d\sigma(z).
$$
Now we deal with $A_2$. For each $y\in 2R$ we have
\begin{multline*}
\iint_{\begin{subarray}{l}
x,z\in 2R\\
|x-y|\geq \ve(\ell_x+\ell_y)\\
|x-z|\geq \ve(\ell_x+\ell_z)\\
|y-z|\leq \ell_y
\end{subarray}}
c(x,y,z)^2d\sigma(x)d\sigma(z) \leq  \iint_{\begin{subarray}{l}
x,z\in 2R\\
|x-y|\geq \ve\ell_y\\
|y-z|\leq \ell_y
\end{subarray}}
\frac{C}{|x-y|^2}d\sigma(x)d\sigma(z)\\
\leq  C\sigma(B(y,\ell_y))\int_{|x-y|\geq \ve\ell_y}
\frac1{|x-y|^2}d\sigma(x) \leq C\ve^{-1}\theta_{\sigma}(R_1)^2.
\end{multline*}
Therefore,
$$A_2 \leq C\ve^{-1}\theta_{\sigma}(R_1)^2\sigma(R_1).$$
Thus, $A$ is bounded above by the right-hand side of \rf{mult44}.

The term $B$ is estimated similarly to $A$. We will not go
through the details. Then we obtain
\begin{multline*}
\iiint_{\begin{subarray}{l}
x,y,z\in 2R\\
|x-y|>\frac1{2L}\ell_x
\\|x-z|,|y-z|\leq|x- 
y|\end{subarray}}c(x,y,z)^2d\sigma(x)d\sigma(y)d\sigma(z)
\\* \leq C\theta_\sigma(R_1)^2\sigma(R_1) + 2
\iiint_{\begin{subarray}{l}
x,y,z\in 2R\\
|x-y|\geq\ve(\ell_x+\ell_y)\\
|x-z|\geq\ve(\ell_x+\ell_z)\\
|y-z|\geq\ve(\ell_y+\ell_z)
\end{subarray}}c(x,y,z)^2d\sigma(x)d\sigma(y)d\sigma(z),
\end{multline*}
The reader may check that the `$2$' preceding the integral on the
right-hand side can be eliminated if one argues a little more
carefully (although this fact will be not needed for the estimates
below).
\Endproof\vskip4pt 

We denote
$$c^2_\ell(\sigma_{|2R}):= \iiint_{\begin{subarray}{l}
x,y,z\in 2R\\
|x-y|\geq\ve(\ell_x+\ell_y)\\
|x-z|\geq\ve(\ell_x+\ell_z)\\
|y-z|\geq\ve(\ell_y+\ell_z)
\end{subarray}}c(x,y,z)^2d\sigma(x)d\sigma(y)d\sigma(z).$$
We also set
$$c^2_{\ell,\sigma_{|2R}}(x) := \iint_{\begin{subarray}{l}
y,z\in 2R\\
|x-y|\geq\ve(\ell_x+\ell_y)\\
|x-z|\geq\ve(\ell_x+\ell_z)\\
|y-z|\geq\ve(\ell_y+\ell_z)
\end{subarray}}c(x,y,z)^2 d\sigma(y)d\sigma(z),$$
and
$$c^2_{\ell,\sigma}(A,B,C) := \iiint_{\begin{subarray}{l}
x\in A,\,y\in B,\,z\in C\\
|x-y|\geq\ve(\ell_x+\ell_y)\\
|x-z|\geq\ve(\ell_x+\ell_z)\\
|y-z|\geq\ve(\ell_y+\ell_z)
\end{subarray}}c(x,y,z)^2d\sigma(x)d\sigma(y)d\sigma(z).$$
Notice that $c^2_{\ell,\sigma}(A,B,C)$ is symmetric with respect to
$A,B,C$.

  By \rf{we31}, \rf{we32} and Lemma \ref{simpler}, to prove
\rf{eqclau} it is enough to show that
\begin{equation}\label{eqclau'}
c^2_\ell(\sigma_{|2R})\leq C\theta_\sigma(R_1)^2\sigma(R_1).
\end{equation}
We set $$G_R:=2R\setminus \bigcup_iQ_i.$$ Observe that
$\sigma$-almost all $G_R$ are  contained in $\Gamma_R$, by (a) of
Lemma \ref{lemafi}.

We split $c^2_\ell(\sigma_{|2R})$ as follows:
\begin{eqnarray} \label{we42}
c^2_\ell(\sigma_{|2R}) & = & c^2_{\ell,\sigma}\bigl(\cup_i Q_i,
\cup_i Q_i, \cup_i Q_i\bigr) + 3c^2_{\ell,\sigma}\bigl(\cup_i Q_i,
\cup_i Q_i, G_R\bigr) \\[5pt] && + 3c^2_{\ell,\sigma}\bigl(\cup_i
Q_i, G_R, G_R\bigr) +
c^2_{\sigma}\bigl(G_R,G_R,G_R\bigr).\nonumber
\end{eqnarray}

\Subsec{Estimate of $c^2_{\sigma}\bigl(G_R,G_R,G_R\bigr)$}
\label{subar3}
The measure $\sigma_{|G_R}$ coincides with $f\,d\hhh$, where $f$
is some function such that $\|f\|_\infty\leq
C\theta_\sigma(R_1)$. Since the Cauchy transform is bounded on
$L^2(\hhh)$ (with $\|\CC\|_{L^2(\hhh),L^2(\hhh)}$ bounded above
by some absolute constant), we have
\begin{eqnarray*}
c^2_{\sigma}\bigl(G_R,G_R,G_R\bigr) & \leq & \|f\|_\infty^3\,
c^2(\hhh) \,\lesssim \,\theta_\sigma(R_1)^3\HH^1(\Gamma_R)\\*
& \lesssim & \theta_\sigma(R_1)^3\diam(R_1) \,=\,
C\theta_\sigma(R_1)^2\sigma(R_1).
\end{eqnarray*}

\Subsec{Estimate of $c^2_{\ell,\sigma}\bigl(\cup_i Q_i,
\cup_i Q_i, \cup_i Q_i\bigr)$} \label{subq3}
We set \begin{equation} \label{dd1} c^2_{\ell,\sigma}\bigl(\cup_i
Q_i, \cup_j Q_j, \cup_k Q_k\bigr) = \sum_{i,j,k}
c^2_{\ell,\sigma}(Q_i,Q_j,Q_k).
\end{equation}
Now we put
\begin{eqnarray*}
Q_i & = & \Bigl[Q_i\cap (6\wt{Q}_j\cup6\wt{Q}_k)\Bigr]\cup
\Bigl[Q_i\setminus (6\wt{Q}_j\cup6\wt{Q}_k)\Bigr],\\
Q_j & = & \Bigl[Q_j\cap (6\wt{Q}_i\cup6\wt{Q}_k)\Bigr]\cup
\Bigl[Q_j\setminus(6\wt{Q}_i\cup6\wt{Q}_k)\Bigr],\\
Q_k & = & \Bigl[Q_k\cap (6\wt{Q}_i\cup6\wt{Q}_j)\Bigr]\cup
\Bigl[Q_k\setminus (6\wt{Q}_i\cup6\wt{Q}_j)\Bigr].
\end{eqnarray*}
We replace $Q_i$, $Q_j$, and $Q_k$ in \rf{dd1} by the right-hand
side of the identities above, and we get
\begin{multline*}
c^2_{\ell,\sigma}\bigl(\cup_i Q_i, \cup_j Q_j, \cup_k Q_k\bigr)\\
\begin{split}
\leq &\;\sum_{i,j,k} c^2_{\ell,\sigma}\bigl(Q_i\setminus
(6\wt{Q}_j\cup6\wt{Q}_k),\,\,Q_j\setminus(6\wt{Q}_i\cup6\wt{Q}_k),\,\,
Q_k\setminus (6\wt{Q}_i\cup6\wt{Q}_j)\bigr)\\
&\mbox{} + \sum_{i,j,k} c^2_{\ell,\sigma}\bigl(Q_i\cap
(6\wt{Q}_j\cup6\wt{Q}_k),\,\,Q_j,\,\,Q_k\bigr)\\
&\mbox{} + \sum_{i,j,k} c^2_{\ell,\sigma}\bigl(Q_i,\,\, Q_j\cap
(6\wt{Q}_i\cup6\wt{Q}_k),\,\,Q_k\bigr)\\
&\mbox{} + \sum_{i,j,k} c^2_{\ell,\sigma}\bigl(Q_i,\,\,Q_j,\,\,
Q_k\cap (6\wt{Q}_i\cup6\wt{Q}_j)\bigr)\\
=: & \;U + V_1 + V_2 + V_3.
\end{split}
\end{multline*}
First we estimate $V_1,\,V_2,\,V_3$. By symmetry, $V_1=V_2=V_3$,
and also
\begin{eqnarray*}
V_1 & \leq & \sum_{i,j,k} c^2_{\ell,\sigma}\bigl(Q_i\cap
6\wt{Q}_j,\,\,Q_j,\,\,Q_k\bigr) + \sum_{i,j,k}
c^2_{\ell,\sigma}\bigl(Q_i\cap 6\wt{Q}_k,\,\,Q_j,\,\,Q_k\bigr)\\
& = & 2\sum_{i,j,k} c^2_{\ell,\sigma}\bigl(Q_i\cap
6\wt{Q}_j,\,\,Q_j,\,\,Q_k\bigr) \\
& \leq & 2\sum_j c^2_{\ell,\sigma}(Q_j,\,6\wt{Q}_j,\,\,2R)\\
& = & 2\sum_j c^2_{\ell,\sigma}(Q_j,\,6\wt{Q}_j,\,6\wt{Q}_j)
+2\sum_j
c^2_{\ell,\sigma}(Q_j,\,6\wt{Q}_j,\,\,2R\setminus6\wt{Q}_j)\\
& =:& V_{1,1} + V_{1,2}.
\end{eqnarray*}
Now we deal with $V_{1,1}$:
\begin{multline*}
c^2_{\ell,\sigma}\bigl(Q_j, 6\wt{Q}_j,6\wt{Q}_j\bigr)  \leq
2\iiint_{\begin{subarray}{l}
x\in Q_j,\,y,z\in 6\wt{Q}_j\\
|x-z|\geq|x-y|\geq \ve\ell(Q_j)
\end{subarray}}c(x,y,z)^2d\sigma(x)d\sigma(y)d\sigma(z)\\
\begin{split}
& \lesssim \int_{x\in Q_j} \int_{\begin{subarray}{l} y\in 6\wt{Q}_j\\
|x-y|\geq \ve\ell(Q_j)
\end{subarray}}
\biggl(\int_{\begin{subarray}{l} z\in 6\wt{Q}_j\\
|x-z|\geq |x-y|
\end{subarray}}
\frac{1}{|x-z|^2}\,d\sigma(z)\biggr)
d\sigma(y)d\sigma(x)\\
& \lesssim \int_{x\in Q_j} \int_{\begin{subarray}{l} y\in 6\wt{Q}_j\\
|x-y|\geq \ve\ell(Q_j)
\end{subarray}}
\frac{\theta_\sigma(R_1)}{|x-y|} d\sigma(y)d\sigma(x)\\
& \lesssim  \int_{x\in Q_j}\theta_\sigma(R_1)^2d\sigma(x) =
\theta_\sigma(R_1)^2\sigma(Q_j).
\end{split}
\end{multline*}
Thus, $V_{1,1}\lesssim\theta_\sigma(R_1)^2\sigma(R_1).$ \smallbreak

The term $V_{1,2}$ is estimated likewise:
\begin{multline*}
c^2_{\ell,\sigma}\bigl(Q_j,6\wt{Q}_j, 2R\setminus6\wt{Q}_j\bigr)\\
\begin{split}
& \lesssim \int_{x\in Q_j} \int_{\begin{subarray}{l} y\in 6\wt{Q}_j\\
|x-y|\geq \ve\ell(Q_j)
\end{subarray}}
\biggl(\int_{\begin{subarray}{l} z\in 2R\\
|x-z|\geq C^{-1}\ell(\wt{Q}_j)
\end{subarray}}
\frac{1}{|x-z|^2}\,d\sigma(z)\biggr)
d\sigma(y)d\sigma(x)\\
& \lesssim \int_{x\in Q_j} \int_{\begin{subarray}{l} y\in 6\wt{Q}_j\\
|x-y|\geq \ve\ell(Q_j)
\end{subarray}}
\frac{\theta_\sigma(R_1)}{|x-y|} d\sigma(y)d\sigma(x)\\
& \lesssim \theta_\sigma(R_1)^2\sigma(Q_j),
\end{split}
\end{multline*}
and so $V_{1,2}\leq C\theta_\sigma(R_1)^2\sigma(R_1).$

It only remains to estimate $U$. Notice that if
\begin{equation}\label{dd4}
Q_i\setminus 6\wt{Q}_j\neq\varnothing\qquad\mbox{and}\qquad
Q_j\setminus 6\wt{Q}_i\neq\varnothing,
\end{equation}
then $2\wt{Q}_i\cap2\wt{Q}_j=\varnothing.$ Otherwise,
$2\wt{Q}_i\cap2\wt{Q}_j\neq\varnothing$ implies that either
$Q_i\subset2\wt{Q}_i\subset 6\wt{Q}_j$ or
$Q_j\subset2\wt{Q}_j\subset 6\wt{Q}_i$, which contradicts
\rf{dd4}. Thus,
\begin{eqnarray*}
U & = & \sum_{i,j,k} c^2_{\ell,\sigma}\bigl(Q_i\setminus
(6\wt{Q}_j\cup6\wt{Q}_k),\,\,Q_j\setminus(6\wt{Q}_i\cup6\wt{Q}_k),\,\,
Q_k\setminus (6\wt{Q}_i\cup6\wt{Q}_j)\bigr) \\
& \leq &
\sum_{\begin{subarray}{l}i,j,k:\,2\wt{Q}_i\cap2\wt{Q}_j=\varnothing,\\
2\wt{Q}_i\cap2\wt{Q}_k=\varnothing,\\2\wt{Q}_j\cap2\wt{Q}_k=\varnothing\end{subarray}}
c^2_\sigma(Q_i,Q_j,Q_k).
\end{eqnarray*}
Next we wish to compare $c^2_\sigma(Q_i,Q_j,Q_k)$ (for
$Q_i,\,Q_j,\,Q_k$ as in the last sum) with the curvature
$c^2_{\hhh}(g_i,g_j,g_k)$, where $g_i,g_j,g_k$ are the bounded
functions constructed in Lemma \ref{repart}, which are supported
on $\wt{Q}_i,\wt{Q}_j,\wt{Q}_k$ respectively. We set
$$c^2_{\hhh}(g_i,g_j,g_k) : = \iiint_{\Gamma_R^3}
c(x,y,z)^2\,g_i(x)\,g_j(y)\,
g_k(z)\,d\HH^1(x)d\HH^1(y)d\HH^1(z).$$ If $x,x'\in \wt{Q}_i$,
$\,y,y'\in \wt{Q}_j$ and $z,z'\in \wt{Q}_k$, then
\begin{eqnarray*}
c(x,y,z)^2 &\leq &2c(x',y',z')^2 \\
&&\mbox{} + \frac{C\ell(\wt{Q}_i)^2}{|x-y|^2|x-z|^2} +
\frac{C\ell(\wt{Q}_j)^2}{|y-x|^2|y-z|^2} +
\frac{C\ell(\wt{Q}_k)^2}{|z-x|^2|z-y|^2},
\end{eqnarray*}
by Lemma \ref{curpert}. If we integrate $x\in Q_i$, $y\in Q_j$,
and $z\in Q_k$ with respect to $\sigma$, and also $x'\in \wt{Q}_j$
with respect to the measure $g_i\, d\hhh$, $y'\in \wt{Q}_j$ with
respect to $g_j\,d\hhh$, and $z'\in \wt{Q}_k$ with respect to
$g_k\,d\hhh$, then we get
\begin{equation*}
\begin{split}
c_\sigma^2(Q_i,Q_j,Q_k) \;\leq \,\;& 2 c^2_{\hhh}(g_i,g_j,g_k)\\
&\mbox{} + \iiint_{\begin{subarray}{l} x\in Q_i\\y\in Q_j\\z\in
Q_k
\end{subarray}}
\frac{C\ell(\wt{Q}_i)^2}{|x-y|^2|x-z|^2}\,d\sigma(x)
d\sigma(y)\,d\sigma(z)\\
&\mbox{} + \iiint_{\begin{subarray}{l} x\in Q_i\\y\in Q_j\\z\in
Q_k
\end{subarray}}
\frac{C\ell(\wt{Q}_j)^2}{|y-x|^2|y-z|^2}
\,d\sigma(x)\,d\sigma(y)\,d\sigma(z) \\
& \mbox{} + \iiint_{\begin{subarray}{l} x\in Q_i\\y\in Q_j\\z\in
Q_k
\end{subarray}}
  \frac{C\ell(\wt{Q}_k)^2}{|z-x|^2|z-y|^2}
\,d\sigma(x)\,d\sigma(y)\,d\sigma(z).
\end{split}
\end{equation*}
Therefore, by symmetry,
\begin{eqnarray} \label{mee1}
U & \leq  & 2\sum_{i,j,k} c^2_{\hhh}(g_i,g_j,g_k) \\ &&+
3 \sum_i \iiint_{\begin{subarray}{l} x\in
Q_i\\|x-y|>C^{-1}\ell(\wt{Q}_i)
\\|x-z|> C^{-1}\ell(\wt{Q}_i)
\end{subarray}}
\frac{C\ell(\wt{Q}_i)^2}{|y-x|^2|y-z|^2}
\,d\sigma(x)\,d\sigma(y)\,d\sigma(z).
\nn
\end{eqnarray}
By Lemma \ref{repart}, $g:=\sum_i g_i \lesssim
\theta_\sigma(R_1)$, and then
\begin{eqnarray} \label{mee2}
  \sum_{i,j,k} c^2_{\hhh}(g_i,g_j,g_k) &= &
c^2(g\,d\hhh) \lesssim \theta_\sigma(R_1)^3 c^2(\hhh) \\
&\lesssim & \theta_\sigma(R_1)^3\HH^1(\Gamma_R) \lesssim
\theta_\sigma(R_1)^2\sigma(R_1).
\nn
\end{eqnarray}
The integral in \rf{mee1} is estimated as follows:
\begin{multline*}
\iiint_{\begin{subarray}{l} x\in Q_i\\|x-y|>C^{-1}\ell(\wt{Q}_i)
\\|x-z|> C^{-1}\ell(\wt{Q}_i)
\end{subarray}}
\frac{\ell(\wt{Q}_i)^2}{|y-x|^2|y-z|^2}
\,d\sigma(x)\,d\sigma(y)\,d\sigma(z) \\
\begin{split} \;= \,&
\ell(\wt{Q}_i)^2 \int_{x\in Q_i}
\left(\int_{|x-y|>C^{-1}\ell(\wt{Q}_i)}
\frac1{|y-x|^2}\,d\sigma(y)\right)^2d\sigma(x)\\
\;\lesssim \;& \ell(\wt{Q}_i)^2 \int_{x\in Q_i}
\frac{\theta_\sigma(R_1)^2}{\ell(\wt{Q}_i)^2} \,d\sigma(x) =
\theta_\sigma(R_1)^2\sigma(Q_i).
\end{split}
\end{multline*}
 From \rf{mee1}, \rf{mee2} and the preceding estimate we get
$$U \lesssim \theta_\sigma(R_1)^2\sigma(R_1).$$
We are done.
 
\Subsec{Estimates of $c^2_{\ell,\sigma}\bigl(\cup_i Q_i,
\cup_i Q_i, G_R\bigr)$ and $c^2_{\ell,\sigma}\bigl(\cup_i Q_i,
G_R, G_R\bigr)$} \label{subalt}
We leave these estimates for the reader. The arguments are similar
to the ones used for $c^2_{\ell,\sigma}\bigl(\cup_i Q_i, \cup_i
Q_i, \cup_i Q_i\bigr)$. In fact, notice that if by convention one
allows the squares $Q_i$ to be points, then $\bigl(\cup_i
Q_i\bigr) \times \bigl(\cup_i Q_i\bigr) \times G_R$ and
$\bigl(\cup_i Q_i\bigr)\times G_R \times G_R$ are subsets of
$\bigl(\cup_i Q_i\bigr)^3$.

\demo{Acknowledgements} I would like to thank Joan Verdera for
pointing out to me the necessity
of the bilipschitz condition on $\vphi$ for  
$\gamma(E)\simeq\gamma(\vphi(E))$ to hold for any compact $E\subset\C$,  
as
stated in Proposition \ref{propconv}.
Also, I would like to thank the referee for his suggestions, which  
helped to improve
the clarity of the paper.

\end{document}